%% file: ms.tex
\documentclass[a4paper,10pt,3p]{elsarticle}

\usepackage{amsmath}
\usepackage{amsthm}
\usepackage{amssymb}
\usepackage{graphicx}
\usepackage{epstopdf} 
\usepackage{bm}
\usepackage[disable]{todonotes} 
\usepackage{booktabs}
\usepackage{caption}
\usepackage{subcaption}
\usepackage{framed}
\usepackage{hyperref}
\usepackage{placeins}
\graphicspath{{figures/}}
\usepackage{listings}
\usepackage{nameref}
\usepackage{etoolbox}

\theoremstyle{remark}

\newcommand{\hl}[1]{\textcolor{black}{#1}}

\newcommand{\rd}{\textrm{d}}
\newcommand{\vt}[1]{\bm{#1}}

\newcommand{\dd}[2]{\frac{\partial #1}{\partial #2}}

\DeclareMathOperator*{\argmin}{arg\,min}

\makeatletter
\def\ps@pprintTitle{%
 \let\@oddhead\@empty
 \let\@evenhead\@empty
 \def\@oddfoot{}%
 \let\@evenfoot\@oddfoot}
\makeatother

\newcommand{\R}[1]{\label{#1}}

\begin{document}

\title{Non-linearly stable reduced-order models for incompressible flow \\with energy-conserving finite volume methods}
\author[1]{B. Sanderse}
\ead{B.Sanderse@cwi.nl}
\address[1]{Centrum Wiskunde \& Informatica, Amsterdam, the Netherlands}

\begin{abstract}
A novel reduced-order model (ROM) formulation for incompressible flows is presented with the key property that it exhibits non-linearly stability, independent of the mesh (of the full order model), the time step, the viscosity, and the number of modes. The two essential elements to non-linear stability are: (1) \textit{first discretise} the full order model, and \textit{then project} the discretised equations, and (2) use spatial and temporal discretisation schemes \hl{for the full order model} that are \textit{globally energy-conserving} (in the limit of vanishing viscosity). For this purpose, as full order model a staggered-grid finite volume method in conjunction with an implicit Runge-Kutta method is employed. In addition, a constrained singular value decomposition is \hl{employed} which enforces global momentum conservation. The resulting \hl{`velocity-only'} ROM is thus globally conserving mass, momentum and kinetic energy. For non-homogeneous boundary conditions, a (one-time) Poisson equation is solved that accounts for the boundary contribution. The stability of the proposed ROM is demonstrated in several test cases. \hl{Furthermore, it is shown that explicit Runge-Kutta methods can be used as a practical alternative to implicit time integration at a slight loss in energy conservation.}
\end{abstract}

\begin{keyword}
incompressible Navier-Stokes equations \sep reduced-order model \sep energy conservation \sep POD-Galerkin \sep finite volume method \sep stability
\end{keyword}

\maketitle


\section{Introduction}\label{sec:introduction}
\todo{Background/use of ROMs} The simulation of turbulent fluid flows is an ongoing challenge in the scientific community. The computational cost of Direct Numerical Simulation (DNS) or Large Eddy Simulation (LES) of turbulent flows quickly becomes imperative when one is interested in \textit{control}, \textit{design}, \textit{optimization} and \textit{uncertainty quantification} \cite{Benner2015,Smith2014}. For these purposes, a reduction in complexity of the full model is required to arrive at a computationally tractable model, a so-called reduced-order model (ROM). Several techniques exist to construct a ROM, such as balanced truncation, Krylov subspace methods, and POD-Galerkin methods \cite{Antoulas2005,Taira2017}. In this work we focus on one of the most popular techniques, the POD-Galerkin method, in which the governing equations of the full model  are projected onto a lower-dimensional space via a Galerkin step, with the projection basis determined from a proper orthogonal decomposition (POD) of snapshots of the full order model (FOM).


\todo{ROM issues related to fluid flow} Projection-based models have been shown to work for a large class of problems, such as diffusion-dominated linear time-invariant (LTI) systems, in which the input-output relation of the full model can be represented by a lower-dimensional model, due to rapid decay of the singular values of the Hankel matrix \cite{Benner2015}. However, in turbulent flow, which is a nonlinear, convection-dominated problem, the construction of accurate and stable ROMs is still an open challenge. There are several (related) reasons why current ROMs have issues with accuracy and stability in case of turbulent flows: the Kolmogorov $N$-width decays too slowly; the non-linear dynamical system is very sensitive to perturbations; the modes with low energy (small scales, dissipation) which are neglected in the POD procedure are relevant for the dynamics of the large scales; the reduced model can have different stability characteristics \cite{Aubry1993,Carlberg2011,Fick2018,Noack2016,Rempfer2000,Stabile2018}. 

\todo{Possible solutons} A number of approaches have been proposed to tackle the issues of accuracy and stability; we summarize the list in \cite{Fick2018}: including dissipation via a closure model (see also \cite{Cazemier1998}); modifying the POD basis by including functions that resolve a range of scales; using a minimum residual formulation \cite{Carlberg2018}; using an inner product different from $L_2$, e.g.\ based on $H_{1}$. \hl{In this work we are interested in tackling the issues of accuracy and stability through} the promising approach of structure-preserving model reduction, in which reduced-order models are developed in such a way that invariants and/or symmetries of the full model are kept \cite{Afkham2017a,Afkham2018c,Carlberg2018,Karasozen2018,Peng2016}. 
An example is a ROM that inherits the symplectic form of a Hamiltonian system, leading to a ROM that is applicable to stable long-time integration \cite{Peng2016}. 

\todo{Structure-preserving for INS} In the incompressible Navier-Stokes equations, \hl{conservation (or dissipation) of kinetic energy is typically used in order to proof non-linear stability, as the kinetic energy is a norm of the solution}. Although the incompressible Navier-Stokes equations do not form a Hamiltonian system, several symmetries are present in the equations which are tightly related to conservation of kinetic energy, such as the skew-symmetry of the convective operator, \hl{and the symmetry of the diffusive operator}. Several adaptations to the classic POD-Galerkin method were developed to take into account symmetry or invariance properties of the Navier-Stokes equations. For example, Balajewicz et al.\ \cite{Balajewicz2013} added a power balance equation for the resolved turbulent kinetic energy when solving for the POD basis functions and coefficients. Mohebujjaman et al.\ \cite{Mohebujjaman2019} employed a combined projection and data-driven approach in a finite-element context and obtain correction terms by solving a constrained minimization problem, using as Ansatz the negative definiteness of the diffusion operator and the energy-conserving property of the convection term. Mohebujjaman et al.\ also \cite{Mohebujjaman2017} investigated conservation of mass and energy of the ROM in the context of a finite element framework and discuss the treatment of non-homogeneous boundary conditions via a Stokes extension in order to mimic the continuous energy balance. Carlberg et al.\ \cite{Carlberg2018} considered conservative model reduction in a finite-volume context by solving a constrained optimization problem at each time step. Rowley et al.\ \cite{Rowley2004} considered the choice of an appropriate inner product and corresponding energy norm for compressible flow. Kalashnikova et al.\ \cite{Kalashnikova2014} considered energy stability in terms of a continuous formulation. \hl{However, to the author's knowledge, a reduced order model for the incompressible Navier-Stokes equations (based on a finite-volume full-order method) that combines global mass, momentum and kinetic energy conservation (in the invisicid limit) has not been developed to date. This will be the subject of this paper.}

\todo{Stability; divergence-free basis; no pressure}

\R{Lit1}\hl{To achieve energy conservation, we require the discretised full-order model to respect several symmetries of the Navier-Stokes equations: skew-symmetry of the convective term, symmetry of the diffusive term, and the compatibility relation between the divergence (of the velocity) and the gradient (of the pressure). With the latter property satisfied, the contribution of the pressure in the kinetic energy equation disappears (in closed or periodic domains) due to the fact that the velocity-field is divergence free \cite{Holmes2012}. Similarly, when a divergence-free basis for the velocity-field is taken in a Galerkin method, the pressure gradient term disappears from the momentum equations. The resulting model is `velocity-only', which has the advantage that it is faster to evaluate in the online phase (compared to a velocity-pressure formulation), and moreover can avoid the inf-sup stability issues associated to a mixed formulation \cite{Fonn2019}, so that stabilization techniques, such as supremizer modes or a pressure Poisson equation \cite{Ballarin2013,Caiazzo2014,Rubino2019,Stabile2018} are not needed. An important issue when deriving such a `velocity-only' formulation is, next to the discrete compatibility between divergence and gradient operators, the requirement of a discretely divergence-free velocity field. The importance of this property for the accuracy of the ROM was also confirmed in \cite{Mohebujjaman2017}. Another issue is that the velocity-only formulation requires a method to recover the pressure from the velocity field \cite{Akhtar2009,Kean2019} and requires adaptation if the pressure is non-negligible on outflow boundaries \cite{Noack2005}. In this work, these issues are effectively handled by choosing as discrete FOM a finite-volume discretisation on a staggered grid. This discretisation is energy-conserving in space \cite{Harlow1965,Verstappen2003} and time \cite{Sanderse2013b}, possesses the convective skew-symmetry, the diffusive symmetry, the compatibility relation between divergence and gradient, and has divergence-free velocity fields. \R{PBC3}Furthermore, it does not require boundary conditions for the pressure on parts of the boundary where the velocity is prescribed, allowing straightforward recovery of the pressure on the ROM level.}\R{Lit2}

\todo{Our approach} Based on this energy-conserving full-order model, we propose a new stable energy-conserving reduced-order model, which possesses non-linear stability, independent of the viscosity, mesh size, time step or number of modes. A constrained SVD method is employed with a weighted inner product, \textit{in such a way that the reduced model is mass-, momentum- and energy-conserving}. \hl{The resulting ROM is `velocity-only' and therefore does not suffer from inf-sup stability issues. As we will show, the pressure can be computed, if desired, as a postprocessing step by solving a Poisson equation at the reduced-order level.} An important \R{DiscProj}\hl{concept in our work is that we first discretise the equations, and then perform the projection step. Although this approach is common in the dynamical systems community \cite{Benner2015}, and has also been used by e.g.\ \cite{Carlberg2018,Stabile2018}, we are explicitly mentioning this choice because it allows us to prove stability on the level of the ROM and because it simplifies the boundary condition treatment. In particular, once the discrete full-order model and its boundary conditions are specified, no additional boundary conditions will be needed on the ROM level. Non-homogeneous boundary conditions are handled in such a way that the ROM stays velocity-only, and such that no additional boundary conditions need to be imposed when recovering the pressure via a pressure Poisson equation.} 

\begin{figure}[hbtp]
\centering
\includegraphics[width=0.55 \textwidth]{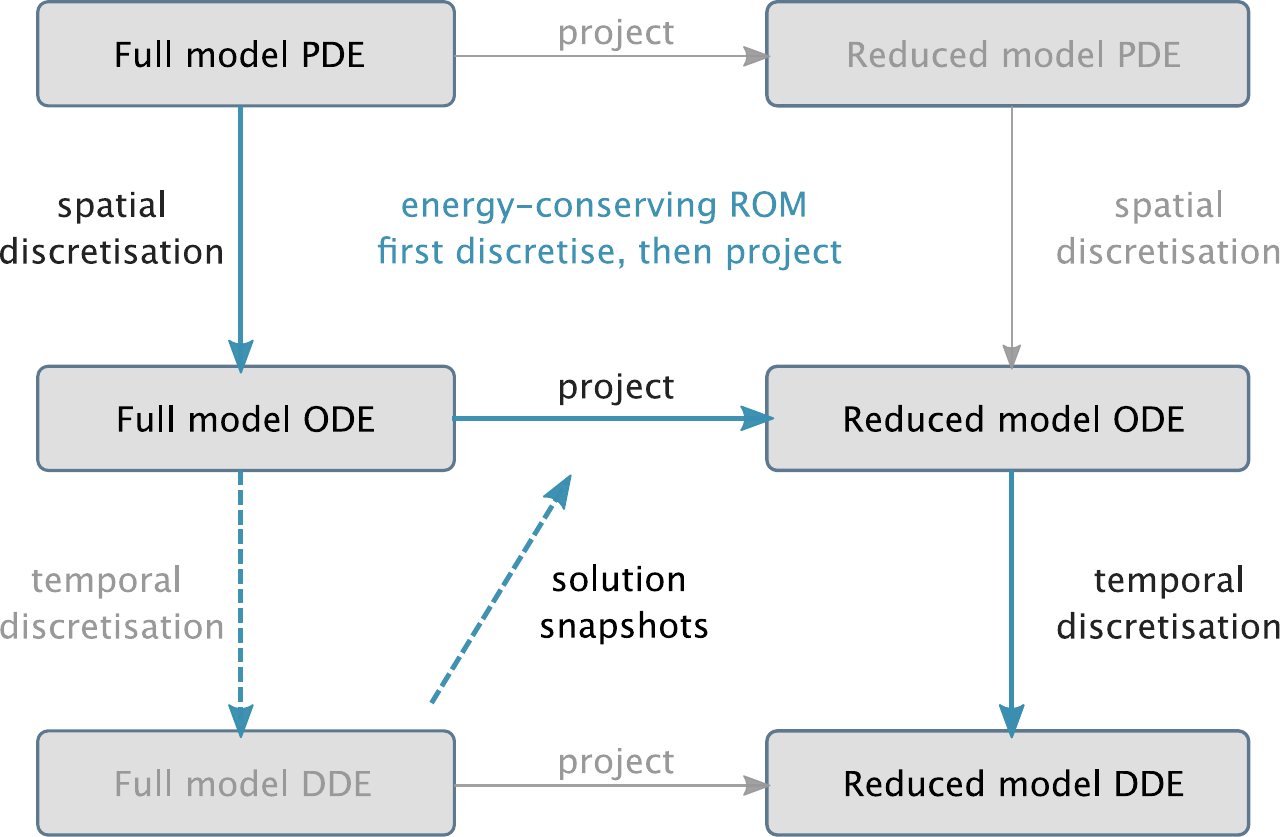}
\caption{Our approach to energy-stable reduced order models follows the blue arrows: first spatial discretisation, then projection (ODE = ordinary differential equation, DDE = discretised differential equation). \label{fig:discretise_projection}}
\end{figure}

A graphical summary of the approach is shown in figure \ref{fig:discretise_projection}: \R{DiscProj3}\hl{first spatial discretization, then projection, and then temporal discretization. As said, the order of the first two steps is especially important, because spatial discretization and projection do generally not commute (amongst others due to the presence of boundary conditions). On the other hand, the last two steps generally commute: for many time integration methods the same ROM would be obtained when first discretizing in time and then projecting the fully discrete equations \cite{Carlberg2017a}. However, this latter sequence does not give us the flexibility to take a different time integration method for the ROM than for the FOM, e.g.\ changing from explicit to implicit time integration. This flexibility will be used in the test cases to make our approach computationally efficient.}

We limit ourselves in the analysis in several important aspects. First, we will consider the so-called \textit{solution reproduction problem}, which is the first step before solving the full \textit{parametric problem} \cite{Fick2018}. Second, we will focus mainly on the \textit{non-linear stability} of the discrete ROM, and refrain from turbulent flows and closure models. Given a non-linearly stable \hl{and non-dissipative} ROM, we have a framework in which we can in future work assess the accuracy of closure models, e.g.\ \cite{Cazemier1998,Mohebujjaman2019,Wang2012a}. Our definition of stability should therefore be interpreted in the classical sense: a certain norm of the solution stays bounded in time. Note that this is different from the `stabilizing' methods that have been proposed in the ROM community, which are required to counteract numerical oscillations \cite{Weller2009}. This latter type of (in)stability is not the focus of this paper.

\todo{Novelty} The novelty of this paper is threefold. First, we derive an energy-conserving reduced-order model, which possesses \textit{nonlinear stability independent of the mesh and time step of the FOM, and independent of the time step and number of modes of the reduced-order model}. Second, we \textit{enforce global momentum conservation via a constrained singular value decomposition}. Third, we propose a procedure to handle non-homogeneous boundary conditions \hl{which keeps the ROM `velocity-only' and the pressure recovery stage straightforward}.

\todo{Paper outline} This paper is organized as follows. First, in section \ref{sec:fullordermodel} we discuss symmetry and energy-conservation properties of the full-order model (the incompressible Navier-Stokes equations) at the continuous, semi-discrete and fully discrete level. In section \ref{sec:ECPOD} we construct the new POD-Galerkin method, which conserves mass, momentum and energy globally. In section \ref{sec:BC} non-homogeneous boundary conditions are discussed. \hl{In section \ref{sec:implementation} we discuss implementation issues, including precomputing the reduced operators, the pressure recovery, and alternative (explicit) time integration strategies}. In section \ref{sec:results} the theoretical results are demonstrated for three cases: a shear layer roll-up, a lid-driven cavity flow, and an actuator disk with unsteady wake. 

\input{fullordermodel.tex}

\input{ECPOD.tex}

\input{BC.tex}

\input{implementation.tex}

\input{results.tex}

\input{conclusions.tex}

\section*{Acknowledgements}
The author would like to thank Michiel Hochstenbach (Eindhoven University of Technology) for the stimulating discussions on the constrained SVD decomposition, Giovanni Stabile (SISSA mathLab) and Kelbij Star (UGent, SCK-CEN) for many stimulating discussions on ROMs, and the anonymous reviewers for their constructive remarks.

\appendix

\input{convective_operator.tex}

\input{BC_FOM.tex}

\input{adapted_norm.tex}

\section*{References}
\bibliographystyle{abbrv}
\bibliography{refs}

\end{document}

%% file: fullordermodel.tex
\section{Energy-conserving discretisation of the incompressible Navier-Stokes equations}\label{sec:fullordermodel}
In order to develop the ROM, the energy (in)equalities of the FOM on the continuous and discrete level are needed. These are derived in this section.

\subsection{Continuous energy estimate}
The incompressible Navier-Stokes equations describe conservation of mass and momentum:
\begin{align}
\nabla \cdot \vt{u} &= 0, \label{eqn:divfree_continuous}\\
\dd{\vt{u}}{t} + \nabla \cdot (\vt{u} \otimes \vt{u}) &= -\nabla p + \nu \nabla \cdot (\nabla \vt{u} + (\nabla \vt{u})^T) + \vt{f},\label{eqn:mom_continuous}
\end{align}
where $\vt{u}(\vt{x},t)$ is the velocity field, $p(\vt{x},t)$ is the \R{ModPres}\hl{modified pressure (pressure scaled by density)}, $\vt{x} \in \Omega \subset \mathbb{R}^d$ ($d=2$ or $3$), $t$ denotes time, $\nu\geq 0$ the kinematic viscosity, \hl{and $\vt{f}(\vt{x},t)$ represents body forces}. The equations are supplemented with an initial condition 
\begin{equation}
\vt{u}(\vt{x},0) = \vt{u}_{0}(\vt{x}),
\end{equation}
and boundary conditions, e.g.\ no-slip conditions
\begin{equation}\label{eqn:bc}
\vt{u} = \vt{0}  \quad \text{on} \quad \partial \Omega.
\end{equation}
In this section and the next (section \ref{sec:ECPOD}), the focus will be on periodic and no-slip conditions. \hl{The extension to more generic boundary conditions (inflow, outflow) will be addressed in section \ref{sec:BC}.} We introduce the convection and diffusion operators $C(\vt{u},\vt{u}):=\nabla \cdot (\vt{u} \otimes \vt{u})$ and $D \vt{u} := \nabla \cdot (\nabla \vt{u} + (\nabla \vt{u})^T)$. Other forms of the convective operator are detailed in \ref{sec:convective_operator}.

To derive the kinetic energy equation, an inner product is needed. We choose the $L_2(\Omega)$ inner product and induced norm  \cite{Berselli2006}:
\begin{equation}\label{eqn:innerproduct_continuous}
(\vt{u},\vt{v}) := \int_{\Omega} \vt{u} \cdot \vt{v} \ \rd \Omega, \qquad \| \vt{u} \| := (\vt{u},\vt{u})^{1/2}.
\end{equation}
The kinetic energy is then defined as $K:=\frac{1}{2} \| \vt{u} \|^2$. An equation for the evolution of $K$ is derived by differentiating $K$ in time and substituting the momentum equation:
\begin{align}
2 \frac{\rd K}{\rd t} &= \frac{\rd (\vt{u},\vt{u})}{\rd t} =-(C(\vt{u},\vt{u}),\vt{u}) - (\vt{u},C(\vt{u},\vt{u})) - (\nabla p, \vt{u}) - (\vt{u}, \nabla p) + (D \vt{u},\vt{u}) + (\vt{u}, D \vt{u}) + (\vt{u},\vt{f}) + (\vt{f},\vt{u}).
\end{align}
The equation simplifies due to three symmetry properties. These symmetry properties will be crucial in developing an energy-stable ROM. First, due to the skew-symmetry of $C(\vt{u},\vt{u})$, we have $(C(\vt{u},\vt{u}),\vt{u})=0$ for periodic or no-slip boundary conditions (see also \ref{sec:convective_operator}). Second, the pressure gradient contribution disappears because $(\nabla p, \vt{u})=(p,\nabla \cdot \vt{u})=0$. Third, due to the symmetry of the diffusive operator we can write \hl{$(D \vt{u},\vt{u})=-(\nabla \vt{u},\nabla \vt{u})$}. The kinetic energy balance then reduces to 
\begin{equation}\label{eqn:energy_continuous}
\frac{\rd K}{\rd t} = -\nu \| \nabla \vt{u} \|^2, 
\end{equation}
in the absence of body forces and boundary conditions. Consequently, in viscous flow the kinetic energy of the flow can only decrease in time through dissipation, and in inviscid flow it is conserved.


\subsection{Spatial discretisation and semi-discrete energy equation}\label{sec:spatial}
In order to construct a non-linearly stable ROM, we first require that the spatial discretisation mimics the energy-conserving properties of the continuous equations just derived. To this end, we consider a finite volume discretisation on a staggered cartesian grid, \hl{also known as the Marker-and-Cell (MAC) method} \cite{Harlow1965,Sanderse2013,Verstappen2003} -- see figure \ref{fig:staggered}. \hl{This specific choice is an important one, as will be explained soon}. For simplicity, we restrict ourselves to a second-order method in two dimensions and partition the domain in $N_{p} = N_{x} \times N_{y}$ finite volumes. We introduce the (time-dependent) solution vectors $u_{h}(t) \in \mathbb{R}^{N_{p}}$, $v_{h}(t) \in \mathbb{R}^{N_{p}}$ and $p_{h}(t) \in \mathbb{R}^{N_{p}}$, which consist of the (time-dependent) unknowns $u_{i+1/2,j}(t)$, $v_{i,j+1/2}(t)$, and $p_{i,j}(t)$, respectively (for $i=1\ldots N_{x}$, $j=1 \ldots N_{y}$). The explicit time-dependence will be suppressed when no confusion can arise. The horizontal and vertical velocity components are gathered in the vector $V_{h} = \begin{pmatrix} u_{h} \\ v_{h} \end{pmatrix} \in \mathbb{R}^{N_{V}}$, with $N_{V}=2N_{p}$.

We integrate the divergence-free constraint \eqref{eqn:divfree_continuous} over a finite volume centred around the unknown $p_{i,j}$, which yields
\begin{equation}
 \bar{u}_{i+1/2,j} - \bar{u}_{i-1/2,j} + \bar{v}_{i,j+1/2} - \bar{v}_{i,j-1/2} = 0.
\end{equation}
The notation $\bar{(.)}$ indicates integration over a face of the finite volume, approximated e.g.\ by $\bar{u}_{i+1/2,j}= u_{i+1/2,j} \Delta y$. In matrix-vector notation, the above equation can be written for all pressure volumes as
\begin{equation}
M_{h} V_{h} = 0, \qquad M_{h} \in \mathbb{R}^{N_{p} \times N_{V}}.
\end{equation}

\begingroup
\begin{figure}[hbtp]
\fontfamily{lmss} 
\fontsize{11pt}{11pt}\selectfont
\centering 
\def\svgwidth{0.4 \textwidth}
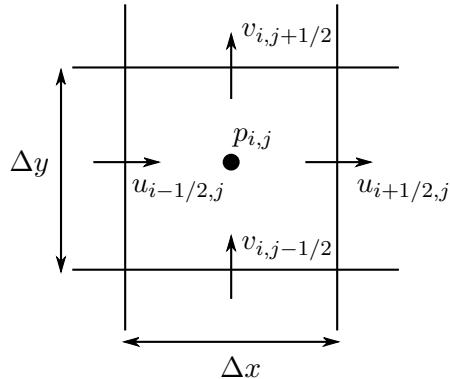 
\caption{Staggered grid with positioning of unknowns around a pressure volume.}
\label{fig:staggered}
\end{figure}
\endgroup

Next we integrate the horizontal component of the momentum equation over a finite volume centred around the unknown $u_{i+1/2,j}$. The convective term in divergence form is discretised by (mesh-independent) interpolation of the neighbouring fluxes. The divergence form ensures that momentum is conserved. When ensuring that the velocity field is discretely divergence free, the divergence form can be rewritten into the following equivalent skew-symmetric form (see \ref{sec:convective_operator_discrete}):
\begin{multline}\label{eqn:convection_rewritten}
[C^{u}_{h} (V_{h},u_{h})]_{i+1/2,j} := \frac{1}{2} u_{i+3/2,j} \frac{1}{2} \left(\bar{u}_{i+1/2,j} + \bar{u}_{i+3/2,j} \right) - \frac{1}{2} u_{i-1/2,j} \frac{1}{2} \left(\bar{u}_{i-1/2,j} + \bar{u}_{i+1/2,j} \right) \\
 + \frac{1}{2} u_{i+1/2,j+1} \frac{1}{2} \left(\bar{v}_{i,j+1/2} + \bar{v}_{i+1,j+1/2} \right) - \frac{1}{2} u_{i+1/2,j-1} \frac{1}{2} \left(\bar{v}_{i,j-1/2} + \bar{v}_{i+1,j-1/2} \right).
\end{multline} 
The skew-symmetric form ensures that kinetic energy is conserved. As a consequence, the discretised convection operator is conserving both momentum and energy. A similar scheme can be derived for the vertical component. The full convection operator then reads
\begin{equation}
C_{h} (V_{h},V_{h}) = \begin{pmatrix}
C^{u}_{h} (V_{h},u_{h}) \\
C^{v}_{h} (V_{h},v_{h})
\end{pmatrix} =: \tilde{C}_{h} (V_{h}) V_{h}.
\end{equation}
The notation $\tilde{C}_{h}(V_h)V_{h}$ is useful for making a distinction between the convect\textit{ing} quantity (the quantity between brackets) and the convect\textit{ed} quantity. The definition of $\tilde{C}_h$ is possible due to the fact that the nonlinearity of the convective term is only quadratic. The skew-symmetry property can be expressed in terms of $\tilde{C}_{h}$ as
\begin{equation}\label{eqn:skew_symmetry_discrete}
\tilde{C}_{h} (V_{h}) = - \tilde{C}_{h} (V_{h})^T.
\end{equation}
The pressure gradient term in \eqref{eqn:mom_continuous} is approximated by $G_{h} p_{h}$, where
\begin{equation}\label{eqn:divgrad_discrete}
\hl{G_{h} = -M_{h}^T.}
\end{equation}
\R{InfSup1}\hl{This relation between the divergence and gradient operator is a key aspect of our approach. Firstly, it allows us to formulate a velocity-only ROM, as will be detailed in section \ref{sec:momentum_conservation}. Secondly, the relation can be seen as a compatibility condition that ensures inf-sup stability \cite{Gallouet2018,Shin1997}, similar to the inf-sup condition typically encountered in finite element methods. Thirdly, it ensures that (like in the continuous case) kinetic energy is conserved by the pressure gradient term, since $V_{h}^T G_{h} p_{h} = -p_{h}^{T} M_{h} V_{h} = 0$. Relation \eqref{eqn:divgrad_discrete} is naturally satisfied by our choice of a finite volume method on a staggered grid, but can also be achieved with different discretisation methods, including collocated grids \cite{Trias2013a} and finite element methods \cite{Gresho1998}.} 

It is important to note that on inflow, no-slip or periodic boundaries, \textit{no boundary conditions are required for the pressure} -- they are implied by the boundary conditions for the velocity \cite{Sanderse2013,Veldman1990}. \R{PBC} \hl{In short, this is because on a staggered grid the pressure gradient is readily available at the interior velocity points, and on the boundaries the velocity is specified. For details we refer to  \ref{sec:BC_FOM}.}

The diffusive operator is discretised by second order central approximations, and can be represented by
\begin{equation}
D_{h} V_{h},
\end{equation}
where $D_{h}$ is a symmetric negative definite matrix, which can be written as $D_{h} = -Q_{h}^{T} Q_{h}$. 

The semi-discrete system then reads
\begin{align}
M_{h} V_{h} &= 0, \label{eqn:mass_semidiscrete} \\
 \Omega_{h}  \frac{\rd V_{h}}{\rd t} & = F_{h}(V_{h},p_{h}),\label{eqn:mom_semidiscrete}
\end{align}
where
$F_{h}(V_{h},p_{h})$ contains the convective, diffusive, pressure gradient, and body force contributions:
\begin{equation}
\hl{F_{h}(V_{h},p_{h}) =  F^{CD}_{h} (V_{h}) - G_{h} p_{h},}
\end{equation}
\begin{equation}
\hl{F^{CD}_{h}(V_{h}) = -\tilde{C}_{h}(V_{h})V_{h}+ \nu D_{h} V_{h} + f_{h}.}
\end{equation}
\R{Fht}\hl{$f_{h}$ is the discrete representation of the body-force (possibly time-dependent, but to ease notation the explicit time-dependence will be left out in $f_{h}$, $F_{h}$ and $F_{h}^{CD}$)}. $\Omega_{h}$ is a (time-independent) diagonal matrix with the finite volume sizes on its diagonal, which is symmetric positive definite. \hl{The pressure Poisson equation (PPE) follows by differentiating \eqref{eqn:mass_semidiscrete} and substituting \eqref{eqn:mom_semidiscrete}:
\begin{equation}\label{eqn:PPE}
L_{h} p_{h}  = M_{h} \Omega_{h}^{-1} F^{CD}_{h} (V_{h}), \qquad L_{h} = M_{h} \Omega_{h}^{-1} G_{h}.
\end{equation}
\R{PBC2}Like the momentum equations, the pressure Poisson equation (PPE) does not require the specification of pressure boundary conditions: the PPE is derived from the \textit{discrete} mass and momentum equations instead of the continuous counterparts, so that boundary conditions on the velocity suffice. This paradigm of working on the discrete level will also be applied when deriving the ROM. Furthermore, equation \eqref{eqn:PPE} will be the basis for the pressure recovery method in section \ref{sec:pressure_recovery}.}

%




To arrive at a semi-discrete energy equation a discrete inner product is needed, i.e.\ a discrete version of \eqref{eqn:innerproduct_continuous}. The natural choice in a finite volume context is (for $V_{h}$, $W_{h} \in \mathbb{R}^{N_{V}}$):
\begin{equation}
(V_{h},W_{h})_{\Omega_{h}} := V_{h}^{T} \Omega_{h} W_{h}, \qquad \|V_{h}\|_{\Omega_{h}}^{2} :=(V_{h},V_{h})_{\Omega_{h}}, \label{eqn:innerproduct_semidiscrete}
\end{equation}
and the discrete energy is defined as $K_{h} := \frac{1}{2} \|V_{h}\|_{\Omega_{h}}^2$. We will also need the unweighted norm $\| V_{h} \|^{2} := V_{h}^{T} V_{h}$.
In absence of boundary contributions, the time evolution of $K_{h}$ is given by
\begin{equation}
\begin{split}
2 \frac{\rd K_{h}}{\rd t} &= \frac{\rd }{\rd t} \left( V_{h}^{T} \Omega_{h} V_{h} \right) \\
&= -V_{h}^{T} (\tilde{C}_{h}(V_{h})^{T} + \tilde{C}_{h} (V_{h}) ) V_{h} - 2 p_{h}^{T} M_{h} V_{h} - 2\nu \|Q_{h} V_{h}\|^2.
\end{split}
\end{equation}
Due to the skew-symmetry property of $\tilde{C}_{h}$ (equation \eqref{eqn:skew_symmetry_discrete}) and the divergence-freeness of $V_{h}$ \eqref{eqn:mass_semidiscrete}, one obtains
\begin{equation}\label{eqn:energy_semidiscrete}
 \frac{\rd K_{h}}{\rd t} = - \nu \|Q_{h} V_{h} \|^2,
\end{equation}
which is the semi-discrete counterpart of equation \eqref{eqn:energy_continuous}.

%

\subsection{Time discretisation and fully discrete energy equation}\label{sec:temporal}
We continue with the temporal discretisation of equations \eqref{eqn:mass_semidiscrete}-\eqref{eqn:mom_semidiscrete} with an implicit $s$-stage Runge-Kutta method \cite{Sanderse2013b}. The stage values follow from
\begin{align}
M_{h} V_{h}^{n,i} &= 0, \label{eqn:mass_RKstage} \\
 \Omega_{h}  \frac{V_{h}^{n,i} - V_{h}^{n}}{\Delta t} &= \sum_{j=1}^{s} a_{ij} (\hl{F^{CD}_{h}}(V_{h}^{n,j}) - G_{h} p_{h}^{n,j}), \label{eqn:mom_RKstage}
\end{align}
and the solution at the next time step is a combination of the stage values:
\begin{align}
M_{h} V_{h}^{n+1} &= 0, \label{eqn:mass_RKstep} \\
 \Omega_{h}  \frac{V_{h}^{n+1} - V_{h}^{n}}{\Delta t} &= \sum_{i=1}^{s} b_{i} (\hl{F^{CD}_{h}}(V_{h}^{n,i}) - G_{h} p_{h}^{n,i}). \label{eqn:mom_RKstep}
\end{align}
Here, $V_{h}^{n}$ and $p_{h}^{n}$ are approximations to $V_{h}(t^{n})$ and $p_{h}(t^{n})$ respectively, which will be collected into a snapshot matrix to be used in the ROM construction (detailed in the next section).

The coefficients $a$ and $b$ of the Runge-Kutta method are chosen such that the temporal discretisation keeps the energy conservation property in the inviscid limit. An example of Runge-Kutta methods that satisfy this property is the family of Gauss methods \cite{Sanderse2013b}. We will employ the lowest order Gauss method, obtained for $s=1$, being the second order implicit midpoint method (corresponding to the following Butcher tableau: $a_{11} = \frac{1}{2}$, $b_{1}=1$).
For the Gauss methods, the fully discrete energy equation can be written as (in absence of body forces):
\begin{equation}\label{eqn:energy_fullydiscrete}
\frac{K_{h}^{n+1} - K_{h}^{n}}{\Delta t} = -\nu \sum_{i=1}^{s} \| Q_{h} V_{h}^{n,i} \|^2, 
\end{equation}
where $K_{h}^{n}$ is an approximation to $K_{h} (t^{n})$. In summary, the careful choice of spatial and temporal discretisation methods has yielded energy equations \eqref{eqn:energy_semidiscrete} and \eqref{eqn:energy_fullydiscrete} that closely mimic the continuous energy estimate \eqref{eqn:energy_continuous}. The fully discrete energy equation shows that, in the absence of boundary contributions, the energy of the solution can only decrease due to viscous dissipation, independent of the mesh, the time step, or the viscosity.

\hl{We note that these implicit methods require the solution of a non-linear system of equations with saddle-point structure, which can be computationally expensive. In section \ref{sec:alternative_time_integration} we will propose a more practical alternative, by using explicit time integration methods at the FOM level, without compromising the energy-conservation property of the ROM.}

%% file: figures/staggered_interior_2.pdf_tex
\begingroup%
  \makeatletter%
  \providecommand\color[2][]{%
    \errmessage{(Inkscape) Color is used for the text in Inkscape, but the package 'color.sty' is not loaded}%
    \renewcommand\color[2][]{}%
  }%
  \providecommand\transparent[1]{%
    \errmessage{(Inkscape) Transparency is used (non-zero) for the text in Inkscape, but the package 'transparent.sty' is not loaded}%
    \renewcommand\transparent[1]{}%
  }%
  \providecommand\rotatebox[2]{#2}%
  \newcommand*\fsize{\dimexpr\f@size pt\relax}%
  \newcommand*\lineheight[1]{\fontsize{\fsize}{#1\fsize}\selectfont}%
  \ifx\svgwidth\undefined%
    \setlength{\unitlength}{69.93857247bp}%
    \ifx\svgscale\undefined%
      \relax%
    \else%
      \setlength{\unitlength}{\unitlength * \real{\svgscale}}%
    \fi%
  \else%
    \setlength{\unitlength}{\svgwidth}%
  \fi%
  \global\let\svgwidth\undefined%
  \global\let\svgscale\undefined%
  \makeatother%
  \begin{picture}(1,0.75737714)%
    \lineheight{1}%
    \setlength\tabcolsep{0pt}%
    \put(0,0){\includegraphics[width=\unitlength,page=1]{staggered_interior_2.pdf}}%
    \put(0.46122931,0.69378898){\color[rgb]{0,0,0}\makebox(0,0)[lt]{\lineheight{1.25}\smash{\begin{tabular}[t]{l}$v_{i,j+1/2}$\end{tabular}}}}%
    \put(0.44452385,0.48308519){\color[rgb]{0,0,0}\makebox(0,0)[lt]{\lineheight{1.25}\smash{\begin{tabular}[t]{l}$p_{i,j}$\end{tabular}}}}%
    \put(0.6917468,0.38024506){\color[rgb]{0,0,0}\makebox(0,0)[lt]{\lineheight{1.25}\smash{\begin{tabular}[t]{l}$u_{i+1/2,j}$\end{tabular}}}}%
    \put(0.46122931,0.26484114){\color[rgb]{0,0,0}\makebox(0,0)[lt]{\lineheight{1.25}\smash{\begin{tabular}[t]{l}$v_{i,j-1/2}$\end{tabular}}}}%
    \put(0,0){\includegraphics[width=\unitlength,page=2]{staggered_interior_2.pdf}}%
    \put(0.41833452,0.00747243){\color[rgb]{0,0,0}\makebox(0,0)[lt]{\lineheight{1.25}\smash{\begin{tabular}[t]{l}$\Delta x$\end{tabular}}}}%
    \put(0,0){\includegraphics[width=\unitlength,page=3]{staggered_interior_2.pdf}}%
    \put(-0.00320838,0.42185649){\color[rgb]{0,0,0}\makebox(0,0)[lt]{\lineheight{1.25}\smash{\begin{tabular}[t]{l}$\Delta y$\end{tabular}}}}%
    \put(0,0){\includegraphics[width=\unitlength,page=4]{staggered_interior_2.pdf}}%
    \put(0.24135156,0.38024506){\color[rgb]{0,0,0}\makebox(0,0)[lt]{\lineheight{1.25}\smash{\begin{tabular}[t]{l}$u_{i-1/2,j}$\end{tabular}}}}%
  \end{picture}%
\endgroup%

%% file: ECPOD.tex
\section{Energy-conserving POD-Galerkin method}\label{sec:ECPOD}
\subsection{Introduction}
We will follow the ODE-based projection approach \cite{Benner2015} in which the POD-Galerkin method is applied to the semi-discrete energy-conserving formulation, i.e.\ we project the FOM given by equations \eqref{eqn:mass_semidiscrete} - \eqref{eqn:mom_semidiscrete}. This is not the only possibility; one can instead first project the continuous equations \eqref{eqn:divfree_continuous} and \eqref{eqn:mom_continuous} and then discretize the resulting system, \hl{common in finite element approaches \cite{Noack2005}}, or project the fully discrete equations \cite{Carlberg2018}. \R{DiscProj2}\hl{First discretising in space, and then performing the projection of the semi-discrete equations (instead of first projecting and then discretising) has several advantages. Firstly, the treatment of velocity boundary conditions is straightforward, and only has to be done once (when spatially discretising the FOM). \R{PBC4}Secondly, no pressure boundary conditions are needed, as they are not present in the discrete FOM. Furthermore, in combination with the staggered grid approach, the pressure disappears from the ROM, so that a velocity-only formulation results, and consequently the inf-sup condition is avoided on the level of the ROM.} 

We make the Ansatz that the velocity field $V_{h}(t) \in \mathbb{R}^{N_V}$ can be approximated by
\begin{equation}\label{eqn:ansatz}
V_{h}(t) \approx V_{r} (t) := \Phi a(t),
\end{equation}
where $\Phi \in \mathbb{R}^{N_V \times M}$, $a(t) \in \mathbb{R}^{M}$, and $M\ll N_{V}$. The subscript $r$ denotes quantities associated to the ROM. Equation \eqref{eqn:ansatz} is substituted into the FOM and then the equations are projected by left-multiplying with $\Phi^{T}$. In the POD approach $\Phi$ is obtained by performing a singular-value decomposition (SVD) of a snapshot matrix $X$ (this will be detailed below). $X$ contains $K$ snapshots of the velocity field $V_{h}$, i.e.\
\begin{equation}\label{eqn:snapshot_matrix}
X = [V_{h}^{1} \ldots V_{h}^{n} \ldots V_{h}^{K} ],
\end{equation}
where the snapshots are obtained from the solution of the fully discretised FOM, equations \eqref{eqn:mass_RKstep}-\eqref{eqn:mom_RKstep}. Each snapshot velocity field is divergence free, so that $M_{h} X_{j} = 0$ for each column $X_{j}$ of $X$. $\Phi$ is subject to the orthonormality condition 
$\Phi^{T} \Phi = I$.

\subsection{Construction of basis via weighted orthonormality condition}
In this work, instead of $\Phi^T \Phi=I$, we use a weighted orthonormality condition, namely
\begin{equation}\label{eqn:weighted_orthonormality}
\Phi^{T} \Omega_{h} \Phi = I.
\end{equation}
This is consistent with equation \eqref{eqn:innerproduct_semidiscrete} and with the form of the ROM momentum equation and the ROM kinetic energy equation, as we will demonstrate next. After substituting \eqref{eqn:ansatz} in \eqref{eqn:mom_semidiscrete} and projecting with $\Phi^{T}$, one obtains the reduced momentum equation
\begin{align}\label{eqn:mom_ROM_step2}
\Phi^{T} \Omega_{h} \frac{\rd \Phi a(t)}{\rd t} = \Phi^{T} F_{h}(\Phi a(t),p_{h}).
\end{align}
\hl{The treatment of the pressure in the right-hand side will be discussed in section \ref{sec:momentum_conservation}.}
In other words, it is natural to require $\Phi^T \Omega_{h} \Phi = I_{M}$ (the identity matrix of dimension $M \times M$), as the ROM  then simplifies to
\begin{equation}\label{eqn:mom_ROM}
\frac{\rd a(t)}{\rd t} = \Phi^{T} F_{h}(\Phi a(t),p_{h}).
\end{equation}
\hl{In practice, the right-hand side $\Phi^{T} F_{h}(\Phi a(t),p_{h})$ is rewritten as $F_{r}(a(t))$, where the ROM operator $F_{r}$ is precomputed offline. This will be discussed in section \ref{sec:offline_decomposition}.}
With condition \eqref{eqn:weighted_orthonormality}, the energy of the FOM is approximated by
\begin{equation}\label{eqn:energy_ROM}
K_{h}(t) \approx K_{r} (t) = \frac{1}{2} (\Phi a(t))^{T} \Omega_{h} \Phi a(t) = \frac{1}{2} a(t)^{T} \Phi^{T} \Omega_{h} \Phi a(t) = \frac{1}{2} a(t)^{T} a(t) = \frac{1}{2} \| a \|^2.
\end{equation}
The choice \eqref{eqn:weighted_orthonormality} thus simplifies the expression for both momentum and energy considerably.

We now specify the construction of $\Phi$. Given that the energy norm is chosen to be based on \eqref{eqn:energy_ROM}, the basis $\Phi$ should be computed from the following minimization problem \cite{Afkham2018a,Quarteroni2016}:
\begin{equation}\label{eqn:minimization_problem}
\Phi = \argmin_{\Phi} \| (I - \Phi \Phi^{T} \Omega_{h}) X \|^{2}_{F} \qquad \text{subject to} \qquad \Phi^T \Omega_{h} \Phi = I_{M},
\end{equation}
instead of the classical `unweighted' minimization problem:
\begin{equation}
\hat{\Phi} = \argmin_{\Phi} \| (I - \Phi \Phi^{T}) X \|^{2}_{F} \qquad \text{subject to} \qquad \hat{\Phi}^T \hat{\Phi} = I_{M}.
\end{equation}
The solution of the weighted problem can be expressed in terms of the solution of the unweighted problem as follows. Let
\begin{equation}
\Phi = \Omega_{h}^{-1/2} \hat{\Phi},
\end{equation}
where $\hat{\Phi}$ follows from the SVD of the scaled snapshot matrix $\hat{X} = \Omega_{h}^{1/2} X$:
\begin{equation}
\hat{X} =  \hat{\Phi} \Sigma \Psi^{*}.
\end{equation}
Since $\Omega_{h}$ is a diagonal matrix, its matrix square root is trivial to compute. The dimensions of the matrices in the SVD are 
\begin{equation}
\hat{\Phi}, \Phi \in \mathbb{R}^{N_V \times N_V}, \qquad \Sigma \in \mathbb{R}^{N_V \times K}, \qquad  \Psi \in \mathbb{R}^{K \times K}.
\end{equation}
The columns of $\hat{\Phi}$, denoted by $\hat{\Phi}_{j}$, are the eigenvectors of the correlation matrix $\hat{X} \hat{X}^{T}$, i.e.
\begin{equation}\label{eqn:eigenvectors}
\hat{X} \hat{X}^{T} \hat{\Phi}_{j} = \lambda_{j} \hat{\Phi}_{j},
\end{equation}
where the eigenvalues $\lambda$ are related to the singular values $\sigma$ (diagonal entries of $\Sigma$) by $\sqrt{\lambda_{i}(\hat{X} \hat{X}^{T})} = \sigma_{i} (\hat{X})$.

The basis for the ROM is obtained by taking (truncating) the first $M$ columns of $\Phi$. $M$ is typically prescribed by analysing the decay of the singular values $\sigma$.  

In summary, the sequence to obtain $\Phi$ is: gather snapshots of the velocity field in $X$; compute $\hat{X}$; compute the SVD of $\hat{X}$ to get $\hat{\Phi}$; compute $\Phi$; truncate $\Phi$. \hl{Note that for the (common) case that $N_{V} \gg K$, one can alternatively use the method of snapshots \cite{Sirovich1987} in which an eigenvalue problem for $\hat{X}^{T} \hat{X}$ is solved, instead of computing the SVD of $\hat{X}$.}

\subsection{Mass conservation of the ROM}
It is well-known that the mass conservation equation is identically satisfied by the ROM approximation, if the boundary conditions are no-slip or periodic \cite{Noack2016}. In a finite volume context, this is shown as follows. The divergence-free condition \eqref{eqn:mass_semidiscrete} becomes
\begin{equation}\label{eqn:ROM_divfree}
M_{h} \Phi a(t) = 0.
\end{equation}
Rewriting equation \eqref{eqn:eigenvectors} yields
\begin{equation}
X X^{T} \Omega_{h} \Phi_{j} = \lambda_{j} \Phi_{j}.
\end{equation}
Left-multiplying this equation with $M_{h}$ and using that the snapshots collected in $X$ are divergence-free ($M_{h} X_{j} = 0$) yields
\begin{equation}
 \lambda_{j} M_{h} \Phi_{j} = 0.
\end{equation}
In other words, the ROM velocity field $V_{r}=\Phi a$ satisfies the divergence-free condition \eqref{eqn:ROM_divfree}, independent of the value of the coefficients $a(t)$, as long as $\lambda_{j}\neq0$. Note that for non-homogeneous boundary conditions this is not the case. We will present a boundary condition treatment in section \ref{sec:BC}.

\subsection{Momentum conservation of the ROM}\label{sec:momentum_conservation}
The ROM momentum equation was given by equation \eqref{eqn:mom_ROM}:
\begin{equation}\label{eqn:mom_ROM_ext}
\begin{split}
\frac{\rd a(t)}{\rd t} &= \Phi^{T} F_{h}(\Phi a(t),p_{h}), \\
&=\Phi^{T} (-\tilde{C}_{h} (\Phi a(t)) \Phi a(t) + \nu D_{h} \hl{\Phi} a(t) + f_{h} - G_{h} p_{h}),  \\
&= \hl{\Phi^{T} (-\tilde{C}_{h} (\Phi a(t)) \Phi a(t) + \nu D_{h} \hl{\Phi} a(t) + f_{h})}, \\
&= \hl{C_{r} (a(t)) + \nu D_{r} a(t) + f_{r}}, \\
&= \hl{F_{r}(a(t))},
\end{split}
\end{equation}
\hl{where the convective and diffusive ROM operators $C_{r}$ and $D_{r}$, and the body force vector $f_{r}$ are precomputed in an offline setting (see section \ref{sec:offline_decomposition})}.

The pressure gradient term disappears because the pressure gradient is linked to the divergence operator according to equation \eqref{eqn:divgrad_discrete}: 
\begin{equation}
\Phi^{T} G_{h} = (G_{h}^{T} \Phi)^{T} = - (M_{h} \Phi)^{T} = 0.
\end{equation}
\R{InfSup2}\hl{We stress that the pressure term is not `neglected' (in the sense of being an approximation); it is really identically equal to zero. This is  accomplished via two crucial properties: (i) the FOM velocity snapshots (and thus the basis $\Phi$) are divergence-free, and (ii) the divergence operator and pressure gradient operator are compatible according to equation \eqref{eqn:divgrad_discrete}. In other words, the ROM velocity field is by construction divergence-free, and the pressure (a Lagrange multiplier) is not needed anymore to enforce this property. The disappearance of the pressure term leads to a `velocity-only' ROM. Consequently, potential issues of inf-sup stability on the reduced order level are avoided. In section \ref{sec:ECROM} it will furthermore be shown that this velocity-only ROM possesses non-linear stability.}

\hl{This `velocity-only' formulation is also possible in the case of non-homogeneous boundary conditions, including the important case of open flows. This will be discussed in section \ref{sec:BC}. If the pressure field is desired, it can be recovered via a post-processing step, as will be detailed in section \ref{sec:pressure_recovery}.}


Momentum is, unlike mass, not a locally conserved quantity. However, momentum \textit{is} globally conserved in case of periodic boundary conditions: integration of the incompressible Navier-Stokes equations over a domain $\Omega$ with periodic boundary conditions yields
\begin{equation}
\frac{\rd \vt{P}(t) }{\rd t} = 0, \qquad \text{where} \qquad \vt{P}(t) = \int_{\Omega} \vt{u} \, \rd \Omega,
\end{equation}
so that momentum $\vt{P}$ is exactly conserved in time. We will require the ROM to satisfy this property. Define the FOM global momentum of each velocity component as
\begin{align}
P_{h}^{u}(t) = e_{u}^{T} \Omega_{h} V_{h}(t),\\
P_{h}^{v}(t) = e_{v}^{T} \Omega_{h} V_{h}(t),
\end{align}
where $e_{u},e_{v} \in \mathbb{R}^{N_{V}}$. $e_{u}$ contains a $1$ for indices associated with the $u$-velocity component, and $e_{v}$ contains a $1$ for indices associated with the $v$-velocity component, such that $e = e_{u} + e_{v} = [1, 1, \ldots, 1]^{T}$.
Evolution of the $u$-component of the FOM global momentum is given by
\begin{align}
\frac{\rd P_{h}^{u}}{\rd t} = e_{u}^{T} \Omega_{h} \frac{\rd V_{h}(t)}{\rd t} = e_{u}^{T} F_{h}(V_{h},p_{h}) = 0,
\end{align}
with a similar expression for $P_{h}^v$. This expression evaluates to zero because of the telescoping property of finite volume methods in combination with a periodic domain.

Evolution of global momentum predicted by the ROM reads
\begin{align}\label{eqn:ROM_dPdt}
\frac{\rd P_{h}^{u}}{\rd t} \approx \frac{\rd P_{r}^{u}}{\rd t} = e_{u}^{T} \Omega_{h} \Phi \frac{\rd a(t)}{\rd t} = e_{u}^{T} \Omega_{h} \Phi \Phi^{T} F_{h}(\Phi a(t),p_{h}).
\end{align}
In order to obtain global conservation of momentum of the ROM, we enforce the basis vectors $\Phi$ to satisfy
\begin{align}
 e_{u}^{T} \Omega_{h} \Phi \Phi^{T} = e_{u}^{T}, \label{eqn:SVD_constraint1}
\end{align}
with a similar expression for $e_v$. In other words, the projection of the vectors $e_{u}$ and $e_{v}$ by $\Omega_{h} \Phi \Phi^{T}$ should be exact. When performing the SVD without truncation this property can be easily achieved by adding these vectors to the snapshot matrix, since a property of the SVD is that the projection of vectors in the snapshot matrix is exact. However, upon truncating the decomposition to arrive at a reduced dimension, this property is generally lost. We therefore propose to enforce property \eqref{eqn:SVD_constraint1} by using an approach \R{constrainedSVD}\hl{that enforces linear constraints in constructing the POD basis, as detailed in \cite{Xiao2014,Xiao2013}. In short, the idea is to change the minimization problem that underlies the SVD (equation \eqref{eqn:minimization_problem}) by a constrained minimization problem, which can be efficiently solved by using a QR-decomposition of the constraint matrix and a projected snapshot set. In the current setting, the constraint matrix is particularly simple, and we can avoid computing the QR-decomposition with the following construction. Furthermore, our analysis is for the weighted norm, instead of the unweighted case in \cite{Xiao2013}.} First, we collect the vectors that should be exactly projected by the truncated SVD in the matrix $E$:
\begin{align}
E = [e_{u} \, e_{v}],
\end{align}
scaled to have norm equal to $1$. \hl{$E$ represents the constraint matrix, which is called $\vt{G}$ in the notation of \cite{Xiao2013}}. Next, we perform an update of the snapshot matrix $X$:
\begin{equation}
\tilde{X} = X - E E^T \Omega_{h} X,
\end{equation}
and we determine its SVD 
\begin{equation}
\tilde{X} = \tilde{\Phi} \tilde{\Sigma} \tilde{\Psi}^{*}.
\end{equation}
Subsequently, we add $E$ to $\tilde{\Phi}$, and then we truncate:
\begin{equation}
\Phi = [E  \, \tilde{\Phi}]_{M},
\end{equation}
where the subscript $M$ indicates that the first $M$ columns are used. The resulting $\Phi$ satisfies equation \eqref{eqn:SVD_constraint1} (and a similar equation for $e_{v}$). The proof is given in \ref{sec:proof_adapted_norm}. Note that, when $M$ is given, enforcing global momentum comes at the price of losing two or three of the modes present in $\tilde{\Phi}$ (depending on the spatial dimension of the problem).

The initial condition for $a$ is given by
\begin{equation}
a(0) = \Phi^{T} \Omega_{h} V_{h}(0).
\end{equation}
Consequently, the initial momentum ($u$-component) is given by
\begin{equation}
P^{u}_{r}(0) = e_{u}^{T} \Omega_{h} \Phi a(0) = e_{u}^{T} \Omega_{h} \Phi \Phi^{T} \Omega_{h} V_{h}(0),
\end{equation}
whereas the initial momentum of the FOM is
\begin{equation}
P^{u}_{h}(0) = e_{u}^{T} \Omega_{h} V_{h}(0).
\end{equation}
The error between the two is
\begin{equation}\label{eqn:mom_projection_error}
P^{u}_{r}(0) - P^{u}_{h}(0) = - e_{u}^{T} (I - \Omega_{h} \Phi \Phi^{T}) \Omega_{h} V_{h}(0) = 0,
\end{equation}
when the constrained SVD is employed: the total momentum of the ROM is constant in time and equals the total momentum of the FOM, for the case of periodic boundary conditions. The same equation holds for the $v$-component. \hl{Upon time discretisation, no additional measures need to be taken to ensure that momentum stays conserved, since it is a linear constraint which is conserved by any consistent time integration method. The fully discrete momentum error is thus:
\begin{equation}\label{eqn:discrete_mom_eqn_ROM}
\epsilon_{P} := P^{u,n}_{r} - P^{u,n}_{h} = \underbrace{P^{u,n}_{r} - P^{u}_{r}(0)}_{\text{time integration error ROM}} + \underbrace{P^{u}_{r}(0) - P^{u}_{h}(0)}_{\text{projection error, \eqref{eqn:mom_projection_error}}}  + \underbrace{P^{u}_{h}(0) - P^{u,n}_{h}}_{\text{time integration error FOM}} = 0.
\end{equation}
We note that if the constrained SVD is \textit{not} employed, the third term remains 0, but the first and second term both become nonzero. The first term scales with the right-hand side of \eqref{eqn:ROM_dPdt} (and the order of accuracy of the ROM time integration). The second term scales with $e_{u}^{T} (I - \Omega_{h} \Phi \Phi^{T}) \Omega_{h} V_{h}(0)$.}

Note that this global momentum-conserving construction is different from the conservative model reduction method from Carlberg et al.\ \cite{Carlberg2018}. We consider global momentum conservation for the case of periodic boundary conditions, by adapting the construction of the SVD, whereas in \cite{Carlberg2018} a constrained optimization problem is considered, which minimizes the residual of the full-order model over subdomains. In our view, conservation means that the integral of a certain quantity (mass, momentum, energy) stays invariant in time, which only holds for particular boundary conditions; in \cite{Carlberg2018} the term conservation is used to indicate the difference between the rate of change of a conserved quantity and the contribution of surface integrals and source terms. Furthermore, our method differs in the fact that we are not only considering primary conserved quantities, but are considering kinetic energy conservation (a so-called secondary or `derived' quantity) which has the advantage that it provides a non-linear stability bound to the solution. This is detailed in the next section.

\subsection{Energy conservation of the ROM}\label{sec:ECROM}
One of the key questions in this paper is whether the kinetic energy of the ROM can be bounded in a similar way as the energy of the FOM (equations \eqref{eqn:energy_continuous}, \eqref{eqn:energy_semidiscrete}, \eqref{eqn:energy_fullydiscrete}). To this end, we differentiate the expression for the energy of the ROM, as given by equation \eqref{eqn:energy_ROM}, and simplify by using equations \eqref{eqn:ROM_divfree} and \eqref{eqn:mom_ROM_ext} (in absence of body forces):
\begin{equation}
\begin{split}
2 \frac{\rd K_{r}(t) }{\rd t} &= \frac{\rd a^{T}}{\rd t} a + a^{T}  \frac{\rd a}{\rd t}, \\
&= -(\Phi{^T} \tilde{C}_{h}(\Phi a) \Phi a)^{T} a - a^{T} (\Phi^{T} \tilde{C}_{h}(\Phi a) \Phi a) + \nu (\Phi^{T} D \Phi a)^{T} a + a^{T} \nu \Phi^{T} D_{h} \Phi a,\\
&= - a^{T} \Phi{^T} (\tilde{C}_{h}(\Phi a)^{T} + \tilde{C}_{h}(\Phi a)) \Phi a - 2 \nu \|Q_{h} \Phi a \|^2, \\
&= - 2 \nu \|Q_{h} \Phi a \|^2,\\
&= \hl{- 2 \nu \|Q_{r} a \|^2.}
\end{split}
\end{equation}
The crucial steps in the derivation are the fact that the basis is divergence-free ($M_{h} \Phi=0$) and that the following properties of the FOM spatial discretisation operators hold: $G_{h}^{T} = -M_{h}$; $\tilde{C}_{h}(\Phi a)$ is skew-symmetric; $D_{h}$ is symmetric negative definite. \hl{In terms of the offline ROM operators (detailed in section \ref{sec:offline_decomposition}) one can equivalently state that each element $C_{r,i}$ of the convective operator $C_r$ is skew-symmetric, and that the diffusive operator $D_{r}$ is symmetric negative definite. This last observation allows us to write $D_{r}$ as $-Q_{r}^T Q_{r}$, with $Q_{r} =Q_{h} \Phi$.}

In summary, the energy evolution of the ROM in absence of boundary contributions and body forces is given by 
\begin{equation}\label{eqn:energy_eqn_ROM}
\boxed{\frac{\rd K_{r}(t)}{\rd t} = - \nu \|\hl{Q_{r}} a \|^{2}.}
\end{equation}
Consequently, \textit{the ROM is non-linearly stable, independent of the number of POD modes used}.

In the inviscid limit, the energy of the ROM stays constant in time:
\begin{equation}\label{eqn:energy_rom_inviscid}
K_{r} (t) = K_{r} (0),
\end{equation}
where
\begin{equation}
K_{r}(0) = \frac{1}{2} a(0)^{T} a(0) = \frac{1}{2} (\Phi^{T} \Omega_{h} V_{h}(0))^{T} \Phi^{T} \Omega_{h} V_{h}(0) = \frac{1}{2} V_{h}(0)^{T} \Omega_{h} \Phi \Phi^{T} \Omega_{h} V_{h}(0).
\end{equation}
The FOM kinetic energy is given by
\begin{equation}
K_{h} (0) = \frac{1}{2} V_{h}(0)^{T} \Omega_{h} V_{h}(0),
\end{equation}
and the error in the ROM (due to truncation) is therefore given by
\begin{equation}\label{eqn:energy_rom_inviscid_initial}
K_{r} (0) - K_{h}(0) = - \frac{1}{2} V_{h}(0)^{T} (I - \Omega_{h} \Phi \Phi^{T}) \Omega_{h} V_{h}(0).
\end{equation}
\R{EnergyConstraint}\hl{Similar to the addition of the global momentum constraint, it is possible to add $V_{h}(0)$ to the truncated SVD so that it is projected exactly by $\Omega_{h} \Phi \Phi^{T}$. This guarantees that in the inviscid case the kinetic energy of the ROM remains equal to that of the FOM. However, for viscous simulations, the ROM will not reproduce the kinetic energy evolution exactly (compare equation \eqref{eqn:energy_eqn_ROM} to \eqref{eqn:energy_semidiscrete}), and so we do not enforce this property.}

The last step in obtaining the ROM is to specify a time discretisation for equation \eqref{eqn:mom_ROM} such that a fully discrete equivalent of \eqref{eqn:energy_ROM} is obtained. The key is, not surprisingly, to use the energy-conserving Runge-Kutta time discretisation methods introduced in section \ref{sec:spatial}. For example, the implicit midpoint method applied to \eqref{eqn:mom_ROM} reads
\begin{equation}\label{eqn:ROM_IM}
\frac{a^{n+1} - a^{n}}{\Delta t} = \hl{F_{r}} (a^{n+1/2}),
\end{equation}
where $a^{n+1/2} = \frac{1}{2} (a^{n} + a^{n+1})$. \hl{This non-linear equation is solved efficiently by using Newton's method; this requires the Jacobian of $F_{r}(a)$, which is given in section \ref{sec:offline_decomposition}.} The corresponding energy evolution is
\begin{equation}\label{eqn:discrete_energy_eqn_ROM}
\boxed{\frac{K_{r}^{n+1} - K_{r}^{n}}{\Delta t} = -\nu \| \hl{Q_{r}} a^{n+1/2} \|^2,}
\end{equation}
which is strictly decreasing in time when the viscosity $\nu$ is nonzero, and hence \textit{the fully discrete ROM solution is stable}.

\hl{The energy equality \eqref{eqn:discrete_energy_eqn_ROM} is derived for homogeneous (no-slip, periodic) boundary conditions and in the absence of body forces. For more generic boundary conditions, such as inflow conditions, the energy of the flow is not strictly decreasing in time. In this case, the ROM estimate should still be such that it mimics the energy estimate of the FOM (see e.g.\ \cite{Sanderse2013} for the boundary contributions to the energy equation). Further details on general boundary conditions are given in section \ref{sec:BC}.}

\hl{The energy-conserving property expressed by equations \eqref{eqn:energy_eqn_ROM} and \eqref{eqn:discrete_energy_eqn_ROM} is \textit{independent} of whether the snapshot matrix has been generated using an energy-conserving discretisation method. The only condition on the snapshot matrix is that the snapshots are divergence-free. Of course, when considering the error of the fully discrete ROM with respect to the fully discrete FOM, the FOM energy error \textit{is} present:
\begin{equation}\label{eqn:energy_error_decomposition}
\epsilon_{K} := K_{r}^{n} - K_{h}^{n}  = \underbrace{K_{r}^{n} - K_{r}(0)}_{\text{time integration error ROM}} + \underbrace{K_{r}(0) - K_{h}(0)}_{\text{projection error}} + \underbrace{K_{h}(0)- K_{h}^{n}}_{\text{time integration error FOM}}.
\end{equation}
The first term of this expansion is zero when employing the energy-conserving ROM satisfying \eqref{eqn:discrete_energy_eqn_ROM} (for $\nu=0$). The second term is generally nonzero and depends on the number of modes used and the decay of the singular values. The last term is zero when an energy-conserving time integration method is used on the FOM level.}

%% file: BC.tex
\section{Non-homogeneous boundary conditions}\label{sec:BC}
In this section we extend the results of the previous section to the more generic case of (stationary) non-homogeneous boundary conditions and forcing terms, which is known to be non-trivial (see for example \cite{Lassila2014,Noack2005,Weller2009}). One of the main issues is the pressure term, whose contribution does not vanish \cite{Noack2005}; in the kinetic energy equation, this gives the term
\begin{equation}
(\vt{u},\nabla p) = \int_{\partial \Omega} p \, \vt{u} \cdot \vt{n} \, \rd  S,
\end{equation}
and in the Galerkin projection a similar term appears. In the spirit of \cite{Balajewicz2013,Couplet2005,Gunzburger2007,Holmes2012}, we will split the velocity field into two components: a time-dependent field $\vt{u}_{\text{hom}}(\vt{x},t)$ that satisfies homogeneous boundary conditions, and a stationary field $\vt{u}_{\text{bc}}(\vt{x})$ that satisfies non-homogeneous boundary conditions, 
\begin{equation}
\vt{u} (\vt{x},t) = \vt{u}_{\text{hom}}(\vt{x},t) +  \vt{u}_{\text{bc}}(\vt{x}).
\end{equation}
The non-homogeneous term is also known as a lifting function \cite{Fonn2019,Hijazi2020}. Both terms satisfy $\nabla \cdot \vt{u}_{\text{hom}}(\vt{x},t) = 0$ and $\nabla \cdot \vt{u}_{\text{bc}}(\vt{x}) = 0$. In the case of Dirichlet boundary conditions we have $\vt{u}_{\text{bc}}(\vt{x}) = \vt{u}_{\partial \Omega}$ and $\vt{u}_{\text{hom}}(\vt{x}) = \vt{0}$ on $\partial \Omega$. The stationary field is found from $\vt{u}_{\text{bc}} = \nabla q$, where $q(\vt{x})$ follows from solving the Poisson equation 
\begin{equation}\label{eqn:Poisson_bc}
\nabla^2 q(\vt{x}) = 0 \qquad \text{with} \qquad \nabla q = \vt{u}_{\partial \Omega} \quad \text{on} \quad \partial \Omega.
\end{equation}
Subsequently, with this construction in the kinetic energy equation for $\vt{u}_{\text{hom}}$ the pressure term contribution still disappears:
\begin{equation}
(\vt{u}_{\text{hom}},\nabla p) = \int_{\partial \Omega} p \, \vt{u}_{\text{hom}} \cdot \vt{n} \, \rd  S = 0.
\end{equation}

With this insight, a ROM formulation incorporating non-homogeneous boundary conditions is constructed, but starting from the semi-discrete instead of the continuous equations. The semi-discrete equations \eqref{eqn:mass_semidiscrete}-\eqref{eqn:mom_semidiscrete} with non-homogeneous boundary conditions read:
\begin{align}
M_{h} V_{h}(t) &= \hl{y_{M}}, \label{eqn:mass_semidiscrete_bc} \\
 \Omega_{h}  \frac{\rd V_{h}(t)}{\rd t} & = F^{CD}_{h} (V_{h}(t)) - (G_{h} p_{h}(t) + y_{G}), \label{eqn:mom_semidiscrete_bc} 
 \end{align}
with
 \begin{equation}
 \hl{F^{CD}_{h} (V_{h}(t)) =- C_{h}(V_{h}(t),V_{h}(t)) + \nu (D_{h} V_{h}(t) + y_{D}) + f_{h}.}
 \end{equation}
The boundary terms \hl{such as $y_{M} \in \mathbb{R}^{N_{p}}$, $y_{G}, y_{D} \in \mathbb{R}^{N_{V}}$} are described for example in \cite{Sanderse2013} and \hl{summarized in \ref{sec:BC_FOM}}. This formulation holds for both inflow, outflow, symmetry, periodic, and no-slip conditions, and also encompasses the case of body forces via the term $f_{h}$. \hl{In order to retain a velocity-only ROM formulation, the two properties mentioned in section \ref{sec:momentum_conservation} need to be satisfied: the basis $\Phi$ should be divergence-free, and the pressure gradient and divergence operator should be compatible.}

\hl{The first property (divergence-free basis) is achieved as follows.} We approximate the FOM velocity field with a ROM with homogeneous boundary conditions and a term that incorporates the boundary conditions \cite{Gunzburger2007}:
\begin{equation}\label{eqn:VROM_bc}
V_{h}(t) \approx V_{r}(t) + V_{bc} = \Phi a(t) + V_{bc}.
\end{equation}
\hl{This is a generalization of the ROM in section \ref{sec:ECPOD}, which corresponds to the case $V_{bc}=0$.} $V_{bc}$ is chosen such that 
\begin{equation}
M_{h} (V_{r}(t) + V_{bc}) = \hl{y_{M}},
\end{equation}
which reduces to $M_{h} V_{bc} = \hl{y_{M}}$, since $V_{r}(t)$ satisfies homogeneous boundary conditions. The solution for $V_{bc}$ is the discrete version of equation \eqref{eqn:Poisson_bc}, namely
\begin{equation}\label{eqn:poisson_Vbc}
V_{bc} = \Omega_{h}^{-1} G_{h} \zeta_{h}, \qquad \text{where} \qquad \hl{L_{h}} \zeta_{h} = \hl{y_{M}}.
\end{equation}
Thus, the solution of one Poisson equation \hl{at the FOM level} is needed to find the boundary terms $V_{bc}$, \hl{which is best performed during the offline FOM stage, when the Poisson equation is already being solved many times}. This $V_{bc}$ is subtracted from the velocity snapshots $V_{h}$ in order to arrive at the snapshot matrix of $V_{r}$:
\begin{equation}
X = [V_{h}^{1} - V_{bc} \ldots V_{h}^{n} - V_{bc} \ldots V_{h}^{K} - V_{bc}],
\end{equation}
which satisfies $M_{h} X_{j} = 0$ for each column $j$ of $X$, and the projection matrix $\Phi$ thus satisfies $M _{h}\Phi_{j} = 0$. 

\hl{The second property (compatibility of divergence and gradient operators) is naturally achieved on a staggered grid when the boundary conditions are periodic or no-slip (the extension to inflow conditions and symmetry conditions is straightforward \cite{Sanderse2013}). \R{PBC5}The case of outflow conditions is more involved. We prescribe the normal component of the stress tensor \cite{Caiazzo2014}:
\begin{equation}\label{eqn:outflow_BC}
(-p \vt{I} + \nu \nabla \vt{u}) \cdot \vt{n} = \begin{pmatrix}
- p_{\infty} \\ 0
\end{pmatrix},
\end{equation}
where $p_{\infty}(\vt{x},t)$ can be space- and time-dependent (it will be taken 0 in our open-flow test case). The finite volumes adjacent to the outflow boundary feature this term. Consequently, the expression for the pressure gradient is
\begin{equation}
G_{h} p_{h} (t) + y_{G},
\end{equation} 
where $y_{G} \in \mathbb{R}^{N_{V}}$ contains the value of $p_{\infty}$ in the entries corresponding to an outflow boundary volume. The divergence-gradient relation $G_{h}=-M_{h}^{T}$ (equation \eqref{eqn:divgrad_discrete}) is then still satisfied (for details, see \cite{Sanderse2013}). In the ROM, the pressure contribution thus becomes
\begin{equation}
\Phi^{T} (G_{h} p_{h}(t) + y_{G}) = - p_{h}^T(t) M_{h} \Phi  + \Phi^{T} y_{G} = \Phi^{T} y_{G},
\end{equation}
which can be easily precomputed. Consequently, the ROM formulation remains `velocity-only', just like in the homogeneous case:
\begin{equation}\label{eqn:mom_ROM_bc}
\Phi^{T} \Omega_{h} \frac{\rd }{\rd t} \left( \Phi a(t) + V_{bc} \right) = \frac{\rd a(t)}{\rd t} = F_{r}( a(t)), 
\end{equation}
where precomputing $F_{r}$ is now slightly more involved due to the presence of $V_{bc}$:
\begin{equation}
F_{r}(a(t)) = \Phi^{T} \left( - C_{h} (\Phi a(t) + V_{bc},\Phi a(t) +V_{bc}) + \nu (D_{h} (\Phi a(t) + V_{bc}) +y_{D})  - y_{G} + f_{h} \right).
\end{equation}
This construction is detailed in \ref{sec:ROM_Vbc}.}
Note that for the case of time-dependent or parameter-dependent boundary conditions appearing in $y_{M}$, the method needs to be extended by incorporating pressure snapshots and a projection of the divergence-free constraint. \hl{Although this is the subject of future work, in section \ref{sec:pressure_recovery} we will give details on the projection of the divergence-free constraint in order to formulate a pressure Poisson equation for the recovery of the pressure.}

%% file: implementation.tex
\section{\hl{Implementation aspects and pressure recovery}}\label{sec:implementation}

\subsection{\hl{Offline decomposition}}\label{sec:offline_decomposition}
\hl{The semi-discrete and fully discrete ROM momentum equations feature the right-hand side (see equations \eqref{eqn:mom_ROM} and \eqref{eqn:ROM_IM}) 
\begin{equation}\label{eqn:mom_ROM2}
F_{r} (a) := \Phi^{T} F_{h}(\Phi a) = \Phi^{T} ( - \tilde{C}_{h}(\Phi a) \Phi a + \nu D_{h} \Phi a + f_{h}).
\end{equation}
To be efficient as a ROM, this term should be evaluated such that it does not depend on operations on the FOM-level. This is performed by precomputing the convective and diffusive operators as follows. For ease of notation, we discuss the case of homogeneous boundary conditions. The non-homogeneous case is discussed in \ref{sec:ROM_bc}.}

\hl{The projection of the (linear) diffusive operator is straightforward and can be written as 
\begin{equation}
\Phi^{T} D_{h} (\Phi a) = D_{r} a,
\end{equation}
where $D_{r} = \Phi^{T} D_{h} \Phi \in \mathbb{R}^{M \times M}$.} 
\hl{The projection of the nonlinear convective operator is more involved, but given that the nonlinearity is only quadratic, it can be written as
\begin{equation}
\Phi^{T} \tilde{C}_{h} (\Phi a) (\Phi a) = C_{r} (a \otimes a) = \begin{bmatrix}
C_{r,1} & \ldots & C_{r,M} 
\end{bmatrix}
(a \otimes a),
\end{equation}
where $C_{r} \in \mathbb{R}^{M \times M^2}$ is a ``matricized" third order tensor, consisting of $M$ components $C_{r,i} \in \mathbb{R}^{M \times M}$:
\begin{equation}
C_{r,i} = \Phi^{T} \tilde{C}_{h} (\Phi_{i}) \Phi.
\end{equation}
\R{hyperreduction}In the current work, similar to \cite{Stabile2018}, the number of modes that we require for an accurate ROM is sufficiently low to warrant that storing and pre-computing this third order tensor is not prohibitive. If many modes would be required, if higher-order or non-polynomial nonlinearities are present in the equations (e.g.\ present in a source term), or if the discretisation of the quadratic term is nonlinear (e.g.\ a flux-limiting scheme), precomputing as high-order tensors is not possible or not efficient. In such cases, one can resort for example to the discrete empirical interpolation method (DEIM) \cite{Chaturantabut2010}, or a lifting procedure \cite{Kramer2019}. Designing a DEIM procedure such that the skew-symmetry of the convective operator is retained will be subject of future work.}

\hl{The entire semi-discrete ROM formulation for both homogeneous and non-homogeneous boundary conditions can be summarized as
\begin{equation}
\frac{\rd a}{\rd t} = F_{r} (a) = F_{r,2} (a \otimes a) + F_{r,1} a + F_{r,0},
\end{equation}
where $F_{r,0} \in \mathbb{R}^{M}$ contains all terms that are independent of the solution (boundary conditions, source terms), $F_{r,1} \in \mathbb{R}^{M \times M}$ contains linear terms (the diffusive operator and the interaction of the convective operator with the inhomogeneous boundary conditions), and $F_{r,2} \in \mathbb{R}^{M \times M^2}$ contains the quadratic terms of the convective operator. The Jacobian of $F_{r}$ with respect to $a$, as needed for implicit time integration (equation \eqref{eqn:ROM_IM}), is given by
\begin{equation}\label{eqn:Jacobian}
J_{r} = \frac{\partial F_{r}}{\partial a} = F_{r,2} (  I_{M} \otimes a + a \otimes I_{M}) + F_{r,1},
\end{equation}
where $I_{M}$ is the $M\times M$ identity matrix.}

\subsection{\hl{Alternative time integration methods and linear stability}}\label{sec:alternative_time_integration}
From an implementation point of view, the implicit time integration methods proposed (on both the FOM and ROM level) might seem daunting at first sight. \hl{However, one should note that the Jacobian \eqref{eqn:Jacobian} is rather simple to form (because the nonlinearity is only quadratic). Furthermore, because the pressure is absent in the ROM, one does \textit{not} need to solve a coupled system of saddlepoint type, as is the case on the FOM level.} 

\hl{In any case, one can think of test cases in which explicit time integration (at FOM and/or ROM level) is more efficient than implicit time integration.} In that case, we propose the following strategy. First, solve the FOM with a time integration method of choice (e.g.\ explicit, IMEX, etc.)\ -- this does \textit{not} affect energy conservation of the ROM as long as the FOM spatial discretization is energy-conserving. Second, use linear stability theory to estimate the eigenvalues of the ROM operator and use this to determine an efficient time integrator and time step for the ROM. \R{ETI}\hl{Since each $C_{r,i}$ is a skew-symmetric matrix, it has purely imaginary eigenvalues; similarly, since $D_{r}$ is a symmetric negative definite matrix, its eigenvalues are negative and real. Gershgorin's theorem can then be used to find a bound on the eigenvalues of the combined operator from which the time step can be determined \cite{Trias2011b}.} Although in this second step the non-linear stability property is strictly speaking lost, linear stability of the ROM in combination with adaptive high-order time stepping can result in a negligible energy error, which could form an efficient alternative to implicit time integration. Of course, this depends on the application under consideration. \hl{In the test cases in section \ref{sec:results}, several simulations with explicit Runge-Kutta methods (for both FOM and ROM) will be performed to show the potential of this alternative. The explicit Runge-Kutta method applied to the FOM reads \cite{Sanderse2012a}:
\begin{align}
\Omega_{h} \hat{V}_{h}^{n,i} &=  \Omega_{h} V_{h}^{n}  + \Delta t \sum_{j=1}^{i} \tilde{a}_{ij} (\hl{F^{CD}_{h}}(V_{h}^{n,j}) - G_{h} p_{h}^{n,j}), \label{eqn:mom_ERKstage} \\
L_{h} p_{h}^{n,i} &= M_{h} \hat{V}_{h}^{n,i} - y_{M}, \label{eqn:Poisosn_ERKstage} \\
\Omega_{h} V_{h}^{n,i} &= \Omega_{h} \hat{V}_{h}^{n,i} - G_{h} p_{h}^{n,i}.
\end{align}
Here $\tilde{a}$ denotes the shifted Butcher tableau, which includes the $b$ coefficients. For details, we refer to \cite{Sanderse2012a}. The application of an explicit Runge-Kutta method to the ROM is straightforward, as the pressure is absent.}

\subsection{\hl{Recovering the pressure}}\label{sec:pressure_recovery}
\hl{The pressure is not part of the ROM formulation \eqref{eqn:mom_ROM_ext}, because the reduced velocity field is by construction divergence-free. In principle one could recover the pressure once the velocity field is known (via $V_{h} \approx V_{r} = \Phi a$) by solving a Poisson equation at the FOM level, but then the computational costs scale with the FOM dimension. Instead, we aim for pressure recovery at the ROM level. Roughly speaking, two approaches exist \cite{Kean2019}: (i) solve a pressure Poisson equation (PPE) on the ROM level, or (ii) use the momentum equation recovery formulation with supremizer stabilization. In this work we follow the PPE strategy, as it naturally fits in the proposed framework.}

\hl{The ROM - PPE can be derived in two ways. Firstly, one can combine the divergence-free constraint and momentum equation on the ROM level: project the divergence equation based on a suitable pressure space, take its time derivative, and substitute the projected momentum equation. However, given that the velocity basis is divergence-free, this leads to a singular pressure equation, for the same reason that the pressure disappeared in equation \eqref{eqn:mom_ROM_ext}. The second approach, which we follow, is to combine the divergence-free constraint and momentum equation on the FOM level, to obtain the corresponding PPE on the FOM level (equation \eqref{eqn:PPE}), and then project it to the ROM level \cite{Akhtar2009}.}


\hl{In order to obtain a pressure estimate on the ROM level, we collect full-order snapshots of the pressure field in order to approximate the pressure in terms of 
\begin{equation}\label{eqn:pressure_expansion}
p_{h}(t) \approx  p_{r}(t) = \Pi q (t),
\end{equation}
where $\Pi \in \mathbb{R}^{N_{p} \times M_{p}}$ is constructed from a truncated SVD of snapshots of $p_{h}$ (similar to the construction of $\Phi$), such that $\Pi^{T} \Omega_{p} \Pi = I_{M_{p}}$, where $M_p$ is the number of pressure modes and $\Omega_{p}$ is a diagonal matrix containing the size of the pressure volumes. The velocity and pressure snapshots and modes are thus considered in a decoupled fashion \cite{Caiazzo2014}. It is possible to include the temporal mean of the snapshots of the pressure in equation \eqref{eqn:pressure_expansion}, like in \cite{Akhtar2009}, but we did not observe significant differences in the results in our test cases. Substituting the approximation for $p_{h}$ into equation \eqref{eqn:PPE} and projecting with $\Pi^{T}$ gives the PPE-ROM:
\begin{equation}\label{eqn:PPE_ROM}
L_{r} q (t) = \Pi^{T} M_{h} \Omega_{h}^{-1} F^{CD}_{h} (\Phi a(t)),
\end{equation}
where $L_{r} \in \mathbb{R}^{M_{p} \times M_{p}}$ is the reduced Poisson operator, given by
\begin{equation}
L_{r} = \Pi^{T} L_{h} \Pi = \Pi^{T}  M_{h} \Omega_{h}^{-1} G_{h} \Pi = - (G_{h} \Pi)^{T} \Omega_{h}^{-1} (G_{h} \Pi ),
\end{equation}
which is symmetric negative definite, like $L_{h}$, although it is not sparse. Like for the full PPE equation, explained in section \ref{sec:spatial}, no pressure boundary conditions are needed in \eqref{eqn:PPE_ROM}, as they were built into the divergence-free constraint. We stress again that this has been accomplished by deriving the PPE based on the spatially discrete mass and momentum equations. Compared to the derivation of the ROM-PPE based on projection at the continuous level, see e.g.\ \cite{Kean2019} (where boundary conditions need to be prescribed for the PPE), this is a strong advantage of our approach.}

\hl{In order to evaluate equation \eqref{eqn:PPE_ROM} using only ROM operators, a second offline decomposition is required, similar to what is proposed in \cite{Caiazzo2014}. This decomposition is constructed exactly the same as the one described in section \ref{sec:offline_decomposition}, with the main difference that instead of pre-multiplying with $\Phi^T$, we pre-multiply with $\Pi^T M_{h} \Omega_{h}^{-1} \in \mathbb{R}^{M_{p} \times N_{V}}$.}



%% file: results.tex
\section{Results}\label{sec:results}
In this section we show the results of three test cases. In the first test case, we demonstrate the stability and momentum and energy conservation properties of the ROM through an inviscid simulation of a shear-layer roll-up. In the \hl{second} test case we consider the simulation of a lid-driven cavity, a common test case used in the ROM community, for which several stabilization techniques have been tested. We will show that, independent of the number of modes, no stabilization method is needed in our approach, \hl{and that explicit time integration methods can be efficiently employed}. In the \hl{third} test case, we demonstrate the treatment of non-homogeneous boundary conditions (including outflow conditions), and \hl{unsteady forcing terms} by simulating the unsteady wake behind an actuator disk. \hl{We will show significant speed ups and accurate pressure recovery.}


\hl{In all test cases, the number of modes $M_{p}$ for the pressure basis is taken equal to the number of velocity modes $M$. Note that the pressure is only computed as a post-processing step at the end of each time step; it is not needed in the ROM formulation to advance the velocity coefficients. \R{errornorm}The natural error norm to measure the accuracy of the ROM with respect to the FOM is the weighted $L_{2}$-norm, see \eqref{eqn:innerproduct_semidiscrete}. This is the vector norm that corresponds to the Frobenius norm being minimized in equation \eqref{eqn:minimization_problem}. This norm is scaled with the norm of a characteristic velocity $V_{\text{ref}}$. The error in the velocity at each time step is thus computed as
\begin{equation}\label{eqn:error_velocity_L2}
\epsilon_{V}^n = \frac{\|V_{r}^{n} - V_{h}^{n} \|_{\Omega_{h}}}{\| V_{\text{ref}} \|_{\Omega_{h}}}.
\end{equation}
We note that the error between ROM and FOM is not only due to projection when different time integration methods are used on the ROM and FOM level. In that case, the error can be decomposed as:
\begin{equation}\label{eqn:velocity_error_decomposition}
V_{r}^{n} - V_{h}^{n}  = \underbrace{V_{r}^{n} - V_{r}(t^{n})}_{\text{time integration error ROM}} + \underbrace{V_{r}(t^{n}) - V_{h}(t^{n})}_{\text{projection error}} + \underbrace{V_{h}(t^{n})- V_{h}^{n}}_{\text{time integration error FOM}}, 
\end{equation}
which illustrates that the projection is performed at the semi-discrete level (see figure \ref{fig:discretise_projection}). Note that in the second and third test case in this paper, these errors are absent, as the ROM and FOM use the same time integration method, so that we effectively project the fully discrete equations. In the first test case we study the error in momentum and energy, equations \eqref{eqn:discrete_mom_eqn_ROM} and \eqref{eqn:discrete_energy_eqn_ROM}.}

\hl{Similar to \eqref{eqn:error_velocity_L2}, the ROM pressure error is defined as
\begin{equation}\label{eqn:error_pressure_L2}
\epsilon_{p}^n =\frac{\|p_{r}^{n} - p_{h}^{n} \|_{\Omega_{p}}}{\| p_{\text{ref}} \|_{\Omega_{p}} }.
\end{equation}
Since the pressure in the incompressible Navier-Stokes equations is determined up to a constant, both the ROM and FOM pressure fields are shifted to have the same spatial mean. We will also report the basis projection error in the second and third test case:
\begin{equation}\label{eqn:basis_projection_error}
\epsilon_{V,\text{best}}^{n} =  \|V_{\text{best}}^{n} - V_{h}^{n} \|_{\Omega_{h}},
\end{equation}
where $V_{\text{best}}$ is obtained by projecting the FOM solution onto the ROM basis, i.e.\ $V_{\text{best}}^{n} = \Phi \Phi^{T} \Omega_{h} V_{h}^{n}$.}

\hl{An important aspect of the ROM is the computational speed-up compared to the FOM. We will split the ROM time in three components: (1) (offline) computation of the basis (via SVD); (2) (offline) precomputing the ROM operators, (3) (online) time-stepping of the ROM. The time needed to compute the error of the ROM with respect to the FOM will be excluded from the online timings, as it can be computed as a postprocessing step (as long as the time-dependent ROM coefficients are stored while marching in time). The basis construction via SVD could be accelerated by using the method of snapshots, but for the relatively small (two-dimensional) problems considered here, the computational gain is not significant.}

\subsection{Shear-layer roll-up}\label{sec:shear_layer}
\begin{figure}[hbtp]
\centering
	\includegraphics[width=0.5\textwidth]{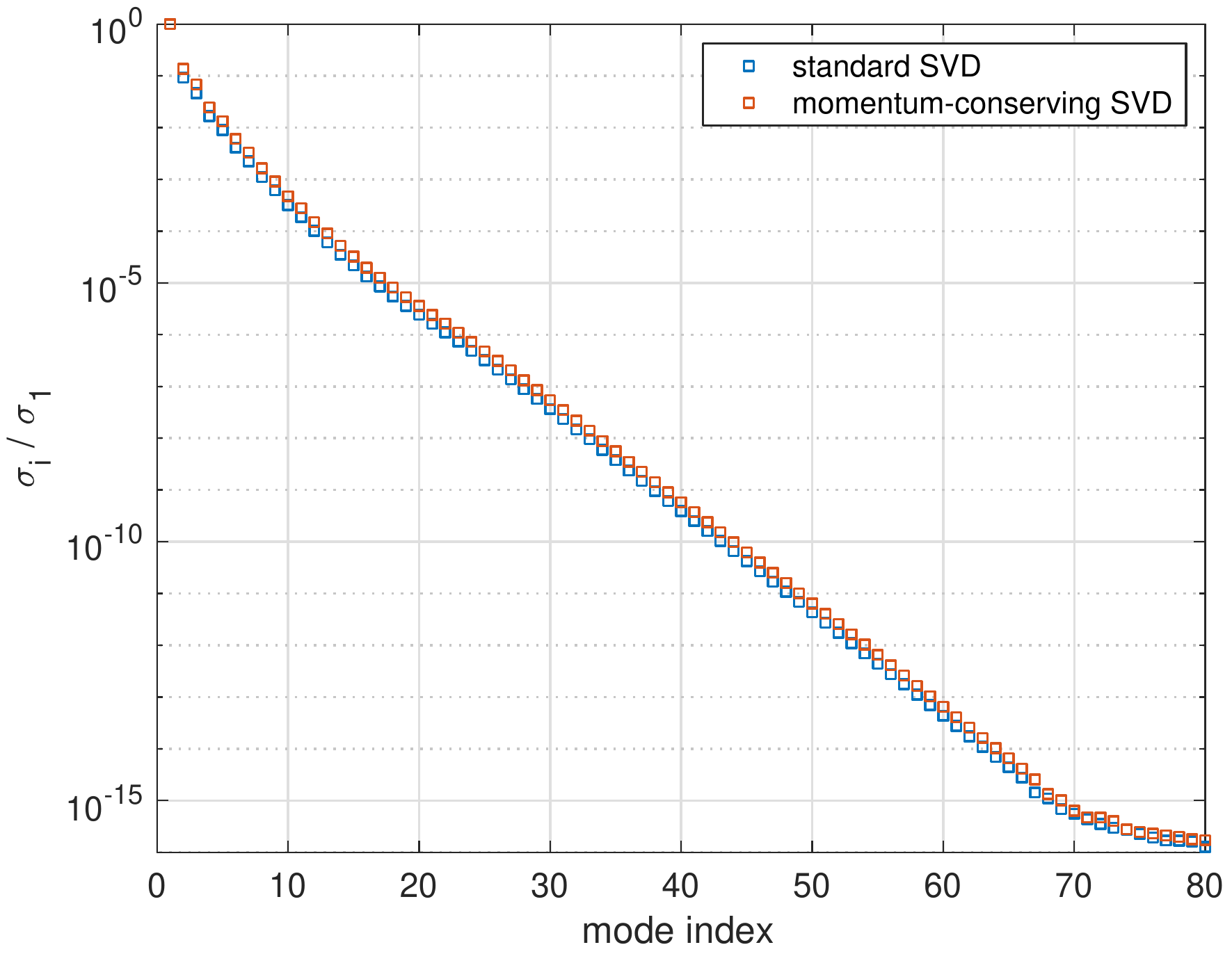}
	\caption{Singular values for inviscid shear-layer roll-up.\label{fig:singular_values_inviscid_shearlayer}}
\end{figure}

We simulate the roll-up of a shear-layer, similar to \cite{Sanderse2013}. The simulation domain is $[0,2\pi] \times [0,2\pi]$, with periodic boundary conditions and the following initial condition:
\begin{equation}
u_{0}(x,y) = 1 + 
\begin{cases}
\tanh ( \frac{y-\pi/2}{\delta}), & y\leq \pi, \\
\tanh ( \frac{3\pi/2 -y}{\delta}), & y>\pi, 
\end{cases}
\qquad 
v_{0}(x,y) = \epsilon \sin (x),
\end{equation}
where $\delta = \pi/15$ and $\epsilon = 1/20$. Compared to \cite{Sanderse2013}, a constant has been added to $u_{0}(x,y)$, in order to ensure that the global momentum of the $u$- and $v$- components differ. In the inviscid case, the energy of the flow should be exactly conserved. The FOM discretisation consists of $200 \times 200$ finite volumes, giving a total of $N_{V} + N_{p} = 1.2 \cdot 10^{5}$ unknowns. Time integration of the FOM is performed with explicit RK4 \cite{Sanderse2012a} with a time step of $\Delta t = 0.01$ from $t=0$ to $t=4$, \hl{resulting in $K=401$ snapshots} (as mentioned in section \ref{sec:alternative_time_integration}, the FOM snapshots need not be energy-conserving). The singular values of the velocity snapshot matrix are shown in figure \ref{fig:singular_values_inviscid_shearlayer}. The ROM basis consists of the first $M$ left singular vectors of the snapshot matrix, where we take $M =2, 4, 8, 16$. The rapid decay in the singular values indicates that the problem is suited for dimension reduction. The effect of using the proposed momentum-conserving SVD of section \ref{sec:momentum_conservation} instead of the standard SVD is a small shift in the singular values.

\begin{figure}[hbtp]
\centering
	\begin{subfigure}[b]{.49\textwidth}
	\centering
	\includegraphics[width=\textwidth]{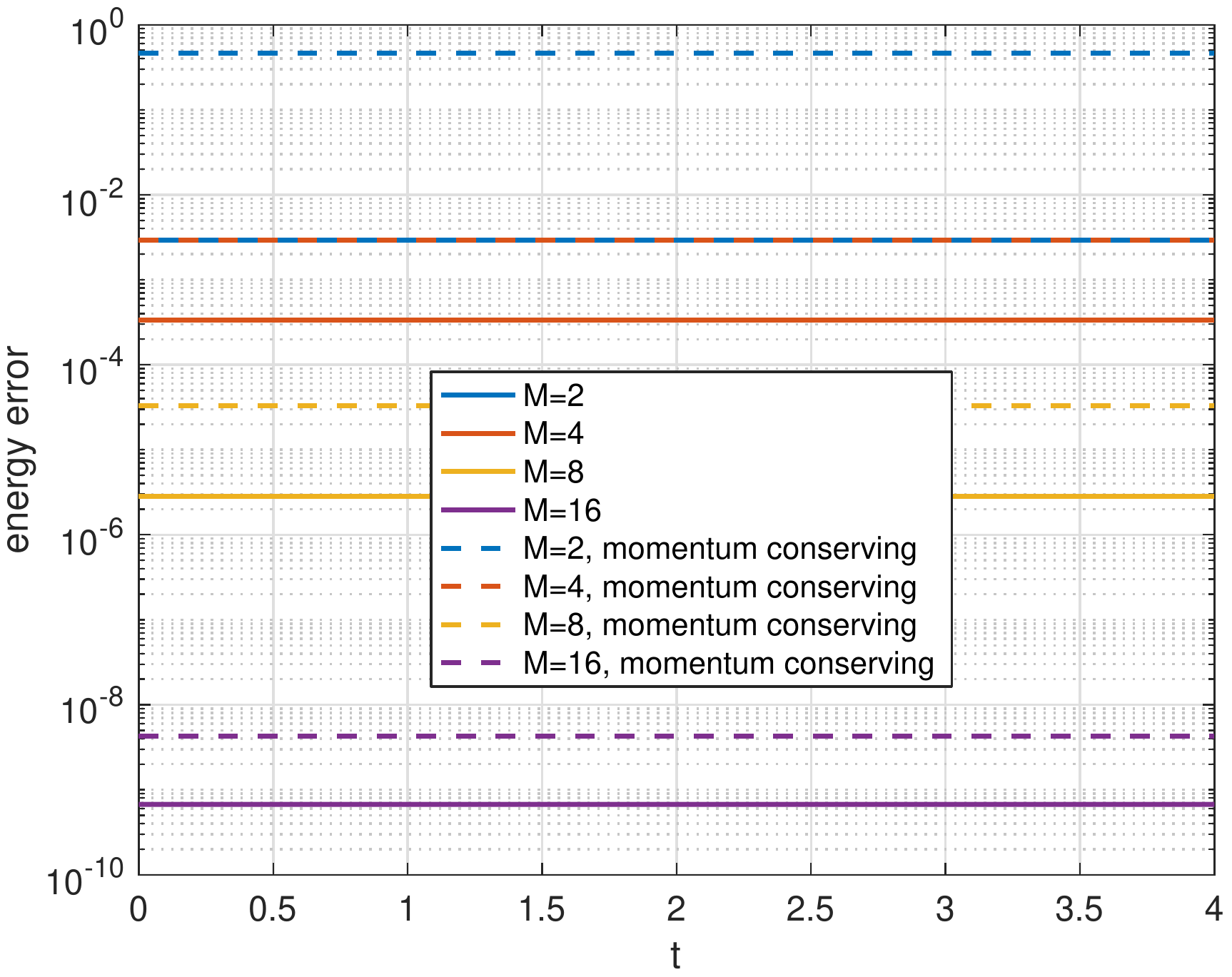}
	\caption{Energy error: $\frac{K_{r}^{n}-K_{h}(0)}{K_{h}(0)}$.\label{fig:energy_error_inviscid_shearlayer}}
	\end{subfigure}
	\hfill
	\begin{subfigure}[b]{.49\textwidth}
	\centering
	\includegraphics[width=\textwidth]{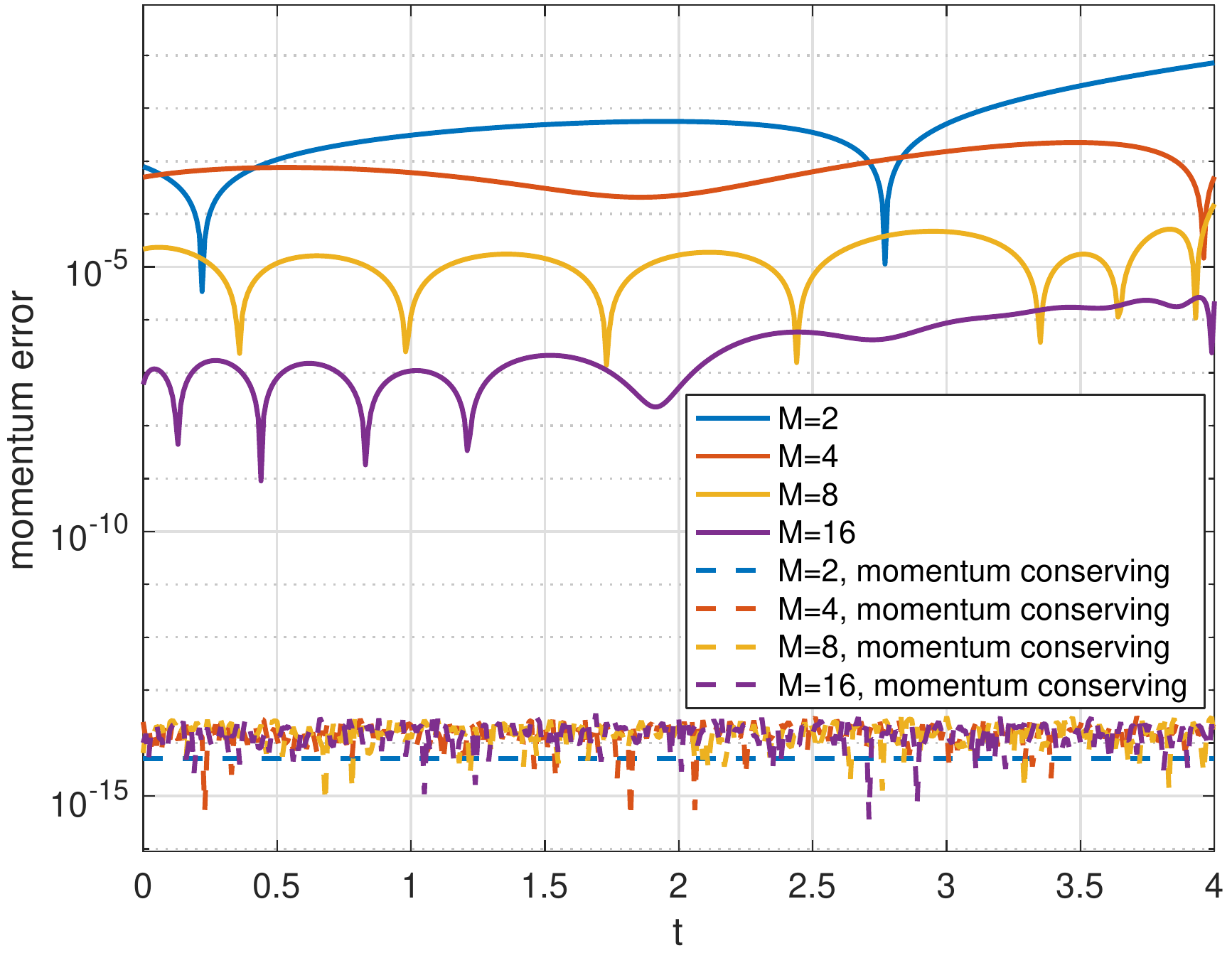}
	\caption{Momentum error: $\frac{P_{r}^{u,n}-P_{h}^{u}(0)}{P^{u}_{h}(0)}$.\label{fig:mom_error_inviscid_shearlayer}}
	\end{subfigure}
\caption{Energy and momentum conservation of ROM with respect to FOM initial condition.}
\end{figure}

\begin{figure}[hbtp]
\centering
		\includegraphics[width=0.5\textwidth]{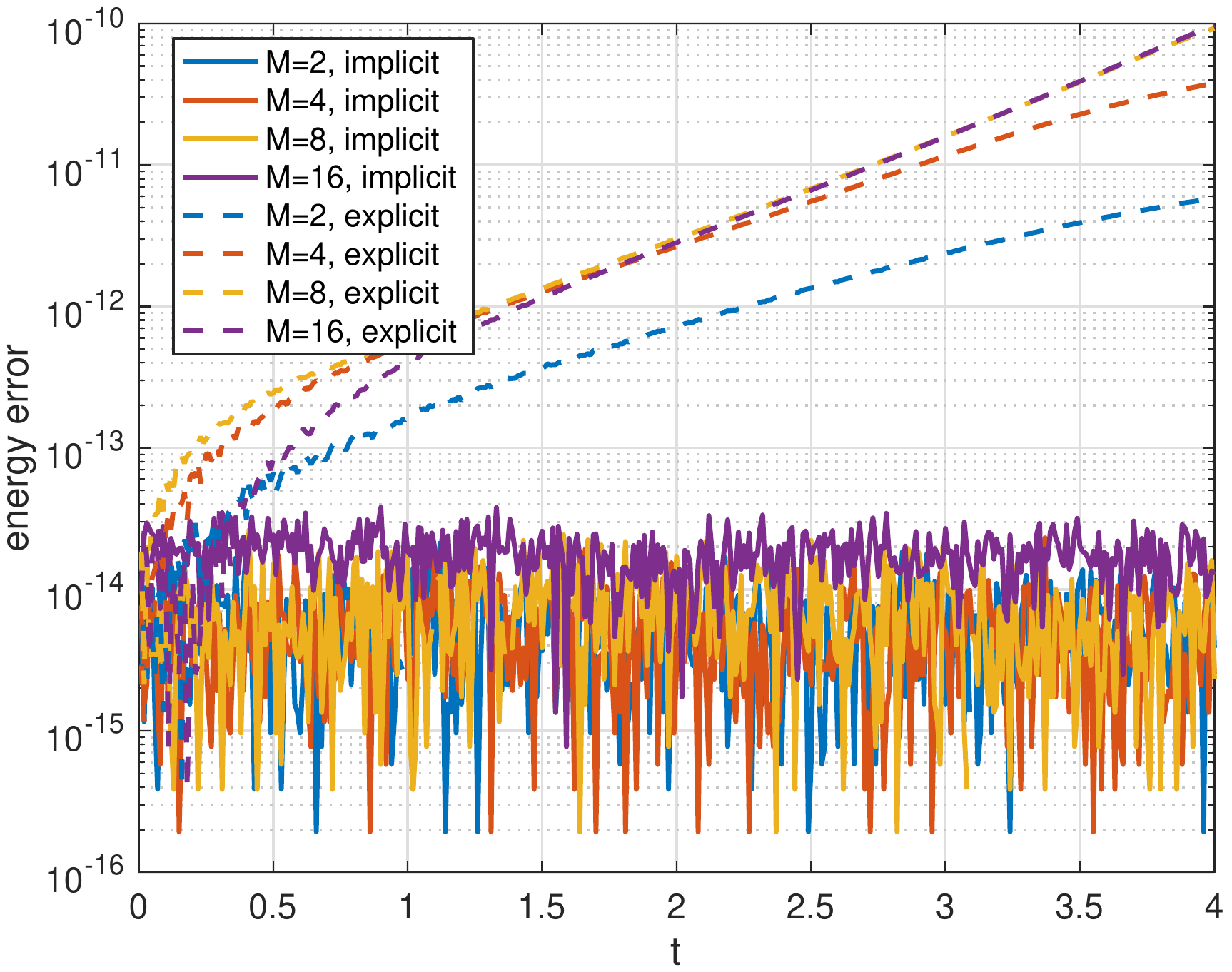}
		\caption{Energy error of the ROM with respect to ROM initial condition, $\frac{K_{r}^{n} - K_{r}(0)}{K_{r}(0)}$, for implicit midpoint and explicit RK4. \label{fig:error_shearlayer_unsteady_IMvsRK}}
\end{figure}

Time integration of the ROM is performed with the implicit midpoint method, with the same $\Delta t$ and end time as used for the FOM. \hl{The energy error between fully discrete ROM and FOM can be decomposed according to equation \eqref{eqn:energy_error_decomposition}:
\begin{equation}\label{eqn:energy_error_decomposition2}
\epsilon_{K} := K_{r}^{n} - K_{h}^{n}  = \underbrace{K_{r}^{n} - K_{r}(0)}_{\text{time integration error ROM}} + \underbrace{K_{r}(0) - K_{h}(0)}_{\text{projection error}} + \underbrace{K_{h}(0)- K_{h}^{n}}_{\text{time integration error FOM}}. 
\end{equation}
We are interested in the first two terms, whose sum is shown in figure \ref{fig:energy_error_inviscid_shearlayer} (the last error term is nonzero, as the FOM time integration is explicit RK4). The first term remains at machine-precision zero, as is further detailed in figure \ref{fig:error_shearlayer_unsteady_IMvsRK}. The ROM is therefore conserving kinetic energy: $K_{r}^{n}=K_{r}(0)$, independent of whether the standard SVD or the momentum-conserving SVD is used. The error displayed in the figure thus corresponds to the second term: $K_{r}(0)-K_{h}(0)$, which is due to projecting the initial FOM velocity field onto the truncated snapshot basis, and decreases with increasing $M$ according to the singular value decay.} We observe that, especially for small $M$, the momentum conserving approach is less accurate in terms of the energy error. This is because two modes have been sacrificed in order to achieve momentum conservation.  For the case $M=2$ this means that momentum is enforced, but that the FOM snapshots are not taken into account in the basis $\Phi$.



In figure \ref{fig:mom_error_inviscid_shearlayer} the momentum error of the $u$-component is plotted as a function of time. For the standard SVD, the error in global momentum increases as a function of time, and decreases when more modes are taken. With the momentum-conserving SVD, the error in momentum stays at machine precision, independent of the number of modes. As noted in section \ref{sec:ECROM}, addition of the initial velocity field as constraint in the SVD can also force the kinetic energy error to zero, independent of $M$. This is not considered here, as it would not generalize to the case of viscous flows.

\begin{figure}[hbtp]
\centering
		\includegraphics[width=0.5\textwidth]{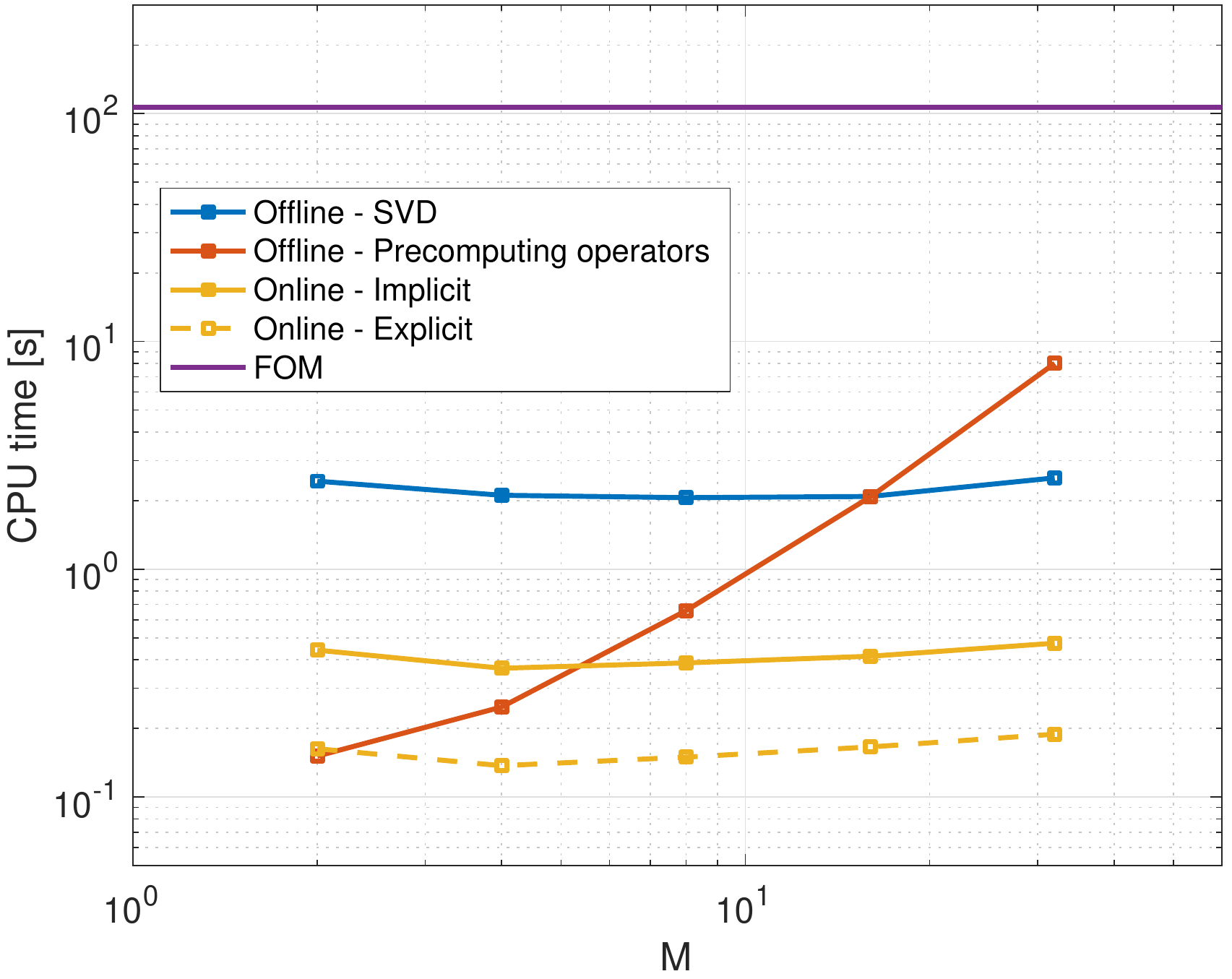}
		\caption{Computational time as a function of number of modes in ROM for inviscid shear-layer roll-up. \label{fig:CPUtimings_shearlayer_unsteady}}
\end{figure}

\hl{Figure \ref{fig:error_shearlayer_unsteady_IMvsRK} further details the energy error of the ROM when excluding the projection error, i.e.\ we plot the first term of equation \eqref{eqn:energy_error_decomposition2}: $K_{r}^{n} - K_{r}(0)$. For implicit midpoint this error remains at machine precision, as expected. When changing the ROM time integration to an explicit (fourth-order Runge-Kutta) method, an energy error is introduced which grows steadily in time. However, given that this error is insignificant compared to the projection error shown in figure \ref{fig:energy_error_inviscid_shearlayer}, RK4 can be a practical alternative given that it is faster to execute on a per-time step basis than the implicit midpoint method (provided that the time-step restriction on the explicit method is not too severe). This is shown in figure \ref{fig:CPUtimings_shearlayer_unsteady}, which compares the computational cost of the ROM to that of the FOM, for both implicit and explicit time integration of the ROM. For $M\leq 8$ we have a computational speed-up factor of at least 50, if we include the offline phase in the speed-up computation. When not considering the offline phase, the speed-up factor is much larger: around 400 for implicit midpoint, and almost 1000 for explicit RK4. 
We note that $M\leq 8$ is realistic from the viewpoint of accuracy: for $M=8$ we have a relative energy error which is less than $10^{-5}$ (figure \ref{fig:energy_error_inviscid_shearlayer}). For $M>10$, the computational costs of the ROM are dominated by the offline phase (precomputing the operators, especially the convective operator); the SVD and the online phase are negligible compared to the precomputing phase.
Given the small energy error with explicit RK4, and the fact that it is significantly faster to execute, the next test cases will be performed with explicit RK4 for both the FOM and ROM.}




\FloatBarrier

\subsection{Lid-driven cavity}\label{sec:LDC}
We perform as second test case a common test used to assess ROMs for the incompressible Navier-Stokes equations (see e.g.\ \cite{Fick2018,Stabile2018}): a lid-driven cavity flow at $\text{Re}=1,000$. The velocity of the lid enters as boundary contribution in \hl{$y_{D}$}, but does not appear in \hl{$y_{M}$ in the divergence-free constraint}, since $y_{M}$ only contains velocity components normal to the boundary. Consequently, in this test case the procedures described in section \ref{sec:BC} are not required ($V_{bc}=0$). No measures need to be taken to ensure stability of the ROM, as it is stable by design.

The full-order model is run on a grid with $100 \times 100$ volumes, resulting in $N_V + N_p = 3\cdot 10^{4}$ unknowns, at a time step of $\Delta t = 0.01$, until a final time $T = 10$. This results in $K=1,001$ snapshots. \hl{Note that, as discussed in sections \ref{sec:alternative_time_integration} and \ref{sec:shear_layer}, both FOM and ROM are now integrated with an explicit RK4 method, in order to show that in practice explicit methods can be an excellent alternative to the implicit methods proposed in section \ref{sec:ECROM}. In contrast to e.g.\ \cite{Stabile2018}, stable and accurate results are obtained without requiring a stabilization method.}

\begin{figure}[hbtp]
\centering
	\begin{subfigure}[b]{.49\textwidth}
		\centering
		\includegraphics[width=\textwidth]{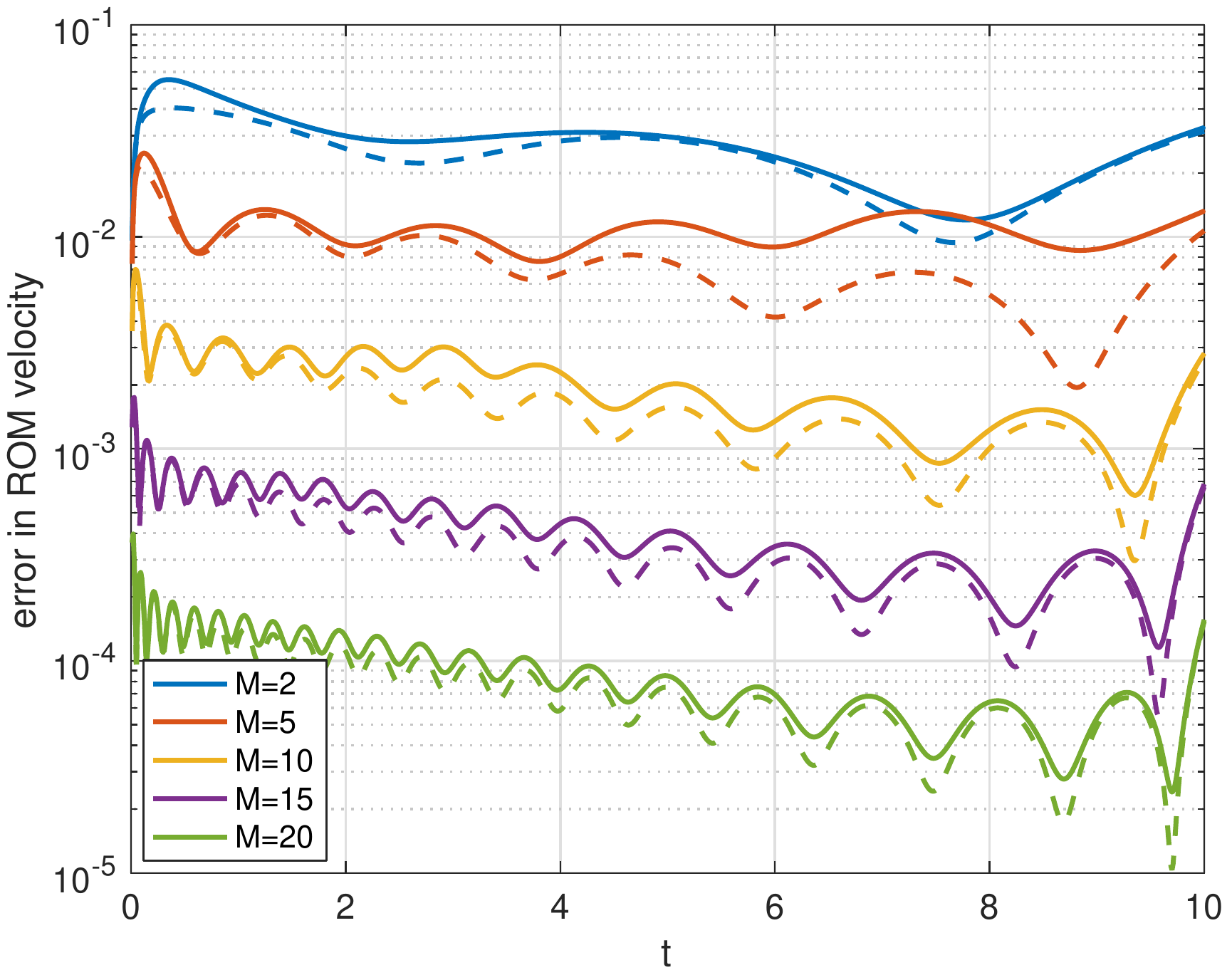}
	\caption{Velocity: $\epsilon_{V}^n$. \label{fig:error_ROM_Re1000_modes}}
	\end{subfigure}
	\hfill
	\begin{subfigure}[b]{.49 \textwidth}
	\centering
			\includegraphics[width=\textwidth]{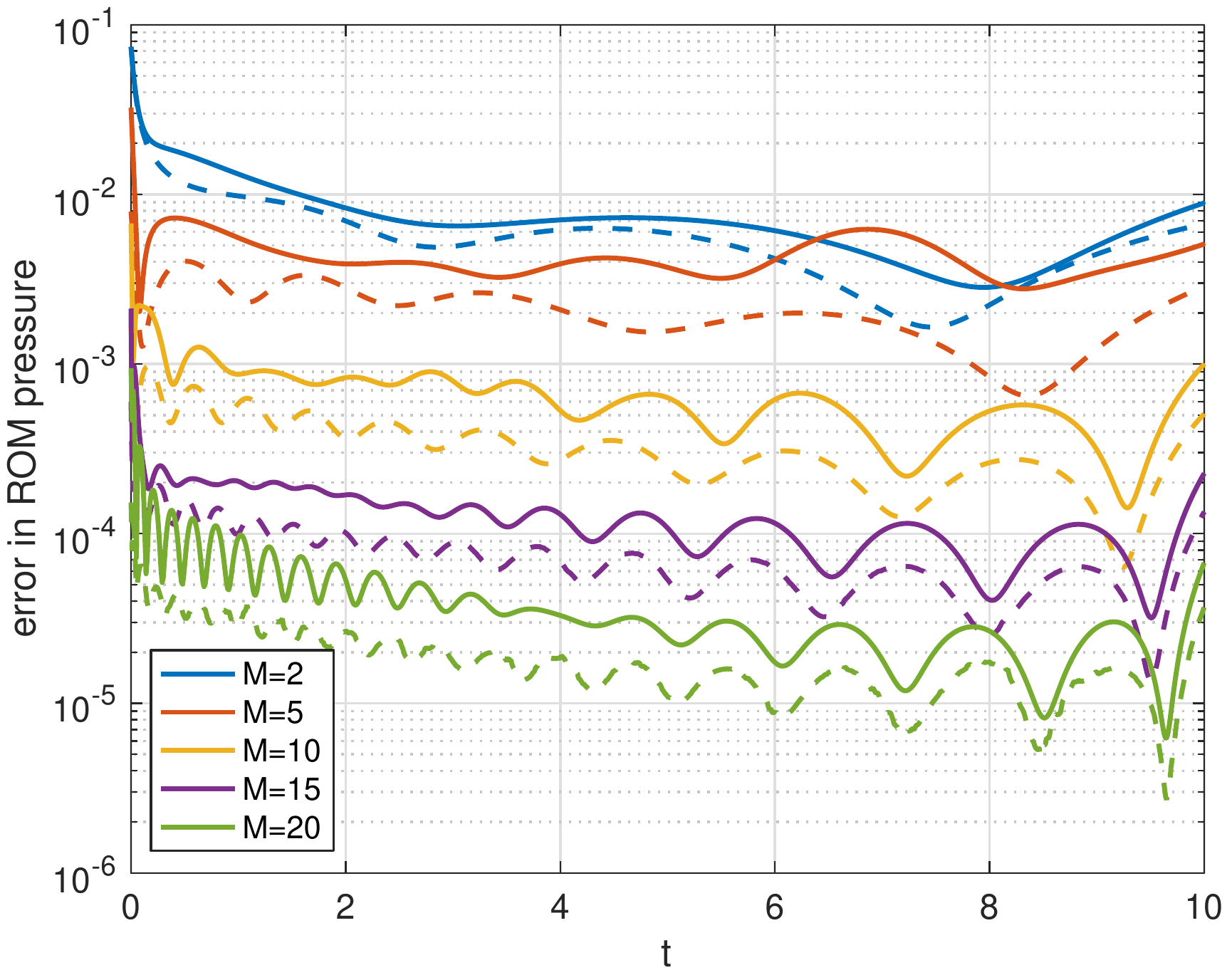}
		\caption{Pressure: $\epsilon_{p}^n$.\label{fig:error_pressure_Re1000_modes}}
	\end{subfigure}
	\caption{Error in ROM with respect to the FOM as a function of time for different number of modes for the lid-driven cavity flow. Dashed lines: basis projection error (projecting snapshots onto truncated basis).}
\end{figure}

\hl{The error with respect to the FOM is given by equation \eqref{eqn:error_velocity_L2} where $V_{\text{ref}}=1$ (the lid velocity).} We study the effect of increasing the number of modes on the accuracy of the velocity field, while using the full snapshot set as basis for the SVD; see figure \ref{fig:error_ROM_Re1000_modes}. We clearly see how the accuracy increases when increasing the number of modes. \hl{Furthermore, the error in the ROM is very close to the `best possible' error, being the basis projection error, equation \eqref{eqn:basis_projection_error}. The behaviour for the accuracy of the pressure, using the pressure recovery method described in section \ref{sec:pressure_recovery}, is very similar, as displayed in figure \ref{fig:error_pressure_Re1000_modes}. The error is computed from \eqref{eqn:error_pressure_L2} with $p_{\text{ref}}=1$, which is chosen according to the typical pressure differential between the upper-right and upper-left corner of the cavity.}

Note that the kinetic energy is not conserved in this test case, but rather increases as a function of time, as kinetic energy is added to the flow through the moving lid, which is initially larger than the energy dissipation in the interior of the cavity. When the flow reaches a steady state, the two effects will balance each other. 

\begin{figure}[hbtp]
\centering
		\includegraphics[width=0.5\textwidth]{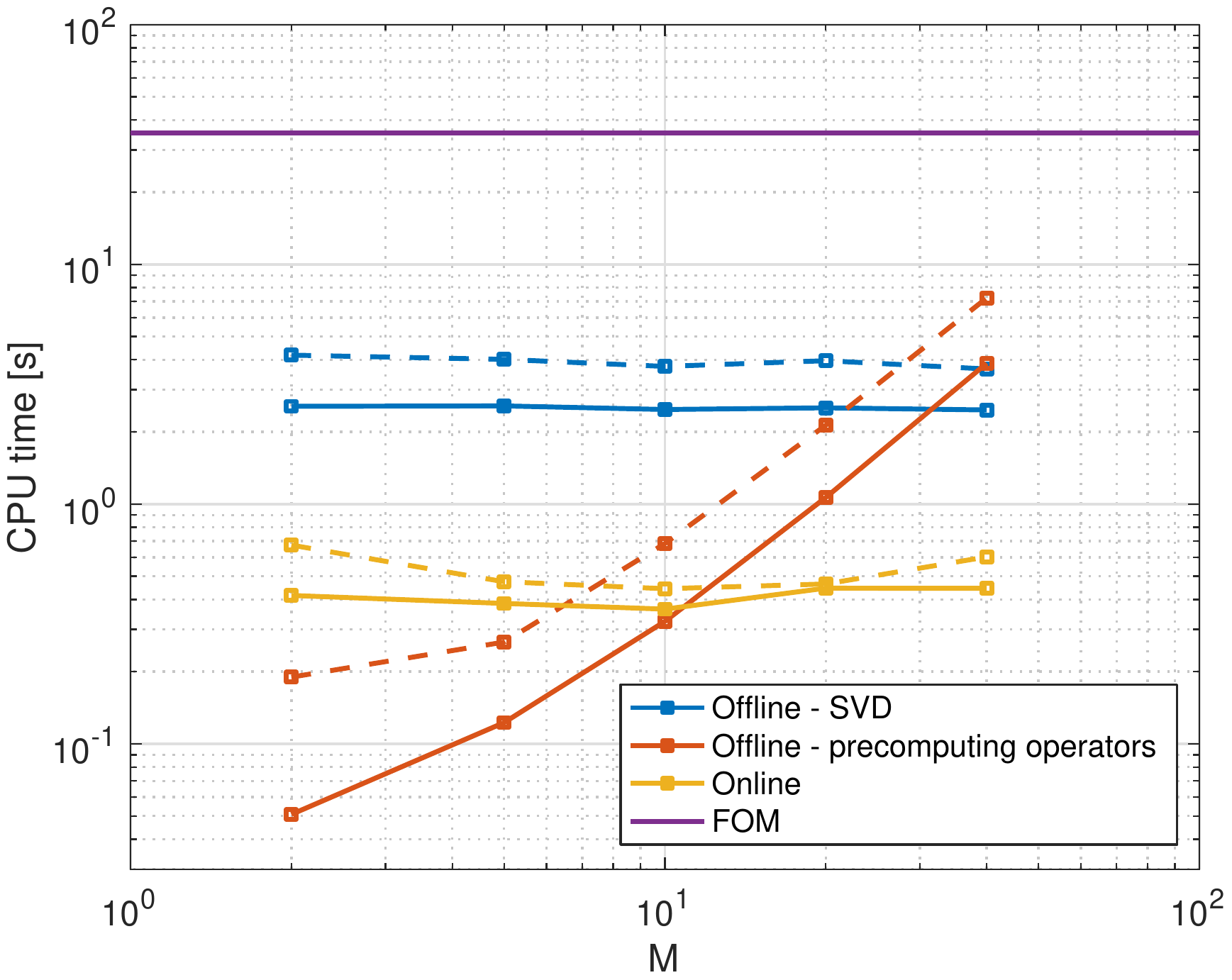}
		\caption{Computational time as a function of number of modes in ROM. Dashed lines: including pressure recovery. \label{fig:CPUtimings_LDC_unsteady}}
\end{figure}

\hl{In figure \ref{fig:CPUtimings_LDC_unsteady} the computational time required by the ROM is compared to the FOM. Like in the previous test case, the speed-up of the ROM compared to the FOM is significant: the online phase is about 100 times faster than the FOM. However, for a large number of modes ($M\geq 40$), precomputing the operators becomes significant and is dominating the ROM cost. In practice, such a large number of modes is not needed, since for $M=15$ the error in velocity and pressure is already smaller than $10^{-3}$ for almost the entire simulation. Figure \ref{fig:CPUtimings_LDC_unsteady} also reveals that the computational costs associated with pressure recovery is significant (dashed lines): it adds roughly a factor of 2 to the offline stage, because it requires an additional SVD and an additional set of precomputed operators.}



\FloatBarrier

\subsection{Actuator in non-uniform inflow}
In this test case we consider an actuator disk in a non-uniform flow field. The actuator disk \hl{is added as an unsteady sink term in the momentum equations} and is typically used to model the flow through wind turbines \cite{Sanderse2011}. This test case features non-homogeneous boundary conditions and we therefore require the method proposed in section \ref{sec:BC}. 

The test case set up is as follows. We consider a simulation domain \hl{$[-4,8]\times [-2,2]$} with the following inflow conditions at $x=-4$:
\begin{equation}
u(x=-4,y) = \frac{3}{4}- \frac{3}{32} (y-2)(y+2).
\end{equation}
This is a parabolic velocity profile with a mean equal to 1. At the domain boundaries \hl{$x=8$}, $y=-2$, and $y=2$ we employ outflow conditions at which the total stress is prescribed, \hl{according to equation \eqref{eqn:outflow_BC}}. For example, at $x=8$ we have, taking $p_{\infty}=0$:
\begin{equation}
-p + \nu \dd{u}{x} = 0, \qquad \dd{v}{x}=0.
\end{equation} 
The initial condition is the parabolic velocity profile and the Reynolds number is \hl{500}. \hl{The actuator is an infinitely thin disk located at $x=0$, $-0.5\leq y\leq 0.5$, with a thrust coefficient of $C_{T}=\frac{1}{2}$, inducing a discontinuity in the pressure (for more details, see \cite{Sanderse2011a}). The additional term in the FOM $u$-momentum equation is
\begin{equation}
f_{h}^{u} (t) = - C_{T} \Delta y (1 + \sin( \pi t)), 
\end{equation}
for the finite volumes intersected by the actuator. The first part of this expression is time-independent and can be precomputed. The second part $(1 + \sin( \pi t))$ is simple to evaluate in online mode and is multiplied with the precomputed forcing term.}

\begin{figure}[hbtp]
\centering
	\begin{subfigure}[t]{.49\textwidth}
		\centering
		\includegraphics[width=\textwidth]{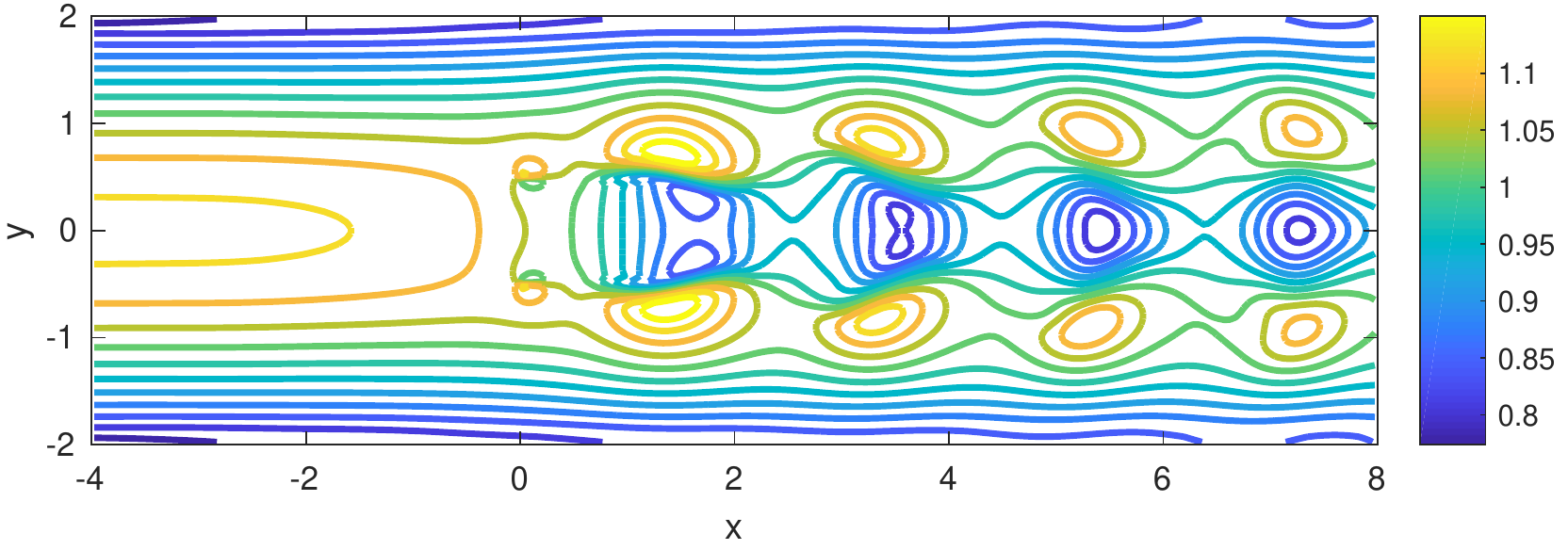}
		\caption{Full-order model, 57,840  degrees of freedom. \label{fig:velocity_actuator_FOM}}
	\end{subfigure}
	\hfill	
	\begin{subfigure}[t]{.49\textwidth}
		\centering		
		\includegraphics[width=\textwidth]{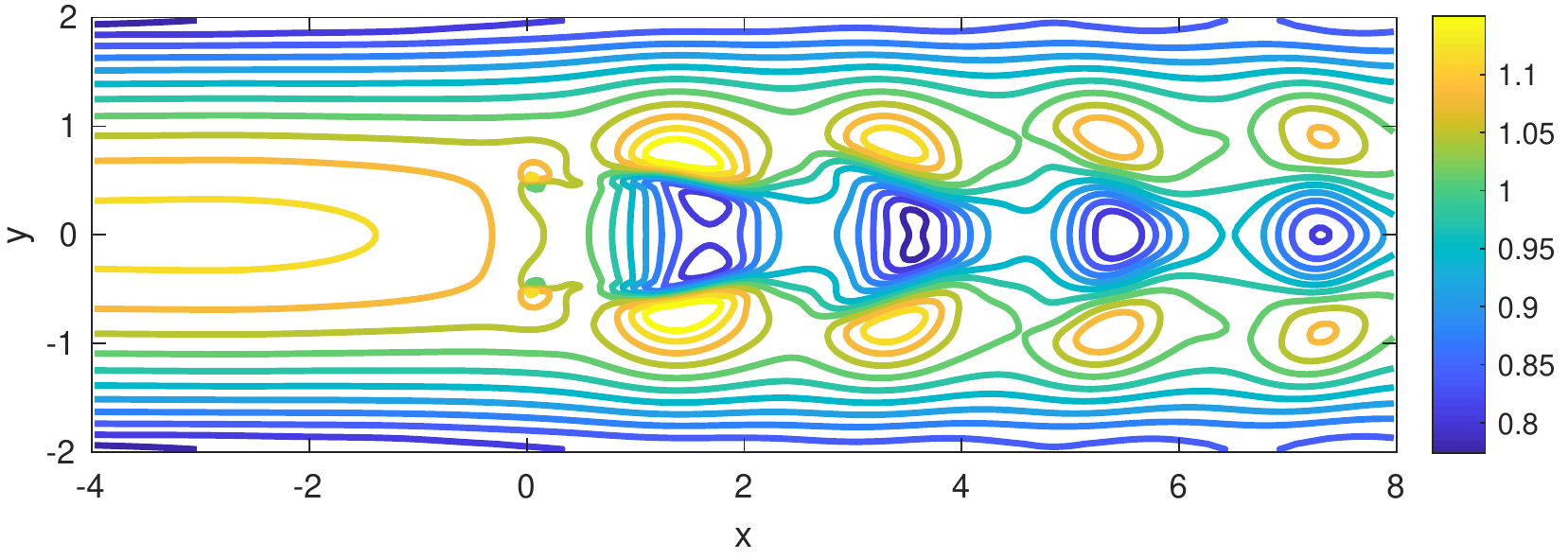}
		\caption{Reduced-order model, 10 degrees of freedom. \label{fig:velocity_actuator_ROM}}
			\end{subfigure}
	\caption{Contours of velocity magnitude for actuator disk test case, at $t=20$. \label{fig:velocity_actuator}}
\end{figure}

\begin{figure}[hbtp]
\centering
	\begin{subfigure}[t]{.49\textwidth}
		\centering
		\includegraphics[width=\textwidth]{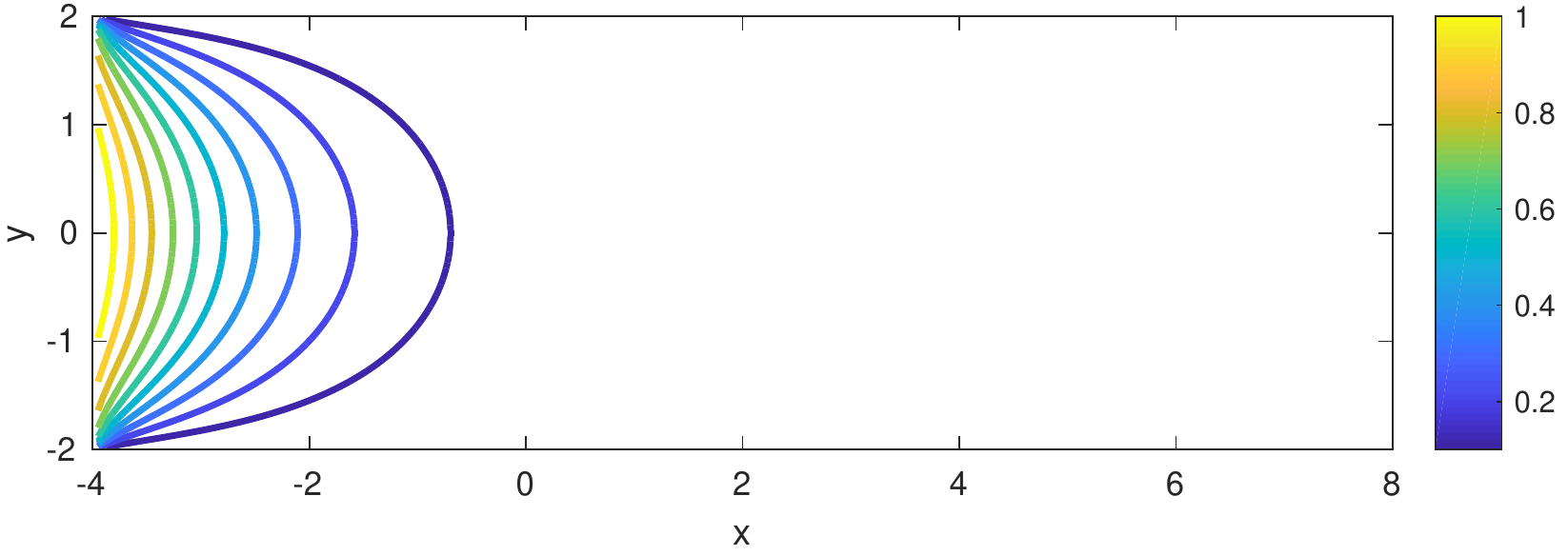}
		\caption{$u_{bc}$ \label{fig:actuator_u_bc}}
	\end{subfigure}
	\hfill	
	\begin{subfigure}[t]{.49\textwidth}
		\centering		
		\includegraphics[width=\textwidth]{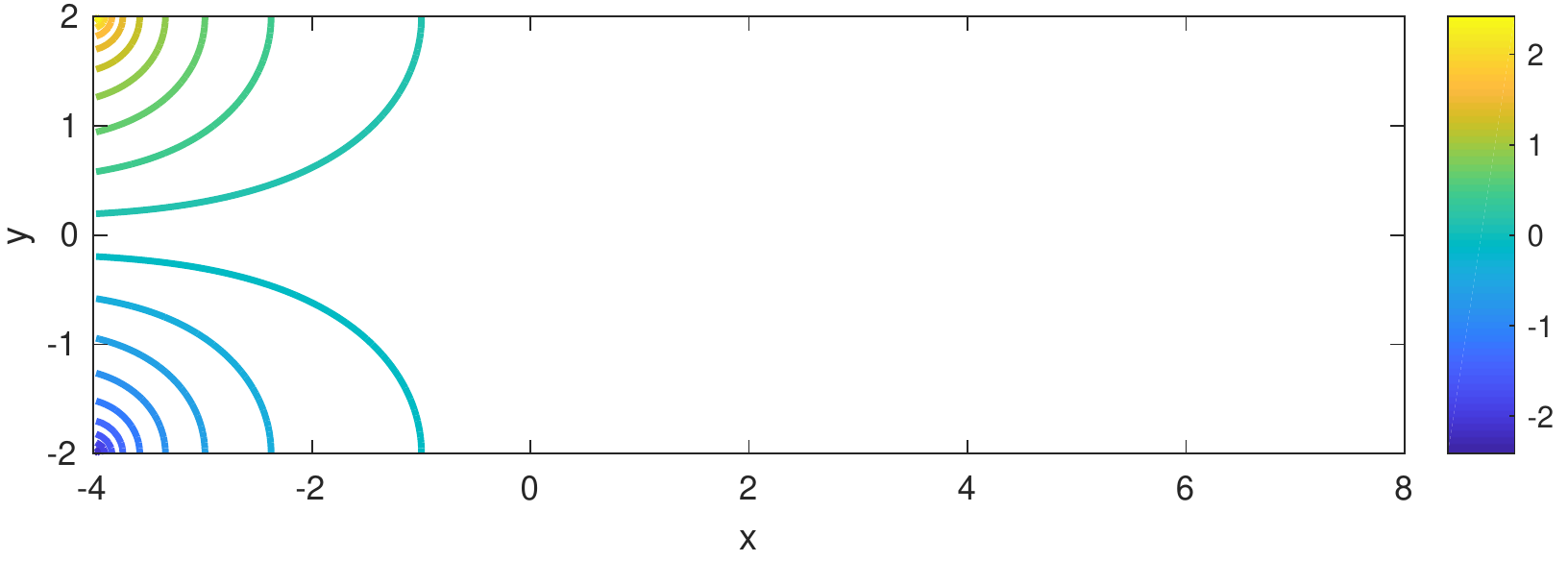}
		\caption{$v_{bc}$ \label{fig:actuator_v_bc}}
	\end{subfigure}
	\caption{Components of boundary condition function $V_{bc}$ for actuator test case.\label{fig:bc_actuator}}
\end{figure}

\begin{figure}[hbtp]
\centering
	\begin{subfigure}[t]{.49\textwidth}
		\centering
		\includegraphics[width=\textwidth]{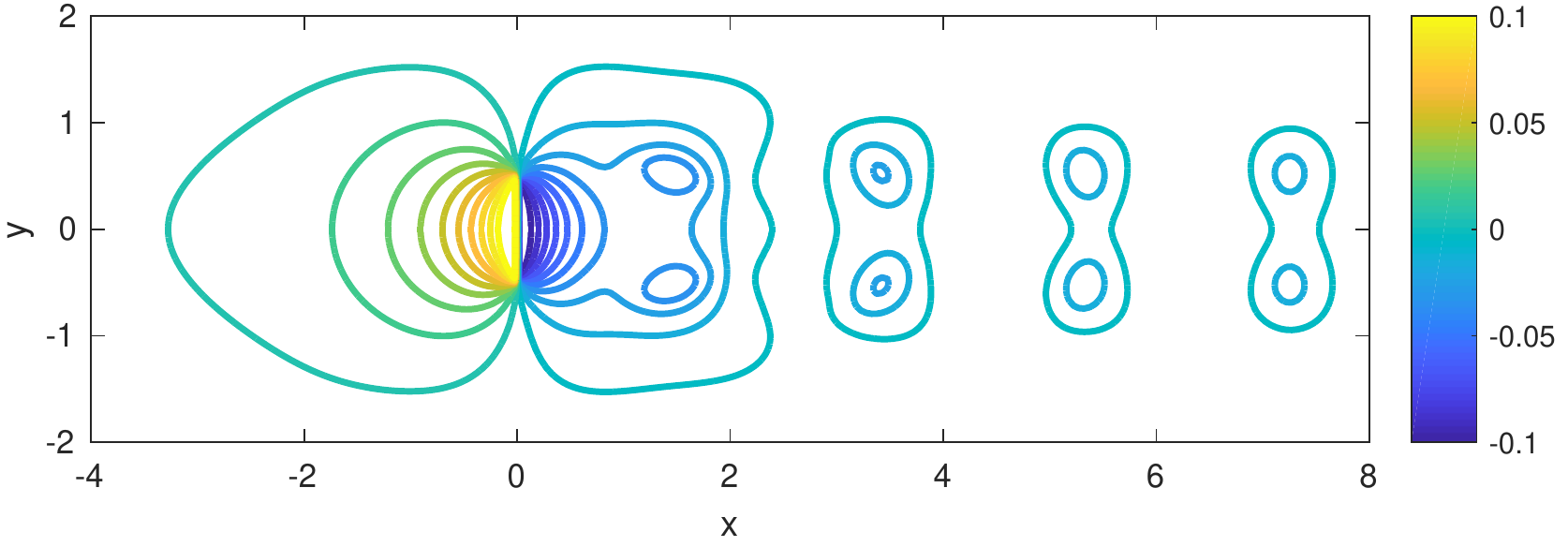}
		\caption{Full-order model, 57,840 degrees of freedom. \label{fig:pressure_actuator_FOM}}
	\end{subfigure}
	\hfill	
	\begin{subfigure}[t]{.49\textwidth}
		\centering		
		\includegraphics[width=\textwidth]{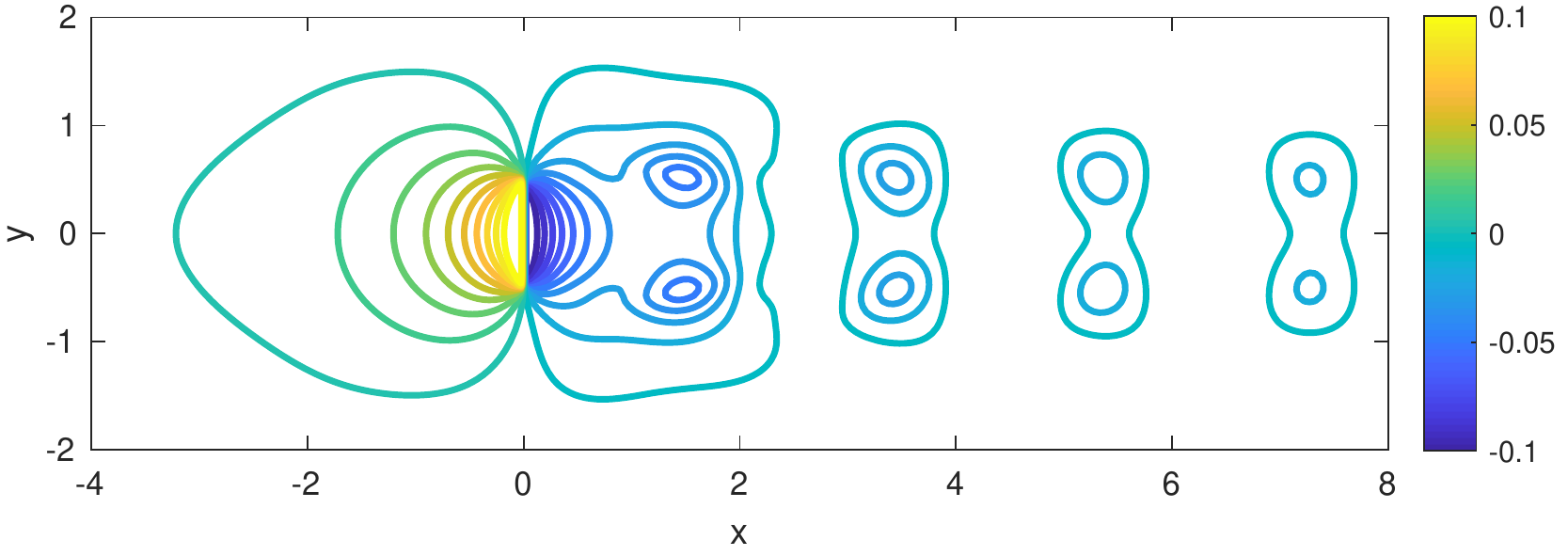}
		\caption{Reduced-order model, 10 degrees of freedom. \label{fig:pressure_actuator_ROM}}
			\end{subfigure}
	\caption{Contours of pressure for actuator disk test case, at $t=20$. \label{fig:pressure_actuator}}
\end{figure}

We first simulate the FOM with \hl{$240\times 80$} finite volumes from $t=0$ to \hl{$t=20$} with $\Delta t = 0.025$ and a classic RK4 scheme, giving $K=801$ snapshots. The resulting velocity field at \hl{$t=20$} is shown in figure \ref{fig:velocity_actuator_FOM}. Based on the FOM, the velocity field due to nonhomogeneous boundary conditions, $V_{bc}$, is computed from equation \eqref{eqn:poisson_Vbc}. Note that $V_{bc}$ is a vector field defined throughout the domain, and not only on the boundary. The components of $V_{bc}$ are shown in figure \ref{fig:bc_actuator}. The $V_{bc}$ field is subtracted from the snapshot matrix. We then simulate the ROM with \hl{$M=10$} modes, and the same time integration method and time step. Figure \ref{fig:velocity_actuator_ROM} shows the ROM velocity field, which is almost identical to the one obtained by the FOM. \hl{Given the velocity field, the pressure can be recovered using the technique from section \ref{sec:pressure_recovery}. This leads to an accurate pressure field, as shown in figure \ref{fig:pressure_actuator}. Similar to the closed flow (lid-driven cavity) case, also for this open flow case we do not need to specify any additional pressure boundary conditions at the ROM level. The boundary conditions were specified on the discrete FOM level, and inherited by the ROM via the projection of the boundary vectors.}

\begin{figure}[hbtp]
\centering
	\begin{subfigure}[t]{.49\textwidth}
		\centering
		\includegraphics[width=\textwidth]{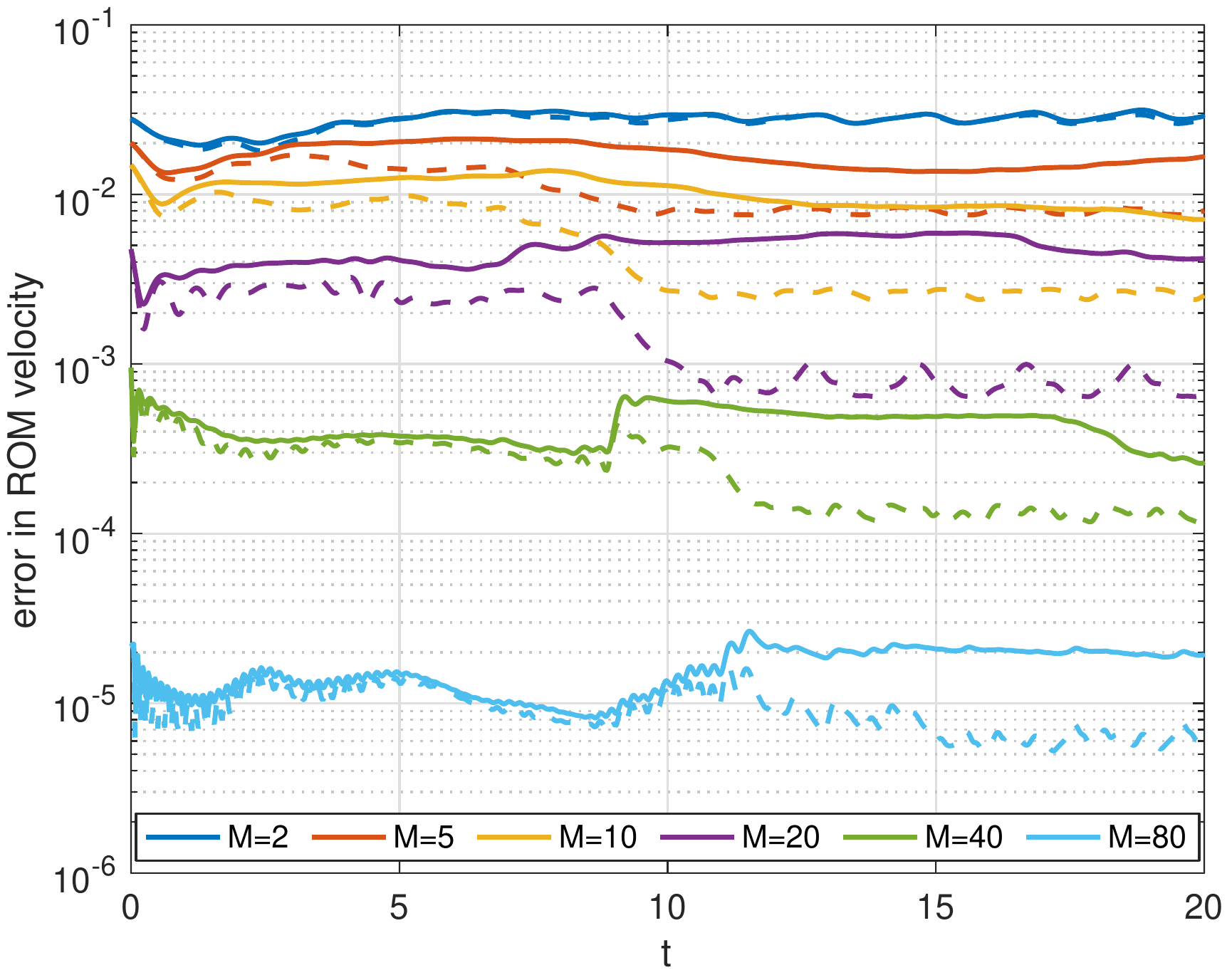}
		\caption{Velocity: $\epsilon_{V}^n$.  \label{fig:error_velocity_actuator}}
	\end{subfigure}
	\hfill	
		\begin{subfigure}[t]{.49\textwidth}
		\centering
		\includegraphics[width=\textwidth]{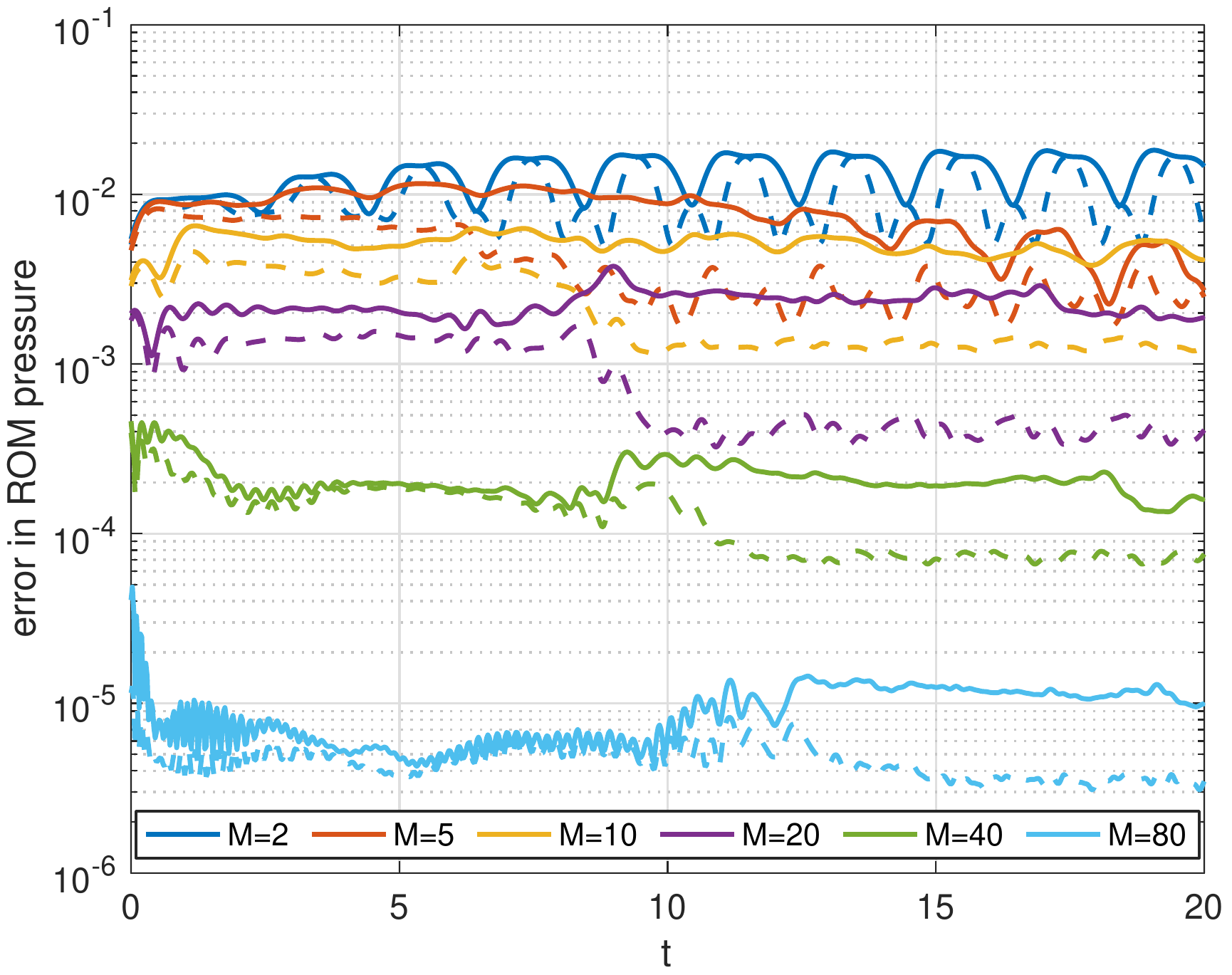}
		\caption{Pressure: $\epsilon_{p}^n$. \label{fig:error_pressure_actuator}}
	\end{subfigure}
\caption{Errors in ROM with respect to FOM for actuator test case. Dashed lines: basis projection error (projecting snapshots onto truncated basis).\label{fig:actuator_errors}}
\end{figure}

\begin{figure}
\centering	
		\includegraphics[width=0.49 \textwidth]{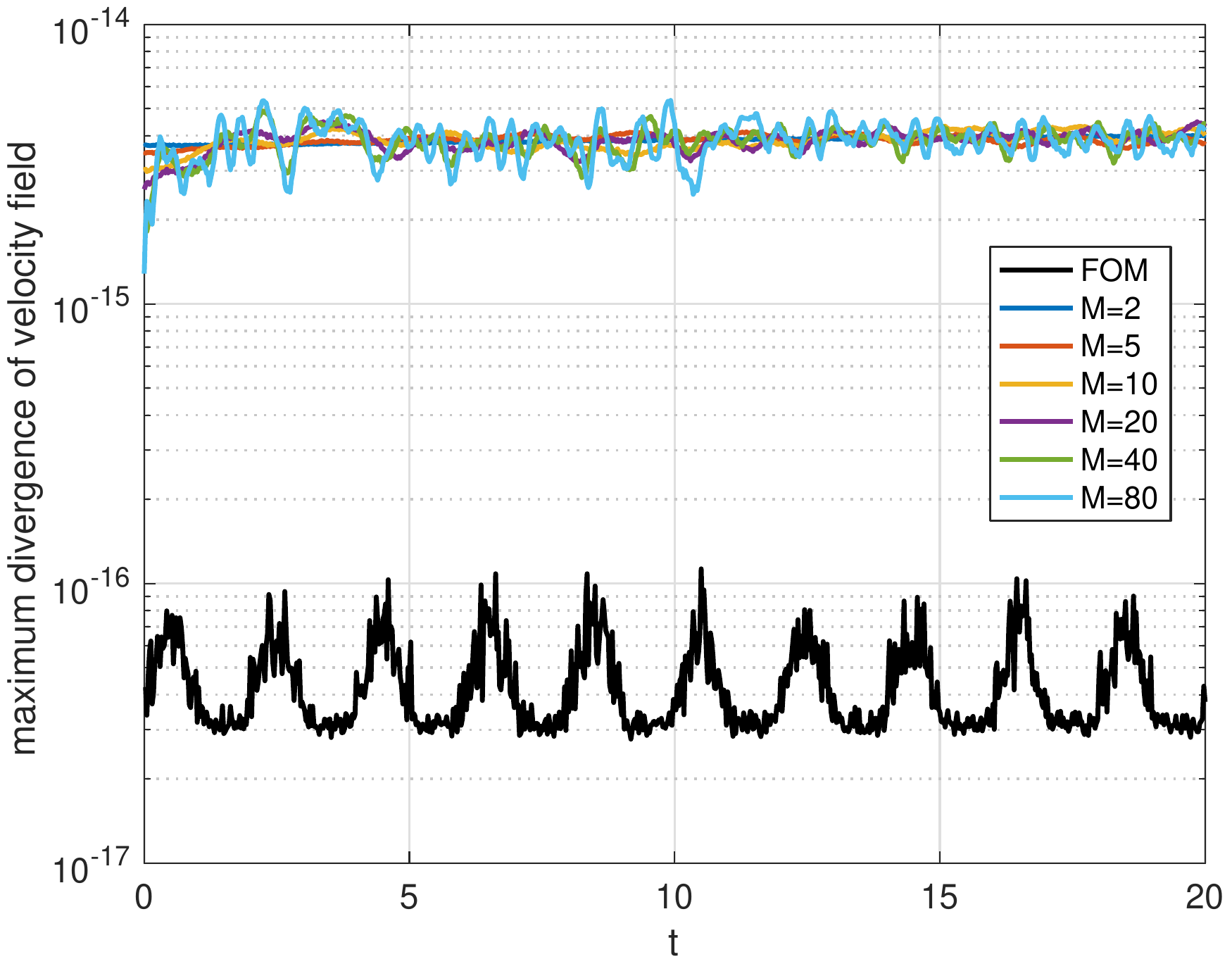}
		\caption{Divergence error. \label{fig:error_divergence_actuator}}
\end{figure}

\hl{A quantitative comparison of the velocity error and the pressure error for different values of $M$ is given in figure \ref{fig:actuator_errors}. The reference velocity used in \eqref{eqn:error_velocity_L2} is the average inflow velocity, $V_{\text{ref}}=1$, whereas the reference pressure used in \eqref{eqn:error_pressure_L2} is the pressure jump over the actuator, $p_{\text{ref}}=\frac{1}{2}C_{T}=\frac{1}{4}$. For all values of $M$, stable results are obtained without requiring stabilization techniques, while clear convergence is achieved upon increasing $M$. For example, for $M=10$, the scaled $L_{2}$ error in the velocity field is around $10^{-2}$. The error in mass conservation of the ROM velocity field, $V_{r} = \Phi a + V_{bc}$, is shown in figure \ref{fig:error_divergence_actuator}. As theoretically derived, our method keeps the velocity field of the ROM divergence-free, also in the case of open flows.}

\begin{figure}[hbtp]
\centering
		\includegraphics[width=0.5\textwidth]{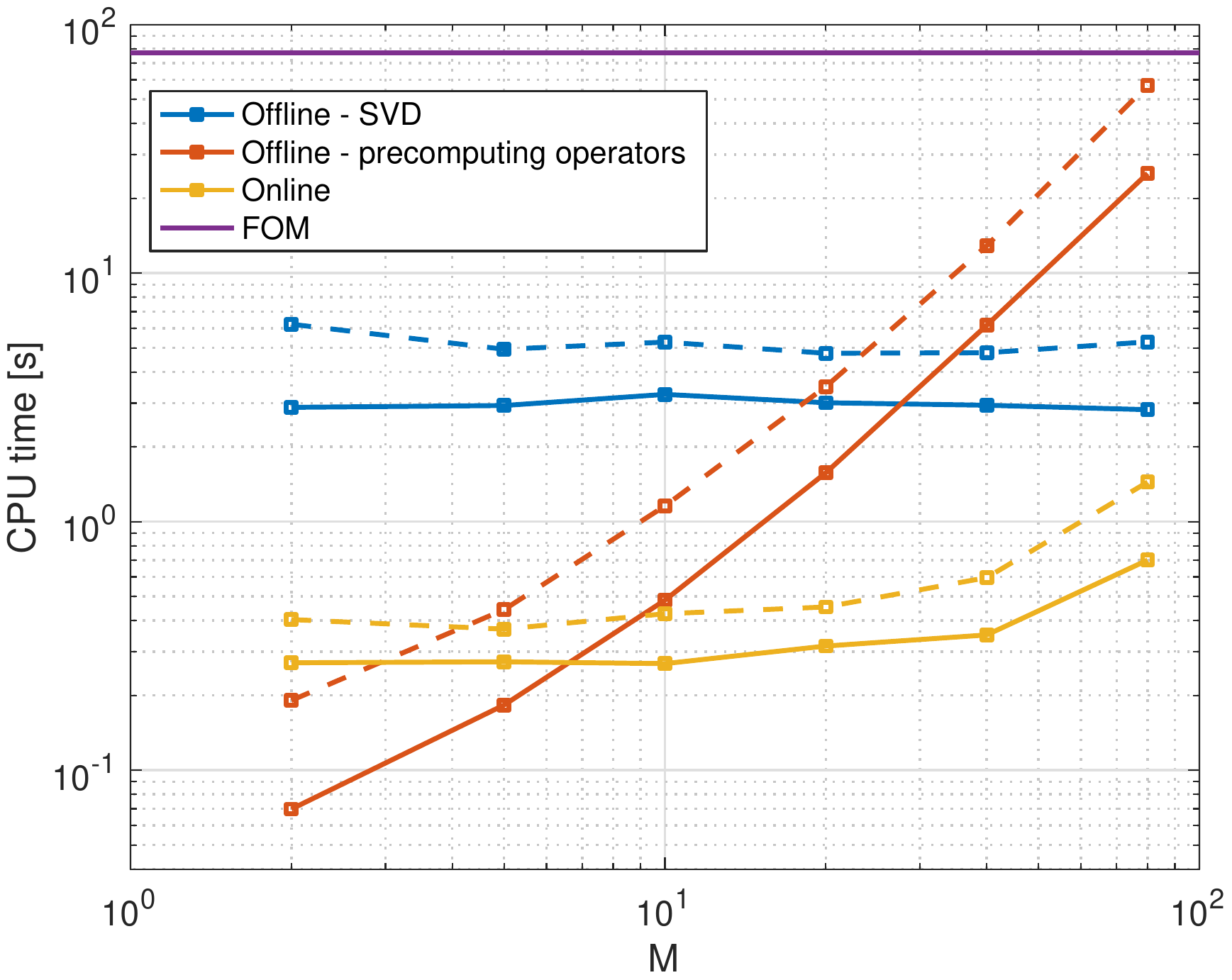}
		\caption{Computational time as a function of number of modes in ROM. Dashed lines: including pressure recovery.\label{fig:CPUtimings_actuator_unsteady}}
\end{figure}

\hl{The computational speed-up of the ROM compared to the FOM is displayed in figure \ref{fig:CPUtimings_actuator_unsteady}. Similar to the other test cases, the $M^{3}$ scaling of the precomputing phase is visible. For $M=10$, the case discussed above, the speed-up factor is around $20$ when including the precomputing phase in the ROM timing, and more than $100$ when not considering the precomputing phase. Like for the lid-driven cavity test case, including pressure recovery roughly doubles the required computational time.}



\FloatBarrier

%% file: conclusions.tex
\section{Conclusions}\label{sec:conclusions}
In this paper we have proposed a novel approach to arrive at an unconditionally stable reduced-order model (ROM) for the incompressible Navier-Stokes equations. The approach hinges on the following four ingredients. First, we have expressed non-linear stability through kinetic energy conservation. Second, we have used a spatially energy-conserving discretisation method as the full order model (FOM), \hl{which has a skew-symmetric convection operator, a divergence-free velocity field, and satisfies a discrete compatibility relation between the divergence and gradient operators}. Third, by performing the projection of the full order model after spatial discretisation an unconditionally stable semi-discrete ROM is obtained. Last, we have used an energy-conserving time integration method. \hl{The resulting ROM is velocity-only} and non-linearly stable. The stability of the method has been \hl{numerically confirmed} for the roll-up of an inviscid shear-layer, for which exact energy conservation was obtained, \hl{independent of the number of modes}.

In addition, we have \hl{used a constrained SVD approach in order to guarantee} momentum conservation on periodic domains. Enforcing momentum conservation comes at the cost of losing a few modes (2 or 3) in the projection matrix, which can be accounted for by taking a few extra modes at a slight increase in computational effort. The constrained SVD approach can be extended to include other linear constraints apart from global momentum.

Furthermore, we have derived a boundary condition treatment for non-homogeneous boundary conditions. The adage of first discretising, then projecting simplifies the boundary condition treatment considerably, as the boundary conditions are built into the discretisation operators and boundary vectors. \hl{The ROM follows from projecting these discrete FOM operators and boundary vectors.} To satisfy the divergence-free constraint, the ROM velocity field is written in terms of a field with homogeneous boundary conditions and a non-homogeneous term, where the latter is obtained by solving once a Poisson equation at the FOM level. In future work, we plan to extend the approach to the case of unsteady and/or parametric boundary conditions. 

\hl{Several test cases, both closed flow and open flow, show the stability, accuracy, and efficiency of the proposed method. Exact energy conservation is shown for an inviscid shear-layer roll-up on a periodic domain, by using implicit time integration. In practice, high-order explicit Runge-Kutta methods have been shown to form a practical alternative, being computationally efficient and suffering only  from a small energy error. The computational speed-up of the ROM compared to the FOM is typically one to two orders of magnitude, depending on the test case and on the number of modes used (and whether the offline phase is included in the computation of the speed-up factor). The pressure is not part of the ROM formulation, but can be recovered once the velocity field is known. This roughly doubles the computational costs, both offline and online (if performed at every time step). With an increasing number of modes, the precomputing phase (in particular assembling the reduced convective operator) becomes the dominant factor in the ROM execution. In our test cases this is not a concern, as the number of modes is typically sufficient before the precomputing phase becomes a dominant factor. \R{ConHyperRed}Nevertheless, one could reduce the complexity of the convective operator (a third order tensor) by using hyper-reduction techniques such as discrete empirical interpolation. In future work we will assess whether hyper-reduction can be performed in such a way that the skew-symmetry of the convective operator is retained. Further efficiency gains can be obtained by making the time step adaptive, e.g.\ by estimating the eigenvalues of the ROM operators.}

As mentioned, this paper has focused mainly on the stability of ROMs, and less on the issue of accuracy in turbulent flows. Our view is that further studies on the accuracy of ROMs (e.g.\ through closure modelling techniques) will benefit significantly from the framework proposed in this paper, \hl{as it is stable and non-dissipative by design}. \hl{Energy conservation is obtained in part because the skew-symmetry property of the convective term is unchanged upon projection. However, this means that the POD-Galerkin method in its basic form cannot be suitable as a reduced model for turbulent flows, as there is no additional energy dissipation coming from the projected convective terms. This is in line with the common understanding that additional measures need to be taken to make POD-Galerkin methods applicable to turbulent flows \cite{Fick2018}.}

%% file: convective_operator.tex
\section{Alternative forms of the convective operator}

\subsection{Continuous}\label{sec:convective_operator}
The convective operator in divergence form can be written such that the role of the advecting velocity $\vt{c}=\vt{u}$ becomes more clear:
\begin{equation}
 C_{\text{div}}(\vt{c},\vt{u}) = \nabla \cdot (\vt{c} \otimes \vt{u}), 
\end{equation}
This distinction allows us to write the divergence form in terms of the advective form $C_{\text{adv}}$ as follows:
\begin{equation}
 C_{\text{div}}(\vt{c},\vt{u}) = C_{\text{adv}} (\vt{c},\vt{u}) +  \vt{u} \left( \nabla \cdot \vt{c}\right). 
\end{equation} 
where 
\begin{equation}
 C_{\text{adv}} (\vt{c},\vt{u}) = (\vt{c} \cdot \nabla) \vt{u}.
\end{equation}
Another commonly used form is the so-called skew-symmetric form, 
\begin{equation}
C_{\text{skew}}(\vt{c},\vt{u}) := \frac{1}{2} C_{\text{div}}(\vt{c},\vt{u}) + \frac{1}{2} C_{\text{adv}} (\vt{c},\vt{u}) = \frac{1}{2} \nabla \cdot (\vt{c} \otimes \vt{u}) + \frac{1}{2} (\vt{c} \cdot \nabla) \vt{u}.
\end{equation}
In case that the advective velocity field is divergence-free ($\nabla \cdot \vt{c}=0$) and the velocity field is sufficiently smooth, the concepts of divergence, advective and skew-symmetric form are equivalent:
\begin{equation}
C(\vt{c},\vt{u}) =  C_{\text{div}}(\vt{c},\vt{u}) =  C_{\text{adv}} (\vt{c},\vt{u}) = C_{\text{skew}}(\vt{c},\vt{u}).
\end{equation}
The notion of skew-symmetry is related to the following property (independent of the divergence-freeness of $\vt{c}$):
\begin{equation}
\begin{split}
C_{\text{skew}} (\vt{c},\vt{u}) \cdot \vt{v} &= \frac{1}{2} (\nabla \cdot (\vt{c} \otimes \vt{u})) \cdot \vt{v} + \frac{1}{2} ((\vt{c} \cdot \nabla) \vt{u}) \cdot \vt{v}, \\
&=\frac{1}{2} \left[ ((\vt{c} \cdot \nabla) \vt{u}) \cdot \vt{v} +  (\vt{u} \cdot \vt{v}) \nabla \cdot \vt{c} \right] +  \frac{1}{2} \left[ \nabla \cdot ( (\vt{u} \cdot \vt{v})  \vt{c}) - (\vt{u} \cdot \vt{v}) \nabla \cdot \vt{c} - ((\vt{c} \cdot \nabla) \vt{v}) \cdot \vt{u} \right], \\
&= \frac{1}{2} ((\vt{c} \cdot \nabla) \vt{u}) \cdot \vt{v} - \frac{1}{2} ((\vt{c} \cdot \nabla) \vt{v}) \cdot \vt{u} + \frac{1}{2} \nabla \cdot ( (\vt{u} \cdot \vt{v})  \vt{c}.
\end{split}
\end{equation}
Upon integration over the entire domain, the contribution of the last term cancels in case of periodic or no-slip boundary conditions, and we obtain
\begin{equation}
(C_{\text{skew}} (\vt{c},\vt{u}), \vt{v}) = \frac{1}{2} ((\vt{c} \cdot \nabla) \vt{u}, \vt{v}) - \frac{1}{2} ((\vt{c} \cdot \nabla) \vt{v}, \vt{u}).
\end{equation}
The convective operator in skew-symmetric form is skew-symmetric `a priori' (i.e.\ without the assumption that $\nabla \cdot \vt{c}=0$),
\begin{equation}
(C_{\text{skew}} (\vt{c},\vt{u}), \vt{v}) = - (\vt{u}, C_{\text{skew}} (\vt{c},\vt{v})).
\end{equation}
The convective operators in advective or divergence form are skew-symmetric provided that $\nabla \cdot \vt{c} = 0$. In that case we have, for example,
\begin{equation}
(C_{\text{div}} (\vt{c},\vt{u}), \vt{v}) = - (\vt{u}, C_{\text{div}} (\vt{c},\vt{v})).
\end{equation}

\subsection{Discrete}\label{sec:convective_operator_discrete}
In two dimensions, the integral of $C^{u}_{\text{div}} = \dd{u^2}{x} +\dd{uv}{y}$ over a finite volume surrounding $u_{i+1/2,j}$ is approximated by
\begin{equation}
C_{h}^{u} (V_{h},u_{h})_{i+1/2,j} :=  \bar{u}_{i+1,j} u_{i+1,j} - \bar{u}_{i,j} u_{i,j} + \bar{v}_{i+1/2,j+1/2} u_{i+1/2,j+1/2} -  \bar{v}_{i+1/2,j-1/2} u_{i+1/2,j-1/2}.
\end{equation}
When interpolating the velocities by mesh-independent weighting of the neighbouring velocities (e.g.\ $u_{i+1,j} = \frac{1}{2} (u_{i+1/2,j} + u_{i+3/2,j})$), leaving the interpolation of the fluxes $\bar{(.)}$ still unspecified, this convective term can be expressed in terms of a matrix-vector product as
\begin{equation}
C_{h}^{u} (V_{h},u_{h}) =\tilde{C}_{h}^{u} (V_{h}) u_{h},
\end{equation}
where (focusing on the $u$-velocities)
\begin{equation}
\tilde{C}_{h}^{u} (V_{h}) = \frac{1}{2}
\begin{pmatrix}
\ddots &  \ddots & \ddots \\
& - \bar{u}_{i,j}  &  \bar{u}_{i+1,j} - \bar{u}_{i,j} &  \bar{u}_{i+1,j} & & \\
 & & - \bar{u}_{i+1,j} &   \bar{u}_{i+2,j} -  \bar{u}_{i+1,j}  &\bar{u}_{i+2,j} & \\
& & & \ddots &  \ddots & \ddots \\
\end{pmatrix}.
\end{equation}
This formulation is possible because the nonlinearity in the convective terms is only quadratic. We note that, apart from the diagonal elements, the matrix is skew-symmetric, independent of the interpolation method for the fluxes $\bar{u}$, $\bar{v}$. Subsequently, the fluxes are computed via mesh-independent weighting of the neighbouring velocities (e.g. $\bar{u}_{i+1,j} = \frac{1}{2} (\bar{u}_{i+1/2,j} + \bar{u}_{i+3/2,j})$). The convective terms can then be rewritten as
\begin{multline}
\frac{1}{2} u_{i+1/2,j} \left[  \frac{1}{2} ( \bar{u}_{i+1/2,j} + \bar{u}_{i+3/2,j} ) -  \frac{1}{2} ( \bar{u}_{i-1/2,j} + \bar{u}_{i+1/2,j} ) +  \frac{1}{2} ( \bar{v}_{i,j+1/2} + \bar{v}_{i+1,j+1/2} ) -  \frac{1}{2} ( \bar{v}_{i,j-1/2} + \bar{v}_{i+1,j-1/2} ) \right] \\
+ \frac{1}{2} u_{i+3/2,j} \frac{1}{2} \left(\bar{u}_{i+1/2,j} + \bar{u}_{i+3/2,j} \right) - \frac{1}{2} u_{i-1/2,j} \frac{1}{2} \left(\bar{u}_{i-1/2,j} + \bar{u}_{i+1/2,j} \right) \\
 + \frac{1}{2} u_{i+1/2,j+1} \frac{1}{2} \left(\bar{v}_{i,j+1/2} + \bar{v}_{i+1,j+1/2} \right) - \frac{1}{2} u_{i+1/2,j-1} \frac{1}{2} \left(\bar{v}_{i,j-1/2} + \bar{v}_{i+1,j-1/2} \right),
\end{multline}
where the term between brackets $\left[ . \right]$ is zero as long as the continuity equation is satisfied.

%% file: BC_FOM.tex
\section{Assembling the ROM operators with boundary conditions}\label{sec:ROM_bc}
\subsection{Summary of boundary conditions of FOM}\label{sec:BC_FOM}
\hl{The semi-discrete equations with non-homogeneous boundary conditions read (equations \eqref{eqn:mass_semidiscrete_bc}-\eqref{eqn:mom_semidiscrete_bc}):
\begin{align}
M_{h} V_{h}(t) &= y_{M}, \label{eqn:mass_conservation_bc}\\
 \Omega_{h}  \frac{\rd V_{h}(t)}{\rd t} & = F^{CD}_{h} (V_{h}(t)) - (G_{h} p_{h}(t) + y_{G}), 
 \end{align}
with
\begin{equation}
F^{CD}_{h} (V_{h}(t)) =- C_{h}(V_{h}(t),V_{h}(t)) + \nu (D_{h} V_{h}(t) + y_{D}) + f_{h}.
\end{equation}
The vector $V_{h}$ does not feature unknowns on the boundaries, but only unknowns associated to interior finite volumes (except in the case of outflow conditions). An example of the positioning of unknowns near a horizontal inflow boundary is given in figure \ref{fig:staggered_boundary}.} 

\begingroup
\begin{figure}[hbtp]
\fontfamily{lmss} 
\fontsize{11pt}{11pt}\selectfont
\centering 
\def\svgwidth{0.6 \textwidth}
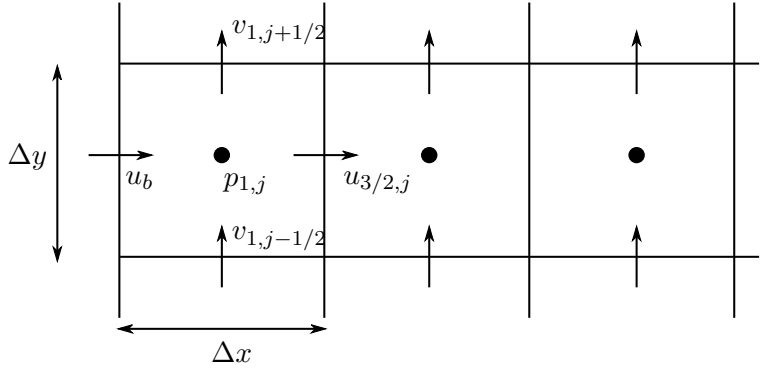 
\caption{\hl{Staggered grid and unknowns near inflow boundary. $u_{b}$ is given and not part of the vector of unknowns.}}
\label{fig:staggered_boundary}
\end{figure}
\endgroup

\hl{\R{PBC6}The divergence-free condition evaluated in a volume adjacent to a boundary then reads:
\begin{equation}\label{eqn:discrete_divergence_bc}
u_{3/2,j} \Delta y  - u_{b} \Delta y + v_{1,j+1/2} \Delta x - v_{1,j-1/2} \Delta x = 0,
\end{equation}
where in terms of equation \eqref{eqn:mass_conservation_bc}, $y_{M}$ contains $u_{b} \Delta y$, and the matrix $M_{h}$ the remaining terms. The momentum equation does not need to be solved for $u_{b}$, as it is known, and is thus solved for the interior points only. For all interior points, the gradient of the pressure is readily available (in this case $(p_{2,j}-p_{1,j})/\Delta x$. Consequently, no pressure boundary conditions need to be specified. It follows directly that for the Poisson equation also no pressure boundary conditions are required: differentiate \eqref{eqn:discrete_divergence_bc} in time, and substitute the momentum equations for $\frac{\rd u_{3/2,j}}{\rd t}$, $\frac{\rd v_{1,j+1/2} }{\rd t}$ and $\frac{\rd v_{1,j-1/2} }{\rd t}$. This gives the discrete Poisson equation without requiring the specification of pressure boundary conditions.}

\hl{The diffusive discretization including boundary conditions reads
\begin{equation}\label{eqn:diffusion_BC}
D_{h} V_{h} + y_{D}.
\end{equation}
On a uniform cartesian grid, $D_{h} \in \mathbb{R}^{N_{V} \times N_{V}}$ contains the standard $\begin{bmatrix} 1 & -2 & 1\end{bmatrix}$ stencil (for each velocity component and each coordinate direction), and $y_{D} \in \mathbb{R}^{N_{V}}$ is a vector containing mostly zeros, except for the finite volumes adjacent to a boundary, which contain the value given by the boundary condition.}

\hl{The convective discretization including boundary conditions is given by
\begin{equation}
C_{h} (V_{h},V_{h}) = K_{h} ( (I_{h} V_{h} + y_{I}) \circ (A_{h} V_{h} + y_{A})),
\end{equation}
where $K_{h} \in \mathbb{R}^{N_V \times N_F}$ is a differencing matrix from faces to volume midpoints (with $-1$ and $1$ as entries), $I_{h} \in \mathbb{R}^{N_F \times N_V}$ interpolates the (convecting) fluxes to the finite volume faces, $A_{h}\in \mathbb{R}^{N_F \times N_V}$ in a similar fashion averages the (convected) velocities to the faces, and $\circ$ is the elementwise (Hadamard) product. $y_{I} \in \mathbb{R}^{N_{F}}$ and $y_{A} \in \mathbb{R}^{N_{F}}$ have nonzero entries for the finite volumes adjacent to boundaries. The simplified form $\tilde{C}(V_{h})V_{h}$ from \ref{sec:convective_operator_discrete}, in which no boundary conditions were present, is obtained as a special case when $y_{I}=y_{A}=0$:
\begin{equation}\label{eqn:convection_BC}
\tilde{C}_{h}(V_{h}) = K_{h} \, \text{diag}(I_{h} V_{h}) \, A_{h},
\end{equation}
where $\text{diag}$ is the operator that transforms a vector to a diagonal matrix with the entries of the vector on its diagonal.}



\subsection{Precomputing the ROM including boundary conditions}\label{sec:ROM_Vbc}
\hl{In this section we explain the projection of the FOM including boundary conditions into a formulation amenable for offline decomposition. We also mention how the ROM operators can be obtained in a rather non-intrusive way from an existing FOM code. There are two sources of boundary conditions when projecting the FOM: (i) the boundary conditions already present in the FOM, explained in the previous section in terms of the vectors $y_{D}$, $y_{I}$, $y_{A}$, and $y_{G}$, acting \textit{only} on finite volumes adjacent to a boundary, and (ii) the boundary field $V_{bc}$ constructed in section \ref{sec:BC}.}

\hl{With both sources present, the ROM formulation of the diffusive term is
\begin{equation}
\Phi^{T} (D_{h} (\Phi a + V_{bc}) + y_{D}) =  (\Phi^{T} D_{h} \Phi) a + \Phi^{T} (D_{h} V_{bc} + y_{D}) =: D_{r} a + y_{D_r}.
\end{equation}
In our FOM code implementation, the diffusive terms (equation \eqref{eqn:diffusion_BC}) are obtained by calling a subroutine \texttt{Diffusion(Vh)}. Given such a subroutine, the ROM terms are simple to obtain, as is best explained using a piece of pseudocode:}
\begin{verbatim}
yDr = Phi'*Diffusion(Vbc)
for i=1:M
   Dr(:,i) = Phi'*(Diffusion(Phi(:,i)) - Diffusion(0))
end
\end{verbatim}

\hl{A similar approach is used to obtain the ROM formulation of the convective terms. These terms read
\begin{equation}
\Phi^{T}  K_{h} ( (I_{h} (\Phi a + V_{bc}) + y_{I}) \circ (A_{h} (\Phi a + V_{bc})+ y_{A})),
\end{equation}
which can be written as
\begin{equation}\label{eqn:convection_ROM}
C_{r,2} (a \otimes a) + C_{r,1} a + y_{C_r}.
\end{equation}
The fact that the ROM convective terms feature not only a purely quadratic term, but also the linear term $C_{r,1}a$, is due to the boundary conditions. This was also observed in e.g.\ \cite{Akhtar2009,Balajewicz2013}.
Similar to the diffusive terms, our FOM code implementation of the convective terms (equation \eqref{eqn:convection_BC}) is given by \texttt{Convection(Vh,Vh)}, where the first argument corresponds to the convecting velocity and the second to the convected velocity, as explained in \ref{sec:convective_operator}. The terms in \eqref{eqn:convection_ROM} are obtained from:}
\begin{verbatim}
yCr = Phi'*Convection(Vbc,Vbc)
for i=1:M
   Cr1(:,i) = Phi'*(Convection(Vbc,Phi(:,i)) + Convection(Phi(:,i),Vbc) 
                    - Convection(Vbc,0) - Convection(0,Vbc)) 
   for j=1:M
      k = (i-1)*M + j
      Cr2(:,k) = Phi'*(Convection(Phi(:,j),Phi(:,i)) - Convection(Phi(:,j),0) 
	                       - Convection(0,Phi(:,i)) + Convection(0,0))
   end
end
\end{verbatim}
\hl{The contribution of pressure boundary conditions (only needed in the case of outflow boundaries) in the ROM simply reads $y_{G_{r}} = \Phi^{T} y_{G}$, and the contribution of body forces similarly reads $f_{r} = \Phi^{T} f_{h}$. The entire ROM formulation then reads
\begin{equation}
\frac{\rd a}{\rd t} = F_{r} (a) = F_{r,2} (a \otimes a) + F_{r,1} a + F_{r,0},
\end{equation}
where
\begin{align}
F_{r,2} &= -C_{r,2}, \\
F_{r,1} &= -C_{r,1} + \nu D_r, \\
F_{r,0} &= -y_{C_r} + \nu y_{D_r} - y_{G_{r}} + f_{r}.
\end{align}}

%% file: figures/staggered_boundary.pdf_tex
\begingroup%
  \makeatletter%
  \providecommand\color[2][]{%
    \errmessage{(Inkscape) Color is used for the text in Inkscape, but the package 'color.sty' is not loaded}%
    \renewcommand\color[2][]{}%
  }%
  \providecommand\transparent[1]{%
    \errmessage{(Inkscape) Transparency is used (non-zero) for the text in Inkscape, but the package 'transparent.sty' is not loaded}%
    \renewcommand\transparent[1]{}%
  }%
  \providecommand\rotatebox[2]{#2}%
  \newcommand*\fsize{\dimexpr\f@size pt\relax}%
  \newcommand*\lineheight[1]{\fontsize{\fsize}{#1\fsize}\selectfont}%
  \ifx\svgwidth\undefined%
    \setlength{\unitlength}{108.51162447bp}%
    \ifx\svgscale\undefined%
      \relax%
    \else%
      \setlength{\unitlength}{\unitlength * \real{\svgscale}}%
    \fi%
  \else%
    \setlength{\unitlength}{\svgwidth}%
  \fi%
  \global\let\svgwidth\undefined%
  \global\let\svgscale\undefined%
  \makeatother%
  \begin{picture}(1,0.48814933)%
    \lineheight{1}%
    \setlength\tabcolsep{0pt}%
    \put(0,0){\includegraphics[width=\unitlength,page=1]{staggered_boundary.pdf}}%
    \put(0.29727432,0.4471651){\color[rgb]{0,0,0}\makebox(0,0)[lt]{\lineheight{1.25}\smash{\begin{tabular}[t]{l}$v_{1,j+1/2}$\end{tabular}}}}%
    \put(0.28650722,0.24224399){\color[rgb]{0,0,0}\makebox(0,0)[lt]{\lineheight{1.25}\smash{\begin{tabular}[t]{l}$p_{1,j}$\end{tabular}}}}%
    \put(0.44584886,0.24507786){\color[rgb]{0,0,0}\makebox(0,0)[lt]{\lineheight{1.25}\smash{\begin{tabular}[t]{l}$u_{3/2,j}$\end{tabular}}}}%
    \put(0,0){\includegraphics[width=\unitlength,page=2]{staggered_boundary.pdf}}%
    \put(0.15555738,0.24466429){\color[rgb]{0,0,0}\makebox(0,0)[lt]{\lineheight{1.25}\smash{\begin{tabular}[t]{l}$u_{b}$\end{tabular}}}}%
    \put(0,0){\includegraphics[width=\unitlength,page=3]{staggered_boundary.pdf}}%
    \put(0.29727432,0.17069702){\color[rgb]{0,0,0}\makebox(0,0)[lt]{\lineheight{1.25}\smash{\begin{tabular}[t]{l}$v_{1,j-1/2}$\end{tabular}}}}%
    \put(0,0){\includegraphics[width=\unitlength,page=4]{staggered_boundary.pdf}}%
    \put(0.26962751,0.00481617){\color[rgb]{0,0,0}\makebox(0,0)[lt]{\lineheight{1.25}\smash{\begin{tabular}[t]{l}$\Delta x$\end{tabular}}}}%
    \put(0,0){\includegraphics[width=\unitlength,page=5]{staggered_boundary.pdf}}%
    \put(-0.00206789,0.27189751){\color[rgb]{0,0,0}\makebox(0,0)[lt]{\lineheight{1.25}\smash{\begin{tabular}[t]{l}$\Delta y$\end{tabular}}}}%
  \end{picture}%
\endgroup%

%% file: adapted_norm.tex
\section{Proof of momentum conservation in weighted norm}\label{sec:proof_adapted_norm}
In this section we prove that the constrained SVD construction in section \ref{sec:momentum_conservation} is conserving momentum for the weighted norm \eqref{eqn:weighted_orthonormality} and keeps $\Phi$ divergence-free. The steps in the construction can be summarized as follows:
\begin{align}
&\text{Form adapted snapshot matrix:}  &\tilde{X} &= X - E E^T \Omega_{h} X, \\
&\text{Transform to include weighted norm:}  &\hat{X} &= \Omega_{h}^{1/2} \tilde{X}, \\
&\text{Perform SVD of $\hat{X}$:}  &\hat{X} &= \hat{\Phi} \Sigma \Psi^{*}, \\
&\text{Transform back to include weighted norm:}  & \tilde{\Phi} &= \Omega_{h}^{-1/2} \hat{\Phi}, \\
&\text{Add $E$ and truncate:} & \Phi &= [E \, \, \tilde{\Phi}]_{M}.
\end{align}
Note that the matrix $E$ is scaled such that 
\begin{equation}
E^{T} \Omega_{h} E = I.
\end{equation}
The proof that $\Phi$ satisfies $\Phi \Phi^{T} \Omega_{h} E = E$ is a substitution exercise:
\begin{equation}
\begin{split}
\Phi \Phi^{T} \Omega_{h} E &= [E \, \, \tilde{\Phi}]_{M} [E \, \, \tilde{\Phi}]_{M}^{T} \Omega_{h} E \\
&= E E^T \Omega_{h} E + \tilde{\Phi}_{M-1} \tilde{\Phi}_{M-1}^{T} \Omega_{h} E \\
&= E +  \Omega_{h}^{-1/2} \hat{\Phi} \hat{\Phi}^{T} \Omega_{h}^{1/2} E \qquad \text{(omitting the truncation subscript)}.
\end{split}
\end{equation}
We proceed to show that the second term equals zero. For this we use equation \eqref{eqn:eigenvectors}: the expression for the modes $\hat{\Phi}$ in terms of the snapshots $\hat{X}$:
\begin{equation}
\lambda_{j} \hat{\Phi}_{j} = \hat{X} \hat{X}^{T} \hat{\Phi}_{j} \quad \rightarrow \quad \lambda_{j} \hat{\Phi}^{T}_{j} \Omega_{h}^{1/2} E = \hat{\Phi}_{j}^{T} \hat{X} \hat{X}^{T} \Omega_{h}^{1/2} E.
\end{equation}
Substituting the expression for the adapted snapshot matrix yields:
\begin{equation}
\begin{split}
\hat{\Phi}^{T}_{j} \Omega_{h}^{1/2} E &= \frac{1}{\lambda_{j}}  \hat{\Phi}_{j} \hat{X} \hat{X}^{T} \Omega_{h}^{1/2} E,\\
&= \frac{1}{\lambda_{j}}  \hat{\Phi}_{j} \Omega_{h}^{1/2} \tilde{X} \tilde{X}^{T} \Omega_{h} E, \\
&= \frac{1}{\lambda_{j}}  \hat{\Phi}_{j} \Omega_{h}^{1/2} (X - E E^T \Omega_{h} X) (X - E E^T \Omega_{h} X)^{T} \Omega_{h} E. \\
\end{split}
\end{equation}
The terms including the snapshot matrix $X$ can be written as
\begin{equation}
\begin{split}
(X - E E^T \Omega_{h} X) (X - E E^T \Omega_{h} X)^{T} \Omega_{h} E =& X X^{T} \Omega_{h} E - E E^T \Omega_{h} X X^{T} \Omega_{h} E \\ 
& \quad - X X^T \Omega_{h} E E^{T} \Omega_{h} E + E E^T \Omega_{h} X X^T \Omega_{h} E E^T \Omega_{h} E,\\
=& X X^{T} \Omega_{h} E - E E^T \Omega_{h} X X^{T} \Omega_{h} E \\ 
& \quad - X X^T \Omega_{h} E + E E^T \Omega_{h} X X^T \Omega_{h} E, \\
=& 0.
\end{split}
\end{equation}
Consequently, as long as $\lambda_{j} \neq 0$, $\hat{\Phi}^{T}_{j} \Omega_{h}^{1/2} E = 0$, and thus $\Phi \Phi^{T} \Omega_{h} E = E$.

It remains to prove that the momentum-conserving construction keeps the basis $\Phi$ divergence-free, in other words, whether
\begin{equation}
\begin{split}
M_{h} \Phi &= M_{h} [E \, \, \tilde{\Phi}]_{M} 
\end{split}
\end{equation}
equals zero. Since $M_{h} E = 0$, we need to only consider $M_{h} \tilde{\Phi}_{j}$ for each column $j$ of $\tilde{\Phi}$.
\begin{equation}
\begin{split}
M_{h} \tilde{\Phi}_{j} &= M_{h}  \Omega_{h}^{-1/2} \hat{\Phi}, \\
&= \frac{1}{\lambda_{j}}  M_{h} \Omega_{h}^{-1/2} \hat{X} \hat{X}^{T} \hat{\Phi}_{j}, \\
&= \frac{1}{\lambda_{j}}  M_{h} \tilde{X} \tilde{X}^{T} \Omega_{h}^{1/2} \hat{\Phi}_{j}, \\
&= \frac{1}{\lambda_{j}}  M_{h} (X - E E^T \Omega_{h} X) (X - E E^T \Omega_{h} X)^{T} \Omega_{h}^{1/2} \hat{\Phi}_{j}, \\
&=0,
\end{split}
\end{equation}
where the last equality follows from $M_{h} E =0$ and $M_{h} X=0$.

%% file: ms.bbl
\begin{thebibliography}{10}

\bibitem{Afkham2018a}
B.~M. Afkham, A.~Bhatt, B.~Haasdonk, and J.~S. Hesthaven.
\newblock {Symplectic Model-Reduction with a Weighted Inner Product}.
\newblock {\em arXiv e-prints arXiv:1803.07799}, 2018.

\bibitem{Afkham2017a}
B.~M. Afkham and J.~S. Hesthaven.
\newblock {Structure Preserving Model Reduction of Parametric Hamiltonian
  Systems}.
\newblock {\em SIAM Journal on Scientific Computing}, 39(6):A2616--A2644, 2017.

\bibitem{Afkham2018c}
B.~M. Afkham and J.~S. Hesthaven.
\newblock {Structure-Preserving Model-Reduction of Dissipative Hamiltonian
  Systems}.
\newblock {\em Journal of Scientific Computing}, 81(1):3--21, 2019.

\bibitem{Akhtar2009}
I.~Akhtar, A.~H. Nayfeh, and C.~J. Ribbens.
\newblock {On the stability and extension of reduced-order Galerkin models in
  incompressible flows: A numerical study of vortex shedding}.
\newblock {\em Theoretical and Computational Fluid Dynamics}, 23(3):213--237,
  2009.

\bibitem{Antoulas2005}
A.~C. Antoulas.
\newblock {\em {Approximation of Large-Scale Dynamical Systems}}.
\newblock Society for Industrial and Applied Mathematics, 2005.

\bibitem{Aubry1993}
N.~Aubry, W.-Y. Lian, and E.~S. Titi.
\newblock {Preserving Symmetries in the Proper Orthogonal Decomposition}.
\newblock {\em SIAM Journal on Scientific Computing}, 14(2):483--505, 1993.

\bibitem{Balajewicz2013}
M.~J. Balajewicz, E.~H. Dowell, and B.~R. Noack.
\newblock {Low-dimensional modelling of high-Reynolds-number shear flows
  incorporating constraints from the Navier–Stokes equation}.
\newblock {\em Journal of Fluid Mechanics}, 729:285--308, 2013.

\bibitem{Ballarin2013}
F.~Ballarin, A.~Manzoni, A.~Quarteroni, and G.~Rozza.
\newblock {Supremizer stabilization of POD-Galerkin approximation of
  parametrized steady incompressible Navier-Stokes equations}.
\newblock {\em International Journal for Numerical Methods in Engineering},
  102(5):1136--1161, 2015.

\bibitem{Benner2015}
P.~Benner, S.~Gugercin, and K.~Willcox.
\newblock {A Survey of Projection-Based Model Reduction Methods for Parametric
  Dynamical Systems}.
\newblock {\em SIAM Review}, 57(4):483--531, 2015.

\bibitem{Berselli2006}
L.~Berselli, T.~Iliescu, and W.~J. Layton.
\newblock {\em {Mathematics of Large Eddy Simulation of Turbulent Flows}}.
\newblock Springer-Verlag, 2006.

\bibitem{Caiazzo2014}
A.~Caiazzo, T.~Iliescu, V.~John, and S.~Schyschlowa.
\newblock {A numerical investigation of velocity-pressure reduced order models
  for incompressible flows}.
\newblock {\em Journal of Computational Physics}, 259:598--616, 2014.

\bibitem{Carlberg2017a}
K.~Carlberg, M.~Barone, and H.~Antil.
\newblock {Galerkin v. least-squares Petrov–Galerkin projection in nonlinear
  model reduction}.
\newblock {\em Journal of Computational Physics}, 330:693--734, 2017.

\bibitem{Carlberg2011}
K.~Carlberg, C.~Bou-Mosleh, and C.~Farhat.
\newblock {Efficient non-linear model reduction via a least-squares
  Petrov-Galerkin projection and compressive tensor approximations}.
\newblock {\em International Journal for Numerical Methods in Engineering},
  86(2):155--181, 2011.

\bibitem{Carlberg2018}
K.~Carlberg, Y.~Choi, and S.~Sargsyan.
\newblock {Conservative model reduction for finite-volume models}.
\newblock {\em Journal of Computational Physics}, 371:280--314, 2018.

\bibitem{Cazemier1998}
W.~Cazemier, R.~Verstappen, and A.~Veldman.
\newblock {Proper orthogonal decomposition and low-dimensional models for
  driven cavity flows}.
\newblock {\em Physics of Fluids}, 10(7):1685--1699, 1998.

\bibitem{Chaturantabut2010}
S.~Chaturantabut and D.~C. Sorensen.
\newblock {Nonlinear Model Reduction via Discrete Empirical Interpolation}.
\newblock {\em SIAM Journal on Scientific Computing}, 32(5):2737--2764, jan
  2010.

\bibitem{Couplet2005}
M.~Couplet, C.~Basdevant, and P.~Sagaut.
\newblock {Calibrated reduced-order POD-Galerkin system for fluid flow
  modelling}.
\newblock {\em Journal of Computational Physics}, 207(1):192--220, 2005.

\bibitem{Fick2018}
L.~Fick, Y.~Maday, A.~T. Patera, and T.~Taddei.
\newblock {A stabilized POD model for turbulent flows over a range of Reynolds
  numbers: Optimal parameter sampling and constrained projection}.
\newblock {\em Journal of Computational Physics}, 371:214--243, 2018.

\bibitem{Fonn2019}
E.~Fonn, H.~van Brummelen, T.~Kvamsdal, and A.~Rasheed.
\newblock {Fast divergence-conforming reduced basis methods for steady
  Navier–Stokes flow}.
\newblock {\em Computer Methods in Applied Mechanics and Engineering},
  346:486--512, 2019.

\bibitem{Gallouet2018}
T.~Gallou{\"{e}}t, R.~Herbin, J.~C. Latch{\'{e}}, and K.~Mallem.
\newblock {Convergence of the Marker-and-Cell Scheme for the Incompressible
  Navier–Stokes Equations on Non-uniform Grids}.
\newblock {\em Foundations of Computational Mathematics}, 18(1):249--289, 2018.

\bibitem{Gresho1998}
P.~M. Gresho and R.~L. Sani.
\newblock {\em {Incompressible Flow {\&} the Finite Element Method -
  Advection-Diffusion {\&} Isothermal Laminar Flow}}.
\newblock John Wiley {\&} Sons, Ltd, 1998.

\bibitem{Gunzburger2007}
M.~D. Gunzburger, J.~S. Peterson, and J.~N. Shadid.
\newblock {Reduced-order modeling of time-dependent PDEs with multiple
  parameters in the boundary data}.
\newblock {\em Computer Methods in Applied Mechanics and Engineering},
  196(4-6):1030--1047, 2007.

\bibitem{Harlow1965}
F.~H. Harlow and J.~E. Welch.
\newblock {Numerical calculation of time-dependent viscous incompressible flow
  of fluid with free surface}.
\newblock {\em Physics of Fluids}, 8(12):2182--2189, 1965.

\bibitem{Hijazi2020}
S.~Hijazi, G.~Stabile, A.~Mola, and G.~Rozza.
\newblock {Data-Driven POD-Galerkin Reduced Order Model for Turbulent Flows}.
\newblock {\em Journal of Computational Physics}, 416:109513, 2019.

\bibitem{Holmes2012}
P.~Holmes, J.~L. Lumley, G.~Berkooz, and C.~W. Rowley.
\newblock {Galerkin projection}.
\newblock {\em Turbulence, Coherent Structures, Dynamical Systems and
  Symmetry}, pages 106--129, 2012.

\bibitem{Kalashnikova2014}
I.~Kalashnikova, M.~F. Barone, S.~Arunajatesan, and B.~G. {van Bloemen
  Waanders}.
\newblock {Construction of energy-stable projection-based reduced order
  models}.
\newblock {\em Applied Mathematics and Computation}, 249:569--596, 2014.

\bibitem{Karasozen2018}
B.~Karas{\"{o}}zen and M.~Uzunca.
\newblock {Energy preserving model order reduction of the nonlinear
  Schr{\"{o}}dinger equation}.
\newblock {\em Advances in Computational Mathematics}, 44(6):1769--1796, 2018.

\bibitem{Kean2019}
K.~Kean and M.~Schneier.
\newblock {Error Analysis of Supremizer Pressure Recovery for POD based Reduced
  Order Models of the time-dependent Navier-Stokes Equations}.
\newblock {\em ArXiv e-print: ArXiv:1909.06022}, sep 2019.

\bibitem{Kramer2019}
B.~Kramer and K.~E. Willcox.
\newblock {Nonlinear model order reduction via lifting transformations and
  proper orthogonal decomposition}.
\newblock {\em AIAA Journal}, 57(6):2297--2307, 2019.

\bibitem{Lassila2014}
T.~Lassila, A.~Manzoni, A.~Quarteroni, and G.~Rozza.
\newblock {Model Order Reduction in Fluid Dynamics: Challenges and
  Perspectives}.
\newblock In {\em Reduced Order Methods for Modeling and Computational
  Reduction}, pages 235--273. Springer International Publishing, 2014.

\bibitem{Mohebujjaman2019}
M.~Mohebujjaman, L.~G. Rebholz, and T.~Iliescu.
\newblock {Physically constrained data-driven correction for reduced-order
  modeling of fluid flows}.
\newblock {\em International Journal for Numerical Methods in Fluids},
  89(3):103--122, 2019.

\bibitem{Mohebujjaman2017}
M.~Mohebujjaman, L.~G. Rebholz, X.~Xie, and T.~Iliescu.
\newblock {Energy balance and mass conservation in reduced order models of
  fluid flows}.
\newblock {\em Journal of Computational Physics}, 346:262--277, 2017.

\bibitem{Noack2016}
B.~R. Noack.
\newblock {From snapshots to modal expansions-bridging low residuals and pure
  frequencies}.
\newblock {\em Journal of Fluid Mechanics}, 802:1--4, 2016.

\bibitem{Noack2005}
B.~R. Noack, P.~Papas, and P.~A. Monkewitz.
\newblock {The need for a pressure-term representation in empirical Galerkin
  models of incompressible shear flows}.
\newblock {\em Journal of Fluid Mechanics}, 523:339--365, 2005.

\bibitem{Peng2016}
L.~Peng and K.~Mohseni.
\newblock {Symplectic Model Reduction of Hamiltonian Systems}.
\newblock {\em SIAM Journal on Scientific Computing}, 38(1):A1--A27, 2016.

\bibitem{Quarteroni2016}
A.~Quarteroni, A.~Manzoni, and F.~Negri.
\newblock {\em {Reduced Basis Methods for Partial Differential Equations}}.
\newblock Springer International Publishing, 2016.

\bibitem{Rempfer2000}
D.~Rempfer.
\newblock {On low-dimensional Galerkin models for fluid flow}.
\newblock {\em Theoretical and Computational Fluid Dynamics}, 14(2):75--88,
  2000.

\bibitem{Rowley2004}
C.~W. Rowley, T.~Colonius, and R.~M. Murray.
\newblock {Model reduction for compressible flows using POD and Galerkin
  projection}.
\newblock {\em Physica D: Nonlinear Phenomena}, 189(1-2):115--129, 2004.

\bibitem{Rubino2019}
S.~Rubino.
\newblock {Numerical analysis of a projection-based stabilized POD-ROM for
  incompressible flows}.
\newblock {\em ArXiv e-print ArXiv:1907.09213}, 2019.

\bibitem{Sanderse2011a}
B.~Sanderse.
\newblock {ECNS: Energy-Conserving Navier-Stokes Solver - Verification of
  steady laminar flows. ECN-E11042}.
\newblock Technical report, Energy research Centre of the Netherlands, 2011.

\bibitem{Sanderse2013b}
B.~Sanderse.
\newblock {Energy-conserving Runge-Kutta methods for the incompressible
  Navier-Stokes equations}.
\newblock {\em Journal of Computational Physics}, 233(1):100--131, 2013.

\bibitem{Sanderse2012a}
B.~Sanderse and B.~Koren.
\newblock {Accuracy analysis of explicit Runge-Kutta methods applied to the
  incompressible Navier-Stokes equations}.
\newblock {\em Journal of Computational Physics}, 231(8):3041--3063, 2012.

\bibitem{Sanderse2011}
B.~Sanderse, S.~Pijl, and B.~Koren.
\newblock {Review of computational fluid dynamics for wind turbine wake
  aerodynamics}.
\newblock {\em Wind Energy}, 14(7):799--819, 2011.

\bibitem{Sanderse2013}
B.~Sanderse, R.~Verstappen, and B.~Koren.
\newblock {Boundary treatment for fourth-order staggered mesh discretizations
  of the incompressible Navier-Stokes equations}.
\newblock {\em Journal of Computational Physics}, 257:1472--1505, 2014.

\bibitem{Shin1997}
D.~Shin and J.~C. Strikwerda.
\newblock {Inf-sup conditions for finite-difference approximations of the
  Stokes equations}.
\newblock {\em Journal of the Australian Mathematical Society Series B-Applied
  Mathematics}, 39(1):121--134, 1997.

\bibitem{Sirovich1987}
L.~Sirovich.
\newblock {Turbulence and the dynamics of coherent structures. I. Coherent
  structures}.
\newblock {\em Quarterly of Applied Mathematics}, 45(3):561--571, 1987.

\bibitem{Smith2014}
R.~Smith.
\newblock {\em {Uncertainty Quantification, Theory, Implementation, and
  Applications}}.
\newblock SIAM, 2014.

\bibitem{Stabile2018}
G.~Stabile and G.~Rozza.
\newblock {Finite volume POD-Galerkin stabilised reduced order methods for the
  parametrised incompressible Navier-Stokes equations}.
\newblock {\em Computers {\&} Fluids}, 173:273--284, 2018.

\bibitem{Taira2017}
K.~Taira, S.~L. Brunton, S.~T. Dawson, C.~W. Rowley, T.~Colonius, B.~J. McKeon,
  O.~T. Schmidt, S.~Gordeyev, V.~Theofilis, and L.~S. Ukeiley.
\newblock {Modal analysis of fluid flows: An overview}.
\newblock {\em AIAA Journal}, 55(12):4013--4041, 2017.

\bibitem{Trias2013a}
F.~Trias, O.~Lehmkuhl, A.~Oliva, C.~P{\'{e}}rez-Segarra, and R.~Verstappen.
\newblock {Symmetry-preserving discretization of Navier–Stokes equations on
  collocated unstructured grids}.
\newblock {\em Journal of Computational Physics}, 258:246--267, 2014.

\bibitem{Trias2011b}
F.~X. Trias and O.~Lehmkuhl.
\newblock {A Self-Adaptive Strategy for the Time Integration of Navier-Stokes
  Equations}.
\newblock {\em Numerical Heat Transfer, Part B: Fundamentals}, 60:116--134,
  2011.

\bibitem{Veldman1990}
A.~E.~P. Veldman.
\newblock {“Missing” boundary conditions? Discretize first, substitute
  next, and combine later}.
\newblock {\em SIAM Journal on Scientific and Statistical Computing},
  11(1):82--91, 1990.

\bibitem{Verstappen2003}
R.~Verstappen and A.~Veldman.
\newblock {Symmetry-preserving discretization of turbulent flow}.
\newblock {\em Journal of Computational Physics}, 187(1):343--368, 2003.

\bibitem{Wang2012a}
Z.~Wang, I.~Akhtar, J.~Borggaard, and T.~Iliescu.
\newblock {Proper orthogonal decomposition closure models for turbulent flows:
  A numerical comparison}.
\newblock {\em Computer Methods in Applied Mechanics and Engineering},
  237-240:10--26, 2012.

\bibitem{Weller2009}
J.~Weller, E.~Lombardi, M.~Bergmann, and A.~Iollo.
\newblock {Numerical methods for low-order modeling of fluid flows based on
  POD}.
\newblock {\em International Journal for Numerical Methods in Fluids},
  63:249--268, 2009.

\bibitem{Xiao2014}
M.~Xiao, P.~Breitkopf, R.~F. Coelho, P.~Villon, and W.~Zhang.
\newblock {Proper orthogonal decomposition with high number of linear
  constraints for aerodynamical shape optimization}.
\newblock {\em Applied Mathematics and Computation}, 247:1096--1112, 2014.

\bibitem{Xiao2013}
M.~Xiao, P.~Breitkopf, R.~{Filomeno Coelho}, C.~Knopf-Lenoir, P.~Villon, and
  W.~Zhang.
\newblock {Constrained Proper Orthogonal Decomposition based on
  QR-factorization for aerodynamical shape optimization}.
\newblock {\em Applied Mathematics and Computation}, 223:254--263, 2013.

\end{thebibliography}
